\documentstyle[amssymb,amsmath]{article}

\newtheorem{th}{Theorem}[section]
\newtheorem{lem}[th]{Lemma}

\newtheorem{cor}[th]{Corollary}
\newtheorem{defn}[th]{Definition}
\newenvironment{defn-new}{\begin{defn} \em}{\end{defn}}
\newtheorem{rem}[th]{Remark}
\newenvironment{rem-new}{\begin{rem} \em}{\end{rem}}
\newtheorem{ex}[th]{Example}
\newenvironment{ex-new}{\begin{ex} \em}{\end{ex}}

\newtheorem{notation}[th]{Notation}
\newenvironment{notation-new}{\begin{rem} \em}{\end{rem}}

\newenvironment{agr-new}{\begin{rem} \em}{\end{rem}}

\makeatletter \@addtoreset{equation}{section} \makeatother

\makeatletter \@addtoreset{figure}{section} \makeatother

\begin{document}

\begin{center}
{\LARGE {\bf On $(N(k),\xi )$-semi-Riemannian manifolds: Pseudosymmetries}}\\%
[0pt]
\bigskip \bigskip {\large {\bf Mukut Mani Tripathi and Punam Gupta}}
\end{center}

{\bf Abstract.} Definition of $({\cal T}_{\!a},{\cal T}_{\!b})$%
-pseudosymmetric semi-Riemannian manifold is given. $({\cal T}_{\!a},{\cal T}%
_{\!b})$-pseudosy-\newline
mmetric $(N(k),\xi )$-semi-Riemannian manifolds are classified. Some results
for ${\cal T}_{\!a}$-pseudosymmetric $(N(k),\xi )$-semi-Riemannian manifolds
are obtained. $({\cal T}_{\!a},{\cal T}_{\!b},S^{\ell })$-pseudosymmetric
semi-Riemannian manifolds are defined. $({\cal T}_{\!a},{\cal T}%
_{\!b},S^{\ell })$-pseudosymmetric $(N(k),\xi )$-semi-Riemannian manifolds
are classified. Some results for $(R,{\cal T}_{\!a},S^{\ell })$%
-pseudosymmetric $(N(k),\xi )$-semi-Riemannian manifolds are obtained. In
particular, some results for $(R,{\cal T}_{\!a},S)$-pseudosymmetric $%
(N(k),\xi )$-semi-Riemannian manifolds are also obtained. After that, the
definition of $({\cal T}_{\!a},S_{{\cal T}_{b}})$-pseudosymmetric
semi-Riemannian manifold is given. $({\cal T}_{\!a},S_{{\cal T}_{b}})$%
-pseudosymmetric $(N(k),\xi )$-semi-Riemannian manifolds are classified. It
is proved that a $(R,S_{{\cal T}_{\!a}})$-pseudosymmetric $(N(k),\xi )$%
-semi-Riemannian manifold is either Einstein or $L=k$ under an algebraic
condition. Some results for $({\cal T}_{\!a},S)$-pseudosymmetric $(N(k),\xi
) $-semi-Riemannian manifolds are also obtained. In last, $({\cal T}%
_{\!a},S_{{\cal T}_{\!b}},S^{\ell })$-pseudosymmetric semi-Riemannian
manifolds are defined and $({\cal T}_{\!a},S_{{\cal T}_{\!b}},S^{\ell })$%
-pseudosymmetric $(N(k),\xi )$-semi-Riemannian manifolds are classified.
\medskip

\medskip

\noindent {\bf Keywords.} ${\cal T}$-curvature tensor; quasi-conformal
curvature tensor; conformal curvature tensor; conharmonic curvature tensor;
concircular curvature tensor; pseudo-projective curvature tensor; projective
curvature tensor; ${\cal M}$-projective curvature tensor; ${\cal W}_{i}$%
-curvature tensors $(i=0,\ldots ,9)$; ${\cal W}_{j}^{\ast }$-curvature
tensors $(j=0,1)$; $\eta $-Einstein manifold; Einstein manifold; $N(k)$%
-contact metric manifold; $\left( \varepsilon \right) $-Sasakian manifold;
Sasakian manifold; Kenmotsu manifold; $\left( \varepsilon \right) $%
-para-Sasakian manifold; para-Sasakian manifold; $(N(k),\xi )$%
-semi-Riemannian manifolds; $\left( {\cal T}_{\!a},{\cal T}_{\!b}\right) $%
-pseudosymmetric semi-Riemannian manifold; $\left( {\cal T}_{\!a},{\cal T}%
_{\!b},S^{\ell }\right) $-pseudosymm\allowbreak etric semi-Riemannian
manifold; $({\cal T}_{\!a},S_{{\cal T}_{b}}) $-pseudosymmetric
semi-Riemannian manifold and $({\cal T}_{\!a},S_{{\cal T}_{b}},S^{\ell })$%
-pseudosymmetric semi-Riemannian manifold. \medskip

\section{Introduction}

Let $M$ be an $n$-dimensional differentiable manifold and ${\frak X}(M)$ the
Lie algebra of vector fields in $M$. We assume that $X,X_{1},\ldots
,X_{s},Y,Z,U,V,W\in {\frak X}(M)$. It is well known that every $(1,1)$
tensor field ${\cal A}$ on a differentiable manifold determines a derivation
${\cal A}\cdot $ of the tensor algebra on the manifold, commuting with
contractions. For example, a $(1,1)$ tensor field ${\cal B}(V,U)$ induces
the derivation ${\cal B}(V,U)\cdot $, thus given a $(0,s)$ tensor field $%
{\cal K}$, we can associate a $(0,s+2)$ tensor ${\cal B}\cdot {\cal K}$,
defined by%
\begin{eqnarray}
({\cal B}\cdot {\cal K})(X_{1},\ldots ,X_{s},X,Y) &=&({\cal B}(X,Y)\cdot
{\cal K})(X_{1},\ldots ,X_{s})  \nonumber \\
&=&-\,{\cal K}({\cal B}(X,Y)X_{1},\ldots ,X_{s})-\cdots  \nonumber \\
&&-\,{\cal K}(X_{1},\ldots ,{\cal B}(X,Y)X_{s}).  \label{eq-RR}
\end{eqnarray}%
Next, for a tensor $\sigma $ of type $(0,2)$, the operators $(X\wedge
_{\sigma }Y)$ and $Q(\sigma ,{\cal K})$are defined by
\[
(X\wedge _{\sigma }Y)Z=\sigma (Y,Z)X-\sigma (X,Z)Y,
\]%
\begin{eqnarray}
Q(\sigma ,{\cal K})(X_{1},\ldots ,X_{s},X,Y) &=&((X\wedge _{\sigma }Y)\cdot
{\cal K)}(X_{1},\ldots ,X_{s})  \nonumber \\
&=&-\,{\cal K}((X\wedge _{\sigma }Y)X_{1},\ldots ,X_{s})-\cdots  \nonumber \\
&&-\,{\cal K}(X_{1},\ldots ,(X\wedge _{\sigma }Y)X_{s}).  \label{eq-pseudo}
\end{eqnarray}

Let $(M,g)$ be an $n$-dimensional semi-Riemannian manifold. Then $(M,g)$ is
said to be

\begin{enumerate}
\item[{\bf (a)}] {\em pseudosymmetric} \cite{Deszcz-Grycak-87} if its
curvature tensor $R$ satisfies
\[
R\cdot R=L_{g}Q(g,R),
\]%
where $L_{g}$ is some smooth function on $M$,

\item[{\bf (b)}] {\em Ricci-generalized pseudosymmetric} \cite{Deszcz-90} if
it satisfies
\[
R\cdot R=L_{S}Q(S,R),
\]%
where $L_{S}$ is some smooth function on $M$ and $S$ is the the Ricci tensor,

\item[{\bf (c)}] {\em Ricci-pseudosymmetric} \cite{Deszcz-89} if on the set $%
{\cal U}=\left\{ x\in M:\left( S-r/n\right) _{x}\not=0\right\} $,
\[
R\cdot S=LQ(g,S),
\]%
where $L$ is some function on ${\cal U}$.
\end{enumerate}

A pseudosymmetric manifold is a generalization of manifold of constant
curvature, symmetric manifold $\left( \nabla R=0\right) $ and semisymmetric
manifold $\left( R\cdot R=0\right) $. Deszcz et al. \cite{DVY} proved that
hypersurfaces in spaces of constant curvature, with exactly two distinct
principal curvatures at every point, are pseudosymmetric.
Ricci-pseudosymmetric manifold is a generalization of manifold of constant
curvature, Einstein manifold, Ricci symmetric manifold $\left( \nabla
S=0\right) $, symmetric manifold, semisymmetric manifold, pseudosymmetric
manifold and Ricci-semisymmetric manifold$\left( R\cdot S=0\right) $.
Similar to pseudosymmetry condition, Deszcz and Grycak \cite%
{Deszcz-90.,Deszcz-91,Deszcz-Grycak-89} and \"{O}zg\"{u}r \cite{Ozgur-2005}
also studied Weyl pseudosymmetric manifolds $\left( R\cdot {\cal C}=L_{g}Q(g,%
{\cal C})\right) $. In 1990, Prvanovi\'{c} \cite{Prvanovic-90} studied the
condition
\[
R\cdot \tilde{T}=LQ(S^{\ell },\tilde{T}),\qquad \ell =0,1,2,\ldots ,
\]%
where $\tilde{T}$ is some $\left( 0,4\right) $-tensor field and $L$ is some
smooth function on $M$. If $\tilde{T}=R$ and $\ell =0$, this condition
becomes the condition for pseudosymmetry and if $\tilde{T}=R$ and $\ell =1$,
this condition becomes the condition for Ricci-generalized pseudosymmetry.

Apart from curvature tensor and the Weyl conformal curvature tensor,
quasi-conformal curvature tensor, concircular curavture tensor, conharmonic
curvature tensor, pseudo-projective curvature tensor, projective curvature
tensor are important curvature tensors in the semi-Riemannian point of view.
It is interesting to study these curvature tensors and other curvature
tensors on $(N(k),\xi )$-semi-Riemannian manifolds. In this paper, we study
several derivation conditions on $(N(k),\xi )$-semi-Riemannian manifolds.
The paper is organized as follows. In Section \ref{sect-GCT}, we give the
definition of ${\cal T}$-curvature tensor. In Section~\ref{sect-NK}, we give
examples and properties of $(N(k),\xi )$-semi-Riemannian manifolds. In
Section~\ref{sect-TTP}, $({\cal T}_{\!a},{\cal T}_{\!b})$-pseudosymmetric
semi-Riemannian manifolds are defined and studied. Some results for ${\cal T}%
_{\!a}$-pseudosymmetric $(N(k),\xi )$-semi-Riemannian manifolds are given.
In Section~\ref{sect-TTSP}, $({\cal T}_{\!a},{\cal T}_{\!b},S^{\ell })$%
-pseudosymmetric semi-Riemannian manifolds are defined and studied. Some
results for $(R,{\cal T}_{\!a},S^{\ell })$-pseudosymmetric $(N(k),\xi )$%
-semi-Riemannian manifolds are given. In particular, some results for $(R,%
{\cal T}_{\!a},S)$-pseudosymmetric $(N(k),\xi )$-semi-Riemannian manifolds
are obtained. In Section~\ref{sect-TSP}, the definition of $({\cal T}%
_{\!a},S_{{\cal T}_{b}})$-pseudosymmetric semi-Riemannian manifold is given.
$({\cal T}_{\!a},S_{{\cal T}_{b}})$-pseudosymmetric $(N(k),\xi )$%
-semi-Riemannian manifolds are classified. It is proved that a $(R,S_{{\cal T%
}_{\!a}})$-pseudosymmetric $(N(k),\xi )$-semi-Riemannian manifold is either
Einstein or $L=k$ under an algebraic condition. Some results for $({\cal T}%
_{\!a},S)$-pseudosymmetric $(N(k),\xi )$-semi-Riemannian manifolds are
obtained. In the last section, the definition of $(T_{\!a},S_{{\cal T}%
_{\!b}},S^{\ell })$-pseudosymmetric semi-Riemannian manifold is given. $%
(T_{\!a},S_{{\cal T}_{\!b}},S^{\ell })$-pseudosymmetric $(N(k),\xi )$%
-semi-Riemannian manifolds are classified.

\section{Preliminaries\label{sect-GCT}}

\begin{defn}
\label{defn-GCT} In an $n$-dimensional semi-Riemannian manifold $\left(
M,g\right) $, {\em ${\cal T}$-curvature tensor} {\rm {\cite{Tripathi-Gupta}}
}is a tensor of type $(1,3)$, which is defined by
\begin{eqnarray}
{\cal T}\left( X,Y\right) Z &=&a_{0}\,R\left( X,Y\right) Z  \nonumber \\
&&+\ a_{1}\,S\left( Y,Z\right) X+a_{2}\,S\left( X,Z\right) Y+a_{3}\,S(X,Y)Z
\nonumber \\
&&+\ a_{4}\,g\left( Y,Z\right) QX+a_{5}\,g\left( X,Z\right)
QY+a_{6}\,g(X,Y)QZ  \nonumber \\
&&+\ a_{7}\,r\left( g\left( Y,Z\right) X-g\left( X,Z\right) Y\right) ,
\label{eq-GCT}
\end{eqnarray}%
where $a_{0},\ldots ,a_{7}$ are real numbers; and $R$, $S$, $Q$ and $r$ are
the curvature tensor, the Ricci tensor, the Ricci operator and the scalar
curvature respectively.
\end{defn}

In particular, the ${\cal T}$-curvature tensor is reduced to

\begin{enumerate}
\item the{\em \ curvature tensor} $R$ if
\[
a_{0}=1,\quad a_{1}=\cdots =a_{7}=0,
\]

\item the {\em quasi-conformal curvature tensor} ${\cal C}_{\ast }$ \cite%
{Yano-Sawaki-68} if
\[
a_{1}=-\,a_{2}=a_{4}=-\,a_{5},\quad a_{3}=a_{6}=0,\quad a_{7}=-\,\frac{1}{n}%
\left( \frac{a_{0}}{n-1}+2a_{1}\right) ,
\]

\item the {\em conformal curvature tensor} ${\cal C}$ \cite[p.~90]%
{Eisenhart-49} if
\[
\qquad a_{0}=1,\; a_{1}=-\,a_{2}=a_{4}=-\,a_{5}=-\,\frac{1}{n-2},\;
a_{3}=a_{6}=0,\; a_{7}=\frac{1}{(n-1)(n-2)},
\]

\item the {\em conharmonic curvature tensor} ${\cal L}$ \cite{Ishii-57} if
\[
a_{0}=1,\quad a_{1}=-\,a_{2}=a_{4}=-\,a_{5}=-\,\frac{1}{n-2},\,\quad
a_{3}=a_{6}=0,\quad a_{7}=0,
\]

\item the {\em concircular curvature tensor} ${\cal V}$ (\cite{Yano-40},
\cite[p. 87]{Yano-Bochner-53}) if
\[
a_{0}=1,\quad a_{1}=a_{2}=a_{3}=a_{4}=a_{5}=a_{6}=0,\quad a_{7}=-\,\frac{1}{%
n(n-1)},
\]

\item the {\em pseudo-projective curvature tensor }${\cal P}_{\ast }$ \cite%
{Prasad-2002} if
\[
a_{1}=-\,a_{2},\quad a_{3}=a_{4}=a_{5}=a_{6}=0,\quad a_{7}=-\,\frac{1}{n}%
\left( \frac{a_{0}}{n-1}+a_{1}\right) ,
\]

\item the {\em projective curvature tensor} ${\cal P}$ \cite[p. 84]%
{Yano-Bochner-53} if
\[
a_{0}=1,\quad a_{1}=-\,a_{2}=-\,\frac{1}{(n-1)}\text{,\quad }%
a_{3}=a_{4}=a_{5}=a_{6}=a_{7}=0,
\]

\item the ${\cal M}${\em -projective curvature tensor }\cite%
{Pokhariyal-Mishra-71} if
\[
a_{0}=1,\quad a_{1}=-\,a_{2}=a_{4}=-\,a_{5}=-\frac{1}{2(n-1)},\quad
a_{3}=a_{6}=a_{7}=0,
\]%
\vspace{-0.5cm}

\item the ${\cal W}_{0}${\em -curvature tensor} \cite[Eq. (1.4)]%
{Pokhariyal-Mishra-71} if
\[
a_{0}=1,\quad a_{1}=-\,a_{5}=-\,\frac{1}{(n-1)},\quad
a_{2}=a_{3}=a_{4}=a_{6}=a_{7}=0,
\]%
\vspace{-0.5cm}

\item the ${\cal W}_{0}^{\ast }${\em -curvature tensor} \cite[Eq. (2.1)]%
{Pokhariyal-Mishra-71} if
\[
a_{0}=1,\quad a_{1}=-\,a_{5}=\frac{1}{(n-1)},\quad
a_{2}=a_{3}=a_{4}=a_{6}=a_{7}=0,
\]

\item the ${\cal W}_{1}${\em -curvature tensor} \cite{Pokhariyal-Mishra-71}
if
\[
a_{0}=1,\quad a_{1}=-\,a_{2}=\frac{1}{(n-1)},\quad
a_{3}=a_{4}=a_{5}=a_{6}=a_{7}=0,
\]

\item the ${\cal W}_{1}^{\ast }${\em -curvature tensor} \cite%
{Pokhariyal-Mishra-71} if
\[
a_{0}=1,\quad a_{1}=-\,a_{2}=-\,\frac{1}{(n-1)},\quad
a_{3}=a_{4}=a_{5}=a_{6}=a_{7}=0,
\]

\item the ${\cal W}_{2}${\em -curvature tensor} \cite{Pokhariyal-Mishra-70}
if
\[
a_{0}=1,\quad a_{4}=-\,a_{5}=-\,\frac{1}{(n-1)},\quad
a_{1}=a_{2}=a_{3}=a_{6}=a_{7}=0,
\]

\item the ${\cal W}_{3}${\em -curvature tensor} \cite{Pokhariyal-Mishra-71}
if
\[
a_{0}=1,\quad a_{2}=-\,a_{4}=-\,\frac{1}{(n-1)},\quad
a_{1}=a_{3}=a_{5}=a_{6}=a_{7}=0,
\]

\item the ${\cal W}_{4}${\em -curvature tensor} \cite{Pokhariyal-Mishra-71}
if
\[
a_{0}=1,\quad a_{5}=-\,a_{6}=\frac{1}{(n-1)},\quad
a_{1}=a_{2}=a_{3}=a_{4}=a_{7}=0,
\]

\item the ${\cal W}_{5}${\em -curvature tensor} \cite{Pokhariyal-82} if
\[
a_{0}=1,\quad a_{2}=-\,a_{5}=-\,\frac{1}{(n-1)},\quad
a_{1}=a_{3}=a_{4}=a_{6}=a_{7}=0,
\]

\item the ${\cal W}_{6}${\em -curvature tensor} \cite{Pokhariyal-82} if
\[
a_{0}=1,\quad a_{1}=-\,a_{6}=-\,\frac{1}{(n-1)},\quad
a_{2}=a_{3}=a_{4}=a_{5}=a_{7}=0,
\]

\item the ${\cal W}_{7}${\em -curvature tensor} \cite{Pokhariyal-82} if
\[
a_{0}=1,\quad a_{1}=-\,a_{4}=-\,\frac{1}{(n-1)},\quad
a_{2}=a_{3}=a_{5}=a_{6}=a_{7}=0,
\]

\item the ${\cal W}_{8}${\em -curvature tensor} \cite{Pokhariyal-82} if
\[
a_{0}=1,\quad a_{1}=-\,a_{3}=-\,\frac{1}{(n-1)},\quad
a_{2}=a_{4}=a_{5}=a_{6}=a_{7}=0,
\]

\item the ${\cal W}_{9}${\em -curvature tensor} \cite{Pokhariyal-82} if
\[
a_{0}=1,\quad a_{3}=-\,a_{4}=\frac{1}{(n-1)},\quad
a_{1}=a_{2}=a_{5}=a_{6}=a_{7}=0.
\]
\end{enumerate}

Denoting
\[
{\cal T}\left( X,Y,Z,V\right) =g({\cal T}\left( X,Y\right) Z,V),
\]%
we write the curvature tensor ${\cal T}$ in its $\left( 0,4\right) $ form as
follows.
\begin{eqnarray}
{\cal T}\left( X,Y,Z,V\right) &=&a_{0}\,R\left( X,Y,Z,V\right)  \nonumber \\
&&+\ a_{1}\,S\left( Y,Z\right) g\left( X,V\right) +a_{2}\,S\left( X,Z\right)
g\left( Y,V\right)  \nonumber \\
&&+\ a_{3}\,S\left( X,Y\right) g\left( Z,V\right) +a_{4}\,S\left( X,V\right)
g\left( Y,Z\right)  \nonumber \\
&&+\ a_{5}\,S\left( Y,V\right) g\left( X,Z\right) +a_{6}\,S\left( Z,V\right)
g\left( X,Y\right)  \nonumber \\
&&+\ a_{7}\,r\left( g\left( Y,Z\right) g\left( X,V\right) -g\left(
X,Z\right) g\left( Y,V\right) \right) .  \label{eq-gen-cur-1}
\end{eqnarray}

\begin{notation}
We will call ${\cal T}$-curvature tensor as ${\cal T}_{a}$-curvature tensor,
whenever it is necessary. If $a_{0},\ldots ,a_{7}$ are replaced by $%
b_{0},\ldots ,b_{7}$ in the definition of ${\cal T}$-curvature tensor, then
we will call ${\cal T}$-curvature tensor as ${\cal T}_{b}$-curvature tensor.
\end{notation}

\section{$(N(k),\protect\xi )$-semi-Riemannian manifolds\label{sect-NK}}

Let $(M,g)$ be an $n$-dimensional semi-Riemannian manifold \cite{ONeill-83}
equipped with a semi-Riemannian metric $g$. If ${\rm index}(g)=1$ then $g$
is a Lorentzian metric and $(M,g)$ a Lorentzian manifold \cite%
{Beem-Ehrlich-81}. If $g$ is positive definite then $g$ is an usual
Riemannian metric and $(M,g)$ a Riemannian manifold. \medskip

The $k${\em -nullity distribution} \cite{Tanno-88} of $(M,g)$ for a real
number $k$ is the distribution%
\[
N(k):p\mapsto N_{p}(k)=\left\{ Z\in
T_{p}M:R(X,Y)Z=k(g(Y,Z)X-g(X,Z)Y)\right\} .
\]%
Let $\xi $ be a non-null unit vector field in $(M,g)$ and $\eta $ its
associated $1$-form. Thus
\[
g(\xi ,\xi )=\varepsilon ,
\]%
where $\varepsilon =1$ or $-\,1$ according as $\xi $ is spacelike or
timelike, and
\begin{equation}
\eta \left( X\right) =\varepsilon g\left( X,\xi \right) ,\qquad \eta \left(
\xi \right) =1.  \label{eq-cond}
\end{equation}

\begin{defn}
An $(N(k),\xi )${\em -semi-Riemannian manifold} {\rm \cite{TG}} consists of
a semi-Riemannian manifold $(M,g)$, a $k$-nullity distribution $N(k)$ on $%
(M,g) $ and a non-null unit vector field $\xi $ in $(M,g)$ belonging to $%
N(k) $.
\end{defn}

Now, we intend to give some examples of $(N(k),\xi )$-semi-Riemannian
manifolds. For this purpose we collect some definitions from the geometry of
almost contact manifolds and almost paracontact manifolds as follows:

\subsection*{Almost contact manifolds}

Let $M$ be a smooth manifold of dimension $n=2m+1$. Let $\varphi $, $\xi $
and $\eta $ be tensor fields of type $(1,1)$, $(1,0)$ and $(0,1)$,
respectively. If $\varphi $, $\xi $ and $\eta $ satisfy the conditions
\begin{equation}
\varphi ^{2}=-I+\eta \otimes \xi ,  \label{eq-str-1}
\end{equation}%
\begin{equation}
\eta (\xi )=1,  \label{eq-str-2}
\end{equation}%
where $I$ denotes the identity transformation, then $M$ is said to have an
almost contact structure $\left( \varphi ,\xi ,\eta \right) $. A manifold $M$
alongwith an almost contact structure is called an {\em almost contact
manifold} \cite{Blair-76}. Let $g$ be a semi-Riemannian metric on $M$ such
that
\begin{equation}
g\left( \varphi X,\varphi Y\right) =g\left( X,Y\right) -\varepsilon \eta
(X)\eta \left( Y\right) ,  \label{eq-str-3}
\end{equation}%
where $\varepsilon =\pm 1$. Then $(M,g)$ is an $\left( \varepsilon \right) $-%
{\em almost contact metric manifold} \cite{Duggal-90-IJMMS} equipped with an
$\left( \varepsilon \right) ${\em -almost contact metric structure} $%
(\varphi ,\xi ,\eta ,g,\varepsilon )$. In particular, if the metric $g$ is
positive definite, then an $(\varepsilon )$-almost contact metric manifold
is the usual {\em almost contact metric manifold }\cite{Blair-76}. \medskip

From (\ref{eq-str-3}), it follows that
\begin{equation}
g\left( X,\varphi Y\right) =-g\left( \varphi X,Y\right)  \label{eq-str-5}
\end{equation}%
and
\begin{equation}
g\left( X,\xi \right) =\varepsilon \eta (X).  \label{eq-str-6}
\end{equation}%
From (\ref{eq-str-2}) and (\ref{eq-str-6}), we have
\begin{equation}
g\left( \xi ,\xi \right) =\varepsilon .  \label{eq-str-7}
\end{equation}

In an $\left( \varepsilon \right) $-almost contact metric manifold, the
fundamental $2$-form $\Phi $ is defined by
\begin{equation}
\Phi (X,Y)=g(X,\varphi Y).  \label{eq-str-4}
\end{equation}%
An $\left( \varepsilon \right) $-almost contact metric manifold with $\Phi
=d\eta $ is an $\left( \varepsilon \right) ${\em -contact metric manifold }%
\cite{Takahashi-69}. For $\varepsilon =1$ and $g$ Riemannian, $M$ is the
usual{\em \ contact metric manifold }\cite{Blair-76}. A contact metric
manifold with $\xi \in N(k)$, is called a {\em $N(k)$-contact metric
manifold }\cite{Blair-Kim-Tripathi-05}. \medskip

An $\left( \varepsilon \right) $-almost contact metric structure $(\varphi
,\xi ,\eta ,g,\varepsilon )$ is called an $\left( \varepsilon \right) ${\em %
-Sasakian structure} if
\[
\left( \nabla _{X}\varphi \right) Y=g(X,Y)\xi -\varepsilon \eta \left(
Y\right) X,
\]%
where $\nabla $ is Levi-Civita connection with respect to the metric $g$. A
manifold endowed with an $\left( \varepsilon \right) $-Sasakian structure is
called an $\left( \varepsilon \right) ${\em -Sasakian manifold }\cite%
{Takahashi-69}. For $\varepsilon =1$ and $g$ Riemannian, $M$ is the usual%
{\em \ Sasakian manifold }\cite{Sasaki-60,Blair-76}. \medskip

An almost contact metric manifold is a {\em Kenmotsu manifold} \cite%
{Kenmotsu-72} if
\begin{equation}
\left( \nabla _{X}\varphi \right) Y=g(\varphi X,Y)\xi -\eta \left( Y\right)
\varphi X.  \label{eq-str-8}
\end{equation}%
By (\ref{eq-str-8}), we have
\begin{equation}
\nabla _{X}\xi =X-\eta (X)\xi .  \label{eq-str-9}
\end{equation}

\subsection*{Almost paracontact manifolds}

Let $M$ be an $n$-dimensional {\em almost paracontact manifold} \cite%
{Sato-76} equipped with an {\em almost paracontact structure} $\left(
\varphi ,\xi ,\eta \right) $, where $\varphi $, $\xi $ and $\eta $ are
tensor fields of type $(1,1)$, $(1,0)$ and $(0,1)$, respectively; and
satisfy the conditions
\begin{equation}
\varphi ^{2}=I-\eta \otimes \xi ,  \label{eq-str-11}
\end{equation}%
\begin{equation}
\eta (\xi )=1.  \label{eq-str-12}
\end{equation}%
Let $g$ be a semi-Riemannian metric on $M$ such that
\begin{equation}
g\left( \varphi X,\varphi Y\right) =g\left( X,Y\right) -\varepsilon \eta
(X)\eta \left( Y\right) ,  \label{eq-str-13}
\end{equation}%
where $\varepsilon =\pm 1$. Then $\left( M,g\right) $ is an $\left(
\varepsilon \right) ${\em -almost paracontact metric manifold} equipped with
an $\left( \varepsilon \right) ${\em -almost paracontact metric structure} $%
(\varphi ,\xi ,\eta ,g,\varepsilon )$. In particular, if ${\rm index}(g)=1$,
then an $(\varepsilon )$-almost paracontact metric manifold is said to be a
{\em Lorentzian almost paracontact manifold}. In particular, if the metric $%
g $ is positive definite, then an $(\varepsilon )$-almost paracontact metric
manifold is the usual {\em almost paracontact metric manifold} \cite{Sato-76}%
. \medskip

The equation (\ref{eq-str-13}) is equivalent to
\begin{equation}
g\left( X,\varphi Y\right) =g\left( \varphi X,Y\right)  \label{eq-str-14}
\end{equation}%
along with
\begin{equation}
g\left( X,\xi \right) =\varepsilon \eta (X).  \label{eq-str-15}
\end{equation}%
From (\ref{eq-str-12}) and (\ref{eq-str-15}), we have
\begin{equation}
g\left( \xi ,\xi \right) =\varepsilon .  \label{eq-str-16}
\end{equation}

An $\left( \varepsilon \right) $-almost paracontact metric structure is
called an $\left( \varepsilon \right) ${\em -para-Sasakian structure} \cite%
{TKYK-09} if
\begin{equation}
\left( \nabla _{X}\varphi \right) Y=-\,g(\varphi X,\varphi Y)\xi
-\varepsilon \eta \left( Y\right) \varphi ^{2}X,  \label{eq-str-17}
\end{equation}%
where $\nabla $ is Levi-Civita connection with respect to the metric $g$. A
manifold endowed with an $\left( \varepsilon \right) $-para-Sasakian
structure is called an $\left( \varepsilon \right) ${\em -para-Sasakian
manifold} \cite{TKYK-09}. For $\varepsilon =1$ and $g$ Riemannian, $M$ is
the usual para-Sasakian manifold \cite{Sato-76}. For $\varepsilon =-1$, $g$
Lorentzian and $\xi $ replaced by $-\xi $, $M$ becomes a Lorentzian
para-Sasakian manifold \cite{Matsumoto-89}. \bigskip

\begin{ex}
{\rm \cite{TG}} The following are some well known examples of $(N(k),\xi )$%
-semi-Riemannian manifolds:

\begin{enumerate}
\item An $N(k)$-contact metric manifold {\rm \cite{Blair-Kim-Tripathi-05}}
is an $(N(k),\xi )$-Riemannian manifold.

\item A Sasakian manifold {\rm \cite{Sasaki-60}} is an $(N(1),\xi )$%
-Riemannian manifold.

\item A Kenmotsu manifold {\rm \cite{Kenmotsu-72}} is an $(N(-1),\xi )$%
-Riemannian manifold.

\item An $(\varepsilon )$-Sasakian manifold {\rm \cite{Takahashi-69}} an $%
(N(\varepsilon ),\xi )$-semi-Riemannian manifold.

\item A para-Sasakian manifold {\rm \cite{Sato-76}} is an $(N(-1),\xi )$%
-Riemannian manifold.

\item An $(\varepsilon )$-para-Sasakian manifold {\rm \cite{TKYK-09}} is an $%
(N(-\varepsilon ),\xi )$-semi-Riemannian manifold.
\end{enumerate}
\end{ex}

In an $n$-dimensional $\left( N(k),\xi \right) $-semi-Riemannian manifold $%
(M,g)$, it is easy to verify that
\begin{equation}
R(X,Y)\xi =\varepsilon k(\eta (Y)X-\eta (X)Y),  \label{eq-curvature}
\end{equation}%
\begin{equation}
R(\xi ,X)Y=\varepsilon k(\varepsilon g(X,Y)\xi -\eta (Y)X),
\label{eq-curvature-2}
\end{equation}%
\begin{equation}
R(\xi ,X)\xi =\varepsilon k(\eta (X)\xi -X),  \label{eq-curvature-3}
\end{equation}%
\begin{equation}
R\left( X,Y,Z,\xi \right) =\varepsilon k(\,\eta \left( X\right) g\left(
Y,Z\right) -\eta \left( Y\right) g\left( X,Z\right) ),
\label{eq-eps-PS-R(X,Y,Z,xi)}
\end{equation}%
\begin{equation}
\eta \left( R\left( X,Y\right) Z\right) =k(\eta \left( X\right) g\left(
Y,Z\right) -\eta \left( Y\right) g\left( X,Z\right) ),
\label{eq-eps-PS-eta(R(X,Y),Z)}
\end{equation}%
\begin{equation}
S(X,\xi )=\varepsilon k(n-1)\eta (X),  \label{eq-ricci}
\end{equation}%
\begin{equation}
Q\xi =k(n-1)\xi ,  \label{eq-Q}
\end{equation}%
\begin{equation}
S(\xi ,\xi )=\varepsilon k(n-1),  \label{eq-S-xi-xi}
\end{equation}%
\begin{equation}
\eta (QX)=\varepsilon g(QX,\xi )=\varepsilon S(X,\xi )=k(n-1)\eta (X).
\label{eq-eta-QX}
\end{equation}%
Moreover, define
\begin{equation}
S^{\ell }(X,Y)=g(Q^{\ell }X,Y)=S(Q^{\ell -1}X,Y),  \label{eq-S-p}
\end{equation}%
where $\ell =0,1,2,\ldots $ and $S^{0}=g$. Using (\ref{eq-eta-QX}) in (\ref%
{eq-S-p}), we get
\begin{equation}
S^{\ell }(X,\xi )=\varepsilon k^{\ell }(n-1)^{\ell }\eta (X).
\label{eq-Sp-QX-xi}
\end{equation}

Now, we give the following Lemma.

\begin{lem}
\label{GCT} {\rm \cite{TG}} Let $M$ be an $n$-dimensional $\left( N(k),\xi
\right) $-semi-Riemannian manifold. Then
\begin{eqnarray}
{\cal T}_{a}(X,Y)\xi &=&(-\varepsilon ka_{0}+\varepsilon
k(n-1)a_{2}-\varepsilon a_{7}\,r)\eta (X)Y  \nonumber \\
&&+\,(\varepsilon ka_{0}+\varepsilon k(n-1)a_{1}+\varepsilon a_{7}\,r)\,\eta
(Y)X  \nonumber \\
&&+\,a_{3}\,S(X,Y)\xi +\varepsilon a_{4}\,\eta (Y)QX  \nonumber \\
&&+\,\varepsilon a_{5}\eta (X)QY+k(n-1)a_{6}g(X,Y)\xi ,  \label{eq-X-Y-xi}
\end{eqnarray}%
\begin{eqnarray}
{\cal T}_{\!a}(\xi ,X)\xi &=&(-\varepsilon ka_{0}\,+\,\varepsilon
k(n-1)a_{2}-\varepsilon a_{7}\,r)X+\varepsilon a_{5}\,QX  \nonumber \\
&&+\left\{ \varepsilon ka_{0}+\varepsilon k(n-1)a_{1}+\varepsilon
k(n-1)a_{3}\right.  \nonumber \\
&&\quad \quad +\left. \varepsilon k(n-1)a_{4}+\varepsilon
k(n-1)a_{6}\,+\varepsilon a_{7}\,r\right\} \eta (X)\xi ,  \label{eq-xi-X-xi}
\end{eqnarray}%
\begin{eqnarray}
{\cal T}_{\!a}(\xi ,Y)Z &=&(ka_{0}+k(n-1)a_{4}+a_{7}r)g\left( Y,Z\right) \xi
\,  \nonumber \\
&&+\,a_{1}\,S\left( Y,Z\right) \xi +\varepsilon k(n-1)a_{3}\eta (Y)Z
\nonumber \\
&&+\,\varepsilon a_{5}\,\eta (Z)QY+\varepsilon a_{6}\,\eta (Y)QZ  \nonumber
\\
&&+\,(-\varepsilon ka_{0}+\varepsilon k(n-1)a_{2}\,-\varepsilon
a_{7}\,r)\eta (Z)Y,  \label{eq-xi-Y-Z}
\end{eqnarray}%
\begin{eqnarray}
\eta ({\cal T}_{\!a}\left( X,Y)\xi \right) &=&\varepsilon
k(n-1)(a_{1}+a_{2}+a_{4}+a_{5})\eta (X)\eta (Y)  \nonumber \\
&&+\,a_{3}\,S(X,Y)+\,k(n-1)a_{6}g(X,Y),  \label{eq-eta-xi-X-Y}
\end{eqnarray}%
\begin{eqnarray}
{\cal T}_{\!a}\left( X,Y,\xi ,V\right) &=&(-\varepsilon ka_{0}\,+\varepsilon
k(n-1)a_{2}\,-\varepsilon a_{7}\,r)\eta (X)g(Y,V)  \nonumber \\
&&+\ (\varepsilon ka_{0}+\varepsilon k(n-1)a_{1}+\varepsilon a_{7}\,r)\,\eta
(Y)g(X,V)  \nonumber \\
&&+\,\varepsilon a_{3}\,S(X,Y)\eta (V)+\varepsilon a_{4}\,\eta (Y)S(X,V)
\nonumber \\
&&+\,\varepsilon a_{5}\,\eta (X)S(Y,V)+\varepsilon k(n-1)a_{6}\,g(X,Y)\eta
(V),  \label{eq-X-Y-xi-V}
\end{eqnarray}%
\begin{eqnarray}
{\cal T}_{\!a}(X,\xi )\xi &=&\left\{ -\varepsilon ka_{0}\,+\varepsilon
k(n-1)a_{2}\,+\varepsilon k(n-1)a_{3}\right.  \nonumber \\
&&\quad +\left. \varepsilon k(n-1)a_{5}+\varepsilon
k(n-1)a_{6}\,\,-\varepsilon a_{7}\,r\right\} \eta (X)\xi  \nonumber \\
&&+\ (\varepsilon ka_{0}+\varepsilon k(n-1)a_{1}+\varepsilon
a_{7}\,r)\,X+\varepsilon a_{4}\,QX,  \label{eq-X-xi-xi}
\end{eqnarray}%
\begin{eqnarray}
S_{{\cal T}_{\!a}}(X,\xi ) &=&\left\{ \varepsilon
k(n-1)(a_{0}+na_{1}+a_{2}+a_{3}+a_{5}+a_{6})\right.  \nonumber \\
&&\left. +\,\varepsilon r(a_{4}+(n-1)a_{7})\right\} \eta (X),
\label{eq-ric-T1}
\end{eqnarray}%
\begin{eqnarray}
S_{{\cal T}_{\!a}}(\xi ,\xi ) &=&\varepsilon
k(n-1)(a_{0}+na_{1}+a_{2}+a_{3}+a_{5}+a_{6})  \nonumber \\
&&+\,\varepsilon r(a_{4}+(n-1)a_{7}).  \label{eq-ric-T2}
\end{eqnarray}
\end{lem}

\begin{rem}
The relations {\rm (\ref{eq-curvature}) -- (\ref{eq-ric-T2})} are true for

\begin{enumerate}
\item a $N(k)$-contact metric manifold {\rm \cite{Blair-Kim-Tripathi-05}}\ ($%
\varepsilon =1$),

\item a Sasakian manifold {\rm \cite{Sasaki-60}} ($k=1$, $\varepsilon =1$),

\item a Kenmotsu manifold {\rm \cite{Kenmotsu-72}} ($k=-1$, $\varepsilon =1$%
),

\item an $(\varepsilon )$-Sasakian manifold {\rm \cite{Takahashi-69}} ($%
k=\varepsilon $, $\varepsilon k=1$),

\item a para-Sasakian manifold {\rm \cite{Sato-76}} ($k=-1$, $\varepsilon =1$%
), and

\item an $(\varepsilon )$-para-Sasakian manifold {\rm \cite{TKYK-09}} ($%
k=-\,\varepsilon $, $\varepsilon k=-\,1$).
\end{enumerate}

\noindent Even, all the relations and results of this paper will be true for
the above six cases.
\end{rem}

\section{$({\cal T}_{\!a},{\cal T}_{\!b})$-pseudosymmetry\label{sect-TTP}}

\begin{defn-new}
A semi-Riemannian manifold $\left( M,g\right) $ is said to be $({\cal T}%
_{\!a},{\cal T}_{\!b})$-pseudosymmetric if
\begin{equation}
{\cal T}_{\!a}\cdot {\cal T}_{\!b}=L_{g}Q(g,{\cal T}_{\!b}),
\label{eq-T.T=Lg}
\end{equation}%
where $L_{g}$ is some smooth function on $M$. In particular, it is said to
be $(R,{\cal T}_{\!a})$-pseudosymmetric or, in brief, ${\cal T}_{\!a}$%
-pseudosymmetric if
\begin{equation}
R\cdot {\cal T}_{\!a}=L_{g}Q(g,{\cal T}_{\!a})  \label{eq-R.T=Lg}
\end{equation}%
holds on the set
\[
{\cal U}_{g}=\left\{ x\in M\,:\left( {\cal T}_{\!a}\right) _{x}\not=0\;{\rm %
at}\;x\right\} ,
\]%
where $L_{g}$ is some smooth function on ${\cal U}_{g}$. In particular, if
in (\ref{eq-R.T=Lg}), {\em ${\cal T}_{\!a}$} is equal to $R$, ${\cal C}%
_{\ast }$, ${\cal C}$, ${\cal L}$, ${\cal V}$, ${\cal P}_{\ast }$, ${\cal P}$%
, ${\cal M}$, ${\cal W}_{0}$, ${\cal W}_{0}^{\ast }$, ${\cal W}_{1}$, ${\cal %
W}_{1}^{\ast }$, ${\cal W}_{2}$, ${\cal W}_{3}$, ${\cal W}_{4}$, ${\cal W}%
_{5}$, ${\cal W}_{6}$, ${\cal W}_{7}$, ${\cal W}_{8}$, ${\cal W}_{9}$, then
it becomes pseudosymmetric, quasi-conformal pseudosymmetric, Weyl
pseudosymmetric, conharmonic pseudosymmetric, concircular pseudosymmetric,
pseudo-projective pseudosymmetric, projective pseudosymmetric, ${\cal M}$%
-pseudosymmetric, ${\cal W}_{0}$-pseudosymmetric, ${\cal W}_{0}^{\ast }$%
-pseudosymmetric, ${\cal W}_{1}$-pseudosymmetric, ${\cal W}_{1}^{\ast }$%
-pseudosymmetric, ${\cal W}_{2}$-pseudosym\allowbreak metric, ${\cal W}_{3}$%
-pseudosymmetric, ${\cal W}_{4}$-pseudosymmetric, ${\cal W}_{5}$%
-pseudosymmetric, ${\cal W}_{6}$-pseudosymmetric, ${\cal W}_{7}$%
-pseu\allowbreak dosymmetric, ${\cal W}_{8}$-pseudosymmetric, ${\cal W}_{9}$%
-pseudosymmetric, respectively.
\end{defn-new}

\begin{th}
\label{th-T-T-11} Let $M$ be an $n$-dimensional $({\cal T}_{\!a},{\cal T}%
_{\!b})$-pseudosymmetric $(N(k),\xi )$-semi-Riemannian manifold. Then
\begin{eqnarray}
&&\varepsilon b_{0}(ka_{0}+\varepsilon
k(n-1)a_{4}+a_{7}r)R(U,V,W,X)+\varepsilon a_{1}b_{0}S(X,R(U,V)W)  \nonumber
\\
&&-2k(n-1)a_{3}(kb_{0}+k(n-1)b_{4}+b_{7}r)\eta (X)\eta (U)g(V,W)  \nonumber
\\
&&-2k(n-1)a_{3}(-kb_{0}+k(n-1)b_{5}-b_{7}r)\eta (X)\eta (V)g(U,W)  \nonumber
\\
&&+\varepsilon a_{1}b_{4}S^{2}(X,U)g(V,W)+\varepsilon
a_{1}b_{5}S^{2}(X,V)g(U,W)  \nonumber \\
&&+\varepsilon a_{1}b_{6}S^{2}(X,W)g(U,V)-a_{5}(b_{1}+b_{3})S^{2}(X,V)\eta
(U)\eta (W)  \nonumber \\
&&-a_{5}(b_{1}+b_{2})S^{2}(X,W)\eta (U)\eta
(V)-a_{5}(b_{2}+b_{3})S^{2}(X,U)\eta (V)\eta (W)  \nonumber \\
&&-2a_{6}b_{1}S^{2}(V,W)\eta (X)\eta (U)-2a_{6}b_{2}S^{2}(U,W)\eta (X)\eta
(V)  \nonumber \\
&&-2a_{6}b_{3}S^{2}(U,V)\eta (X)\eta (W)-2k^{2}(n-1)a_{3}b_{6}g(U,V)\eta
(X)\eta (W)  \nonumber \\
&&-2\left( k(n-1)a_{3}b_{1}+a_{6}(kb_{0}+k(n-1)b_{4}+b_{7}r)\right) \eta
(X)\eta (U)S(V,W)  \nonumber \\
&&-2\left( k(n-1)a_{3}b_{2}+a_{6}(-kb_{0}+k(n-1)b_{5}-b_{7}r)\right) \eta
(X)\eta (V)S(U,W)  \nonumber \\
&&-2k(n-1)(a_{3}b_{3}+a_{6}b_{6})S(U,V)\eta (X)\eta (W)  \nonumber \\
&&+\varepsilon \left(
b_{4}(ka_{0}+k(n-1)a_{4}+a_{7}r)-a_{1}(kb_{0}+k(n-1)b_{4})\right)
S(X,U)g(V,W)  \nonumber \\
&&+\varepsilon \left(
b_{5}(ka_{0}+k(n-1)a_{4}+a_{7}r)-a_{1}(-kb_{0}+k(n-1)b_{5})\right)
S(X,V)g(U,W)  \nonumber \\
&&+\varepsilon b_{6}(ka_{0}+k(n-1)(a_{4}-a_{1})+a_{7}r)S(X,W)g(U,V)
\nonumber \\
&&-\,\varepsilon (kb_{0}+k(n-1)b_{4})(ka_{0}+k(n-1)a_{4}+a_{7}r)g(X,U)g(V,W)
\nonumber \\
&&-\,\varepsilon (-kb_{0}+k(n-1)b_{5})(ka_{0}+k(n-1)a_{4}+a_{7}r)g(U,W)g(X,V)
\nonumber \\
&&-\,\varepsilon k(n-1)b_{6}(ka_{0}+k(n-1)a_{4}+a_{7}r)g(X,W)g(U,V)
\nonumber \\
&&-\,k(n-1)\left( (b_{2}+b_{3})(ka_{0}+k(n-1)a_{4}+a_{7}r)\right.  \nonumber
\\
&&\left. +(a_{2}+a_{4})(-kb_{0}+k(n-1)(b_{5}+b_{6})-b_{7}r)\right)
g(X,U)\eta (V)\eta (W)  \nonumber \\
&&-\,k(n-1)\left( (b_{1}+b_{3})(ka_{0}+k(n-1)a_{4}+a_{7}r)\right.  \nonumber
\\
&&\left. +(a_{2}+a_{4})(kb_{0}+k(n-1)(b_{4}+b_{6})+b_{7}r)\right) g(X,V)\eta
(U)\eta (W)  \nonumber \\
&&-\,\left( (b_{1}+b_{3})(-ka_{0}+k(n-1)(a_{1}+a_{2})-a_{7}r)\right.
\nonumber \\
&&\left. +(a_{1}+a_{5})(kb_{0}+k(n-1)(b_{4}+b_{6})+b_{7}r)\right) S(X,V)\eta
(U)\eta (W)  \nonumber \\
&&-\,\left( (b_{2}+b_{3})(-ka_{0}+k(n-1)(a_{1}+a_{2})-a_{7}r)\right.
\nonumber \\
&&\left. +(a_{1}+a_{5})(kb_{0}+k(n-1)(b_{5}+b_{6})+b_{7}r)\right) S(X,U)\eta
(V)\eta (W)  \nonumber \\
&&-\,\left( (b_{1}+b_{2})(-ka_{0}+k(n-1)(a_{1}+a_{2})-a_{7}r)\right.
\nonumber \\
&&\left. +k(n-1)(b_{4}+b_{5})(a_{1}+a_{5})\right) S(X,W)\eta (U)\eta (V)
\nonumber \\
&&-\,k(n-1)\left( k(n-1)(b_{4}+b_{5})(a_{2}+a_{4})\right.  \nonumber \\
&&\left. +(b_{1}+b_{2})(ka_{0}+k(n-1)a_{4}+a_{7}r)\right) g(X,W)\eta (U)\eta
(V)  \nonumber \\
&=&\,L_{g}(\varepsilon b_{0}R(U,V,W,X)+\varepsilon
b_{4}S(X,U)g(V,W)+\varepsilon b_{5}S(X,V)g(U,W)  \nonumber \\
&&+\varepsilon b_{6}S(X,W)g(U,V)-\varepsilon k(n-1)b_{6}g(U,V)g(X,W)
\nonumber \\
&&-k(n-1)(b_{2}+b_{3})g(X,U)\eta (V)\eta (W)  \nonumber \\
&&-k(n-1)(b_{1}+b_{3})g(X,V)\eta (U)\eta (W)  \nonumber \\
&&-k(n-1)(b_{1}+b_{2})g(X,W)\eta (U)\eta (V)  \nonumber \\
&&+(b_{2}+b_{3})S(X,U)\eta (V)\eta (W)  \nonumber \\
&&+(b_{1}+b_{3})S(X,V)\eta (U)\eta (W)  \nonumber \\
&&+(b_{1}+b_{2})S(X,W)\eta (U)\eta (V)  \nonumber \\
&&-\varepsilon (-kb_{0}+k(n-1)b_{5})g(U,W)g(X,V)  \nonumber \\
&&-\varepsilon (kb_{0}+k(n-1)b_{4})g(V,W)g(X,U)).  \label{eq-T-T-2-i}
\end{eqnarray}
\end{th}

\noindent {\bf Proof.} Let $M$ be an $n$-dimensional $({\cal T}_{\!a},{\cal T%
}_{\!b})$-pseudosymmetric $(N(k),\xi )$-semi-Riemannian manifold. Then
\begin{equation}
{\cal T}_{\!a}(Z,X)\cdot {\cal T}_{\!b}(U,V)W=L_{g}Q(g,{\cal T}%
_{\!b})(U,V,W;Z,X).  \label{eq-txxi-1}
\end{equation}%
Taking $Z=\xi $ in (\ref{eq-txxi-1}), we get
\[
{\cal T}_{\!a}(\xi ,X)\cdot {\cal T}_{\!b}(U,V)W=L_{g}Q(g,{\cal T}%
_{\!b})(U,V,W;\xi ,X),
\]%
which gives
\begin{eqnarray*}
&&{\cal T}_{\!a}(\xi ,X){\cal T}_{\!b}(U,V)W-{\cal T}_{\!b}({\cal T}%
_{\!a}(\xi ,X)U,V)W \\
&&-\ {\cal T}_{\!b}(U,{\cal T}_{\!a}(\xi ,X)V)W-{\cal T}_{\!b}(U,V){\cal T}%
_{\!a}(\xi ,X)W \\
&=&L_{g}((\xi \wedge X){\cal T}_{\!b}(U,V)W-{\cal T}_{\!b}((\xi \wedge
X)U,V)W \\
&&-{\cal T}_{\!b}(U,(\xi \wedge X)V)W-{\cal T}_{\!b}(U,V)(\xi \wedge X)W),
\end{eqnarray*}%
that is,
\begin{eqnarray}
&&{\cal T}_{\!a}(\xi ,X){\cal T}_{\!b}(U,V)W-{\cal T}_{\!b}({\cal T}%
_{\!a}(\xi ,X)U,V)W  \nonumber \\
&&\quad -\ {\cal T}_{\!b}(U,{\cal T}_{\!a}(\xi ,X)V)W-{\cal T}_{\!b}(U,V)%
{\cal T}_{\!a}(\xi ,X)W  \nonumber \\
&=&L_{g}({\cal T}_{\!b}(U,V,W,X)\xi -{\cal T}_{\!b}(U,V,W,\xi )X  \nonumber
\\
&&-g(X,U){\cal T}_{\!b}(\xi ,V)W+\varepsilon \eta (U){\cal T}_{\!b}(X,V)W
\nonumber \\
&&-g(X,V){\cal T}_{\!b}(U,\xi )W+\varepsilon \eta (V){\cal T}_{\!b}(U,X)W
\nonumber \\
&&-g(X,W){\cal T}_{\!b}(U,V)\xi +\varepsilon \eta (W){\cal T}_{\!b}(U,V)X).
\label{eq-T-T-i}
\end{eqnarray}%
Taking the inner product of (\ref{eq-T-T-i}) with $\xi $, we get
\begin{eqnarray}
&&{\cal T}_{\!a}(\xi ,X,{\cal T}_{\!b}(U,V)W,\xi )-{\cal T}_{\!b}({\cal T}%
_{\!a}(\xi ,X)U,V,W,\xi )  \nonumber \\
&&\quad -\ {\cal T}_{\!b}(U,{\cal T}_{\!a}(\xi ,X)V,W,\xi )-{\cal T}%
_{\!b}(U,V,{\cal T}_{\!a}(\xi ,X)W,\xi )  \nonumber \\
&=&L_{g}(\varepsilon {\cal T}_{\!b}(U,V,W,X)-\varepsilon \eta (X){\cal T}%
_{\!b}(U,V,W,\xi )  \nonumber \\
&&-g(X,U){\cal T}_{\!b}(\xi ,V,W,\xi )+\varepsilon \eta (U){\cal T}%
_{\!b}(X,V,W,\xi )  \nonumber \\
&&-g(X,V){\cal T}_{\!b}(U,\xi ,W,\xi )+\varepsilon \eta (V){\cal T}%
_{\!b}(U,X,W,\xi )  \nonumber \\
&&-g(X,W){\cal T}_{\!b}(U,V,\xi ,\xi )+\varepsilon \eta (W){\cal T}%
_{\!b}(U,V,X,\xi )).  \label{eq-T-T-1-i}
\end{eqnarray}%
By using (\ref{eq-X-Y-xi}),\ldots ,(\ref{eq-X-xi-xi}) in (\ref{eq-T-T-1-i}),
we get (\ref{eq-T-T-2-i}). $\blacksquare $

\begin{th}
Let $M$ be an $n$-dimensional $({\cal T}_{\!a},{\cal T}_{\!b})$%
-pseudosymmetric $(N(k),\xi )$-semi-Riemannian manifold. Then
\begin{eqnarray}
&&-\varepsilon \left\{
a_{5}b_{0}+na_{5}b_{1}+a_{5}b_{2}+a_{5}b_{6}+a_{5}b_{3}+a_{5}b_{5}\right\}
S^{2}(V,W)  \nonumber \\
&&+\ (\left\{ (nb_{1}+b_{2}+b_{3}+b_{5}+b_{6}+b_{0})(\varepsilon
ka_{0}+\varepsilon b_{7}r)\right.  \nonumber \\
&&\qquad -\ \varepsilon
k(n-1)(2a_{5}b_{6}+a_{2}b_{3}+a_{1}b_{6}+a_{1}b_{3}+a_{1}b_{5}  \nonumber \\
&&\qquad \left. +\
a_{1}b_{1}+a_{1}b_{2}+a_{2}b_{2}+a_{2}b_{6}+na_{2}b_{1}+a_{1}b_{0}+a_{2}b_{0})\right\}
\nonumber \\
&&-\ \varepsilon (n-1)a_{1}b_{7}r-\varepsilon na_{5}b_{7}r-\varepsilon
b_{4}a_{5}r-\varepsilon a_{1}b_{4}r)S(V,W)  \nonumber \\
&&+\ \left\{ -\varepsilon
k(n-1)(nb_{1}+b_{2}+b_{3}+b_{5}+b_{6}+b_{0})(a_{7}r+ka_{0}+k(n-1)a_{4})%
\right.  \nonumber \\
&&\qquad \left. -\ \varepsilon
k(n-1)r((n-1)b_{7}a_{2}+(n-1)b_{7}a_{4}+a_{2}b_{4}+a_{4}b_{4})\right\} g(V,W)
\nonumber \\
&&+\ (a_{1}+a_{2}+2a_{3}+a_{4}+a_{5}+2a_{6})\left\{
-k^{2}(n-1)^{2}(nb_{1}+b_{2}+b_{3}+b_{5}+b_{6}+b_{0})\right.  \nonumber \\
&&\qquad \left. -\ k(n-1)^{2}b_{7}r-k(n-1)b_{4}r\right\} \eta (V)\eta (W)
\nonumber \\
&=&L_{g}(\varepsilon (b_{0}+b_{5}+b_{6})S(V,W)+\varepsilon
(b_{4}r-k(n-1)(b_{0}+nb_{4}+b_{5}+b_{6}))g(V,W)  \nonumber \\
&&+(b_{2}+b_{3})(r-kn(n-1))\eta (V)\eta (W)).  \label{eq-semi-TS-i}
\end{eqnarray}
\end{th}

\begin{th}
\label{GCT-ss-11} Let $M$ be an $n$-dimensional ${\cal T}_{\!a}$%
-pseudosymmetric $(N(k),\xi )$-semi-Riemannian manifold. Then
\begin{eqnarray}
&&\varepsilon a_{0}kR(U,V,W,X)+\varepsilon ka_{4}S(X,U)g(V,W)+\varepsilon
ka_{5}S(X,V)g(U,W)  \nonumber \\
&&+\,\varepsilon ka_{6}S(X,W)g(U,V)-\varepsilon k^{2}(n-1)a_{6}g(X,W)g(U,V)
\nonumber \\
&&-\,\varepsilon k(ka_{0}+k(n-1)a_{4})g(V,W)g(X,U)  \nonumber \\
&&-\,\varepsilon k(-ka_{0}+k(n-1)a_{5})g(U,W)g(X,V)  \nonumber \\
&&-\,k^{2}(n-1)(a_{2}+a_{3})g(X,U)\eta (V)\eta (W)  \nonumber \\
&&-\,k^{2}(n-1)(a_{1}+a_{3})g(X,V)\eta (U)\eta (W)  \nonumber \\
&&-\,k^{2}(n-1)(a_{1}+a_{2})g(X,W)\eta (U)\eta (V)  \nonumber \\
&&+\,k(a_{2}+a_{3})S(X,U)\eta (V)\eta (W)+k(a_{1}+a_{3})S(X,V)\eta (U)\eta
(W)  \nonumber \\
&&+\,k(a_{1}+a_{2})S(X,W)\eta (U)\eta (V)  \nonumber \\
&&=\, L_{g}(\varepsilon a_{0}R(U,V,W,X)+\varepsilon
a_{4}S(X,U)g(V,W)+\varepsilon a_{5}S(X,V)g(U,W)  \nonumber \\
&&+\varepsilon a_{6}S(X,W)g(U,V)-\varepsilon k(n-1)a_{6}g(U,V)g(X,W)
\nonumber \\
&&-k(n-1)(a_{2}+a_{3})g(X,U)\eta (V)\eta (W)  \nonumber \\
&&-k(n-1)(a_{1}+a_{3})g(X,V)\eta (U)\eta (W)  \nonumber \\
&&-k(n-1)(a_{1}+a_{2})g(X,W)\eta (U)\eta (V)  \nonumber \\
&&+(a_{2}+a_{3})S(X,U)\eta (V)\eta (W)  \nonumber \\
&&+(a_{1}+a_{3})S(X,V)\eta (U)\eta (W)  \nonumber \\
&&+(a_{1}+a_{2})S(X,W)\eta (U)\eta (V)  \nonumber \\
&&-\varepsilon (-ka_{0}+k(n-1)a_{5})g(U,W)g(X,V)  \nonumber \\
&&-\varepsilon (ka_{0}+k(n-1)a_{4})g(V,W)g(X,U)).  \label{eq-semi-sym-R-i}
\end{eqnarray}
\end{th}

\begin{rem-new}
Here two cases arise. First is that when $L_{g}=0$. In this case, it is $%
{\cal T}_{\!a}$-semisymmetric. We exclude this case, because it is already
studied in \cite{TG}. And the second case is that when $L_{g}\not=0 $. In
the following Theorem, we consider the result for $L_{g}\not=0$.
\end{rem-new}

\begin{th}
\label{th-TPS} Let $M$ be an $n$-dimensional ${\cal T}_{\!a}$%
-pseudosymmetric $\left( N(k),\xi \right) $-semi-Riemannian manifold such
that $a_{0}+a_{5}+a_{6}\not=0$.

\begin{enumerate}
\item If $a_{0}+a_{2}+a_{3}+na_{4}+a_{5}+a_{6}\not=0$, then either it is
Einstein manifold or $L_{g}=k$.

\item If $a_{0}+a_{2}+a_{3}+na_{4}+a_{5}+a_{6}=0$, then either it is $\eta $%
-Einstein manifold or $L_{g}=k$.
\end{enumerate}
\end{th}

\noindent {\bf Proof.} Let $M$ be an $n$-dimensional ${\cal T}_{\!a}$%
-pseudosymmetric $\left( N(k),\xi \right) $-semi-Riemannian manifold. On
contracting (\ref{eq-semi-sym-R-i}), we get
\begin{eqnarray*}
&&k(\varepsilon (a_{0}+a_{5}+a_{6})S(V,W)+\varepsilon
(a_{4}r-k(n-1)(a_{0}+na_{4}+a_{5}+a_{6}))g(V,W) \\
&&+(a_{2}+a_{3})(r-kn(n-1))\eta (V)\eta (W)) \\
&=&L_{g}(\varepsilon (a_{0}+a_{5}+a_{6})S(V,W)+\varepsilon
(a_{4}r-k(n-1)(a_{0}+na_{4}+a_{5}+a_{6}))g(V,W) \\
&&+(a_{2}+a_{3})(r-kn(n-1))\eta (V)\eta (W)),
\end{eqnarray*}%
which can be rewritten as
\begin{eqnarray}
&&(L_{g}-k)(\varepsilon (a_{0}+a_{5}+a_{6})S(V,W)  \nonumber \\
&&+\varepsilon (a_{4}r-k(n-1)(a_{0}+na_{4}+a_{5}+a_{6}))g(V,W)  \nonumber \\
&&+(a_{2}+a_{3})(r-kn(n-1))\eta (V)\eta (W)).  \label{eq-TPS}
\end{eqnarray}%
On contracting above equation, we get
\[
(L_{g}-k)(a_{0}+a_{2}+a_{3}+na_{4}+a_{5}+a_{6})(r-kn(n-1)).
\]%
{\bf Case 1.} If $a_{0}+a_{2}+a_{3}+na_{4}+a_{5}+a_{6}\not=0$, then either $%
L_{g}=k $ or $r=kn(n-1)$. If $r=kn(n-1)$, then from (\ref{eq-TPS}), we get
\[
S=k(n-1)g.
\]%
{\bf Case 2.} If $a_{0}+a_{2}+a_{3}+na_{4}+a_{5}+a_{6}=0$, then by (\ref%
{eq-TPS}), we get either $L_{g}=k$ or
\begin{eqnarray*}
-\varepsilon (a_{0}+a_{5}+a_{6})S(V,W) &=&\varepsilon
(a_{4}r-k(n-1)(a_{0}+na_{4}+a_{5}+a_{6}))g(V,W) \\
&&+(a_{2}+a_{3})(r-kn(n-1))\eta (V)\eta (W)).
\end{eqnarray*}
This proves the result. $\blacksquare$

\begin{cor}
An $n$-dimensional ${\cal T}_{\!a}$-pseudosymmetric $\left( N(k),\xi \right)
$-semi-Riemannian manifold is of the form $R\cdot {\cal T}_{\!a}=kQ(g,{\cal T%
}_{\!a})$.
\end{cor}

In view of Theorem~\ref{th-TPS}, we have the following

\begin{cor}
\label{cor-pseudo-1} Let $M$ be an $n$-dimensional ${\cal T}_{\!a}$%
-pseudosymmetric $\left( N(k),\xi \right) $-semi-Riemannian manifold such
that
\[
{\cal T}_{\!a}\in \left\{ R,{\cal V},{\cal P},{\cal M},{\cal W}_{0},{\cal W}%
_{0}^{\ast },{\cal W}_{1},{\cal W}_{1}^{\ast },{\cal W}_{3},\ldots ,{\cal W}%
_{8}\right\}
\]%
Then we have the following table\/{\rm :}~
\[
\begin{tabular}{|l|l|l|}
\hline
${\boldmath M}$ & $L_{g}=$ & ${\boldmath S=}$ \\ \hline
$N(k)$-contact metric & $k$ & $k(n-1)g$ \\ \hline
Sasakian & $1$ & $(n-1)g$ \\ \hline
Kenmotsu & $-\,1$ & $-\,(n-1)g$ \\ \hline
$(\varepsilon )$-Sasakian & $\varepsilon $ & $\varepsilon (n-1)g$ \\ \hline
para-Sasakian & $-\,1$ & $-\,(n-1)g$ \\ \hline
$(\varepsilon )$-para-Sasakian & $-\,\varepsilon $ & $-\,\varepsilon (n-1)g$
\\ \hline
\end{tabular}%
\
\]
\end{cor}

\begin{cor}
{\rm \cite{Ozgur-06} }Let $M$ be an $n$-dimensional, $n\geq 3$, Kenmotsu
manifold. If $M$ is pseudosymmetric then either it is locally isometric to
the hyperbolic space $H^{n}(-1)$ or $L_{g}=-1$ holds on $M$.
\end{cor}

\begin{cor}
{\rm \cite{Ozgur-06} }Every Kenmotsu manifold $M$, $n\geq 3$, is a
pseudosymmetric manifold of the form $R\cdot R=-Q(g,R)$.
\end{cor}

\begin{cor}
Let $M$ be an $n$-dimensional quasi-conformal pseudosymmetric $(N(k),\xi )$%
-semi-Riemannian manifold such that $a_{0}-a_{1}\not=0$\ and $%
a_{0}+(n-2)a_{1}\not=0$. Then we have the following table\/{\rm :}~
\[
\begin{tabular}{|l|l|l|}
\hline
${\boldmath M}$ & $L_{g}=$ & ${\boldmath S=}$ \\ \hline
$N(k)$-contact metric & $k$ & $k(n-1)g$ \\ \hline
Sasakian & $1$ & $(n-1)g$ \\ \hline
Kenmotsu & $-\,1$ & $-\,(n-1)g$ \\ \hline
$(\varepsilon )$-Sasakian & $\varepsilon $ & $\varepsilon (n-1)g$ \\ \hline
para-Sasakian & $-\,1$ & $-\,(n-1)g$ \\ \hline
$(\varepsilon )$-para-Sasakian & $-\,\varepsilon $ & $-\,\varepsilon (n-1)g$
\\ \hline
\end{tabular}%
\
\]
\end{cor}

\begin{cor}
Let $M$ be an $n$-dimensional pseudo-projective pseudosymmetric $\left(
N(k),\xi \right) $-semi-Riemannian manifold such that $a_{0}\not=0$ and $%
a_{0}-a_{1}\not=0$. Then we have the following table\/{\rm :}~
\[
\begin{tabular}{|l|l|l|}
\hline
${\boldmath M}$ & $L_{g}=$ & ${\boldmath S=}$ \\ \hline
$N(k)$-contact metric & $k$ & $k(n-1)g$ \\ \hline
Sasakian & $1$ & $(n-1)g$ \\ \hline
Kenmotsu & $-\,1$ & $-\,(n-1)g$ \\ \hline
$(\varepsilon )$-Sasakian & $\varepsilon $ & $\varepsilon (n-1)g$ \\ \hline
para-Sasakian & $-\,1$ & $-\,(n-1)g$ \\ \hline
$(\varepsilon )$-para-Sasakian & $-\,\varepsilon $ & $-\,\varepsilon (n-1)g$
\\ \hline
\end{tabular}%
\
\]
\end{cor}

\begin{cor}
Let $M$ be an $n$-dimensional Weyl pseudosymmetric $\left( N(k),\xi \right) $%
-semi-Riemannian manifold. Then we have the following table\/{\rm :}~\
\[
\begin{tabular}{|l|l|l|}
\hline
${\boldmath M}$ & $L_{g}=$ & ${\boldmath S=}$ \\ \hline
$N(k)$-contact metric & $k$ & $S=\left( \dfrac{r}{n-1}-k\right) g+\left( nk-%
\dfrac{r}{n-1}\right) \eta \otimes \eta $ \\ \hline
Sasakian & $1$ & $S=\left( \dfrac{r}{n-1}-1\right) g+\left( n-\dfrac{r}{n-1}%
\right) \eta \otimes \eta $ \\ \hline
Kenmotsu {\rm \cite{Ozgur-06}} & $-\,1$ & $S=\left( \dfrac{r}{n-1}+1\right)
g-\left( n+\dfrac{r}{n-1}\right) \eta \otimes \eta $ \\ \hline
$(\varepsilon )$-Sasakian & $\varepsilon $ & $S=\left( \dfrac{r}{n-1}%
-\varepsilon \right) g+\varepsilon \left( n\varepsilon -\dfrac{r}{n-1}%
\right) \eta \otimes \eta $ \\ \hline
para-Sasakian {\rm \cite{Ozgur-2005}} & $-\,1$ & $S=\left( \dfrac{r}{n-1}%
+1\right) g-\left( n+\dfrac{r}{n-1}\right) \eta \otimes \eta $ \\ \hline
$(\varepsilon )$-para-Sasakian & $-\,\varepsilon $ & $S=\left( \dfrac{r}{n-1}%
+\varepsilon \right) g-\varepsilon \left( n\varepsilon +\dfrac{r}{n-1}%
\right) \eta \otimes \eta $ \\ \hline
\end{tabular}%
\]
\end{cor}

\begin{cor}
Let $M$ be an $n$-dimensional conharmonic pseudosymmetric $\left( N(k),\xi
\right) $-semi-Riemannian manifold. Then we have the following table\/{\rm :}%
~
\[
\begin{tabular}{|l|l|l|}
\hline
${\boldmath M}$ & $L_{g}=$ & ${\boldmath S=}$ \\ \hline
$N(k)$-contact metric & $k$ & $S=\left( \dfrac{r}{n-1}-k\right) g+\left( nk-%
\dfrac{r}{n-1}\right) \eta \otimes \eta $ \\ \hline
Sasakian & $1$ & $S=\left( \dfrac{r}{n-1}-1\right) g+\left( n-\dfrac{r}{n-1}%
\right) \eta \otimes \eta $ \\ \hline
Kenmotsu & $-\,1$ & $S=\left( \dfrac{r}{n-1}+1\right) g-\left( n+\dfrac{r}{%
n-1}\right) \eta \otimes \eta $ \\ \hline
$(\varepsilon )$-Sasakian & $\varepsilon $ & $S=\left( \dfrac{r}{n-1}%
-\varepsilon \right) g+\varepsilon \left( n\varepsilon -\dfrac{r}{n-1}%
\right) \eta \otimes \eta $ \\ \hline
para-Sasakian & $-\,1$ & $S=\left( \dfrac{r}{n-1}+1\right) g-\left( n+\dfrac{%
r}{n-1}\right) \eta \otimes \eta $ \\ \hline
$(\varepsilon )$-para-Sasakian & $-\,\varepsilon $ & $S=\left( \dfrac{r}{n-1}%
+\varepsilon \right) g-\varepsilon \left( n\varepsilon +\dfrac{r}{n-1}%
\right) \eta \otimes \eta $ \\ \hline
\end{tabular}%
\
\]
\end{cor}

\begin{cor}
Let $M$ be an $n$-dimensional ${\cal W}_{2}$-pseudosymmetric $\left(
N(k),\xi \right) $-semi-Riemannian manifold. Then we have the following
table\/{\rm :}~
\[
\begin{tabular}{|l|l|l|}
\hline
${\boldmath M}$ & $L_{g}=$ & ${\boldmath S=}$ \\ \hline
$N(k)$-contact metric & $k$ & $\dfrac{r}{n}g$ \\ \hline
Sasakian & $1$ & $\dfrac{r}{n}g$ \\ \hline
Kenmotsu & $-\,1$ & $\dfrac{r}{n}g$ \\ \hline
$(\varepsilon )$-Sasakian & $\varepsilon $ & $\dfrac{r}{n}g$ \\ \hline
para-Sasakian & $-\,1$ & $\dfrac{r}{n}g$ \\ \hline
$(\varepsilon )$-para-Sasakian & $-\,\varepsilon $ & $\dfrac{r}{n}g$ \\
\hline
\end{tabular}%
\
\]
\end{cor}

\begin{cor}
Let $M$ be an $n$-dimensional ${\cal W}_{9}$-pseudosymmetric $\left(
N(k),\xi \right) $-semi-Riemannian manifold. Then we have the following
table\/{\rm :}~
\[
\begin{tabular}{|l|l|l|}
\hline
${\boldmath M}$ & $L_{g}=$ & ${\boldmath S=}$ \\ \hline
$N(k)$-contact metric & $k$ & $S=\left( \dfrac{r}{n-1}-k\right) g+\left( nk-%
\dfrac{r}{n-1}\right) \eta \otimes \eta $ \\ \hline
Sasakian & $1$ & $S=\left( \dfrac{r}{n-1}-1\right) g+\left( n-\dfrac{r}{n-1}%
\right) \eta \otimes \eta $ \\ \hline
Kenmotsu & $-\,1$ & $S=\left( \dfrac{r}{n-1}+1\right) g-\left( n+\dfrac{r}{%
n-1}\right) \eta \otimes \eta $ \\ \hline
$(\varepsilon )$-Sasakian & $\varepsilon $ & $S=\left( \dfrac{r}{n-1}%
-\varepsilon \right) g+\varepsilon \left( n\varepsilon -\dfrac{r}{n-1}%
\right) \eta \otimes \eta $ \\ \hline
para-Sasakian & $-\,1$ & $S=\left( \dfrac{r}{n-1}+1\right) g-\left( n+\dfrac{%
r}{n-1}\right) \eta \otimes \eta $ \\ \hline
$(\varepsilon )$-para-Sasakian & $-\,\varepsilon $ & $S=\left( \dfrac{r}{n-1}%
+\varepsilon \right) g-\varepsilon \left( n\varepsilon +\dfrac{r}{n-1}%
\right) \eta \otimes \eta $ \\ \hline
\end{tabular}%
\
\]
\end{cor}

\section{$({\cal T}_{\!a},{\cal T}_{\!b},S^{\ell })$-pseudosymmetry\label%
{sect-TTSP}}

\begin{defn-new}
A semi-Riemannian manifold is said to be $({\cal T}_{\!a},{\cal T}%
_{\!b},S^{\ell })$-pseudosymmetric if it satisfies
\begin{equation}
{\cal T}_{\!a}\cdot {\cal T}_{\!b}=L_{S^{\ell }}Q(S^{\ell },{\cal T}_{\!b}),
\label{eq-T.T=LSl}
\end{equation}%
where $L_{S^{\ell }}$ is some smooth function on $M$. In particular, it is
said to be $(R,{\cal T}_{\!a},S^{\ell })$-pseudosymmetric if
\begin{equation}
R\cdot {\cal T}_{\!a}=L_{S^{\ell }}Q(S^{\ell },{\cal T}_{\!a})
\label{eq-R.T=LSl}
\end{equation}%
holds on the set ${\cal U}_{S^{\ell }}=\left\{ x\in M:Q(S^{\ell },{\cal T}%
_{\!a})\not=0\right\} $, where $L_{S^{\ell }}$ is some smooth function on $%
{\cal U}_{S^{\ell }}$.
\end{defn-new}

For $\ell =1$, we can give the following definition.

\begin{defn-new}
A semi-Riemannian manifold is called $({\cal T}_{\!a},{\cal T}_{\!b},S)$%
-pseudosymmetric if it satisfies
\begin{equation}
{\cal T}_{\!a}\cdot {\cal T}_{\!b}=L_{S}Q(S,{\cal T}_{\!b}),
\label{eq-T.T=LS}
\end{equation}%
where $L_{S}$ is some smooth function on $M$. In particular, it is said to
be $(R,{\cal T}_{\!a},S)$-pseudosymmetric if
\begin{equation}
R\cdot {\cal T}_{\!a}=L_{S}Q(S,{\cal T}_{\!a})  \label{eq-R.T=LS}
\end{equation}%
holds on the set ${\cal U}_{S}=\left\{ x\in M:Q(S,{\cal T}%
_{\!a})\not=0\right\} $, where $L_{S}$ is some smooth function on ${\cal U}%
_{S}$.
\end{defn-new}

\begin{rem-new}
A semi-Riemannian manifold is said to be $(R,R,S)$-pseudosymmetric or in
short, Ricci-generalized pseudosymmetric if
\[
R\cdot R=L_{S}Q(S,R)
\]%
holds on the set ${\cal U}_{S}=\left\{ x\in M:Q(S,R)\not=0\right\} $, where $%
L_{S}$ is some smooth function on ${\cal U}_{S}$. It is known \cite%
{Deszcz-90} that every $3$-dimensional semi-Riemannian manifold is
Ricci-generalized pseudosymmetric along with $L_{S}=1$, that is, $R\cdot
R=Q(S,R)$.
\end{rem-new}

\begin{th}
\label{th-T-T-111} Let $M$ be an $n$-dimensional $({\cal T}_{\!a},{\cal T}%
_{\!b},S^{\ell })$-pseudosymmetric $(N(k),\xi )$-semi-Riemannian manifold.
Then
\begin{eqnarray}
&&\varepsilon b_{0}(ka_{0}+\varepsilon
k(n-1)a_{4}+a_{7}r)R(U,V,W,X)+\varepsilon a_{1}b_{0}S(X,R(U,V)W)  \nonumber
\\
&&-2k(n-1)a_{3}(kb_{0}+k(n-1)b_{4}+b_{7}r)\eta (X)\eta (U)g(V,W)  \nonumber
\\
&&-2k(n-1)a_{3}(-kb_{0}+k(n-1)b_{5}-b_{7}r)\eta (X)\eta (V)g(U,W)  \nonumber
\\
&&+\varepsilon a_{1}b_{4}S^{2}(X,U)g(V,W)+\varepsilon
a_{1}b_{5}S^{2}(X,V)g(U,W)  \nonumber \\
&&+\varepsilon a_{1}b_{6}S^{2}(X,W)g(U,V)-a_{5}(b_{1}+b_{3})S^{2}(X,V)\eta
(U)\eta (W)  \nonumber \\
&&-a_{5}(b_{1}+b_{2})S^{2}(X,W)\eta (U)\eta
(V)-a_{5}(b_{2}+b_{3})S^{2}(X,U)\eta (V)\eta (W)  \nonumber \\
&&-2a_{6}b_{1}S^{2}(V,W)\eta (X)\eta (U)-2a_{6}b_{2}S^{2}(U,W)\eta (X)\eta
(V)  \nonumber \\
&&-2a_{6}b_{3}S^{2}(U,V)\eta (X)\eta (W)-2k^{2}(n-1)a_{3}b_{6}g(U,V)\eta
(X)\eta (W)  \nonumber \\
&&-2\left( k(n-1)a_{3}b_{1}+a_{6}(kb_{0}+k(n-1)b_{4}+b_{7}r)\right) \eta
(X)\eta (U)S(V,W)  \nonumber \\
&&-2\left( k(n-1)a_{3}b_{2}+a_{6}(-kb_{0}+k(n-1)b_{5}-b_{7}r)\right) \eta
(X)\eta (V)S(U,W)  \nonumber \\
&&-2k(n-1)(a_{3}b_{3}+a_{6}b_{6})S(U,V)\eta (X)\eta (W)  \nonumber \\
&&+\varepsilon \left(
b_{4}(ka_{0}+k(n-1)a_{4}+a_{7}r)-a_{1}(kb_{0}+k(n-1)b_{4})\right)
S(X,U)g(V,W)  \nonumber \\
&&+\varepsilon \left(
b_{5}(ka_{0}+k(n-1)a_{4}+a_{7}r)-a_{1}(-kb_{0}+k(n-1)b_{5})\right)
S(X,V)g(U,W)  \nonumber \\
&&+\varepsilon b_{6}(ka_{0}+k(n-1)(a_{4}-a_{1})+a_{7}r)S(X,W)g(U,V)
\nonumber \\
&&-\varepsilon (kb_{0}+k(n-1)b_{4})(ka_{0}+k(n-1)a_{4}+a_{7}r)g(X,U)g(V,W)
\nonumber \\
&&-\varepsilon (-kb_{0}+k(n-1)b_{5})(ka_{0}+k(n-1)a_{4}+a_{7}r)g(U,W)g(X,V)
\nonumber \\
&&-\varepsilon k(n-1)b_{6}(ka_{0}+k(n-1)a_{4}+a_{7}r)g(X,W)g(U,V)  \nonumber
\\
&&-k(n-1)\left( (b_{2}+b_{3})(ka_{0}+k(n-1)a_{4}+a_{7}r)\right.  \nonumber \\
&&\left. +(a_{2}+a_{4})(-kb_{0}+k(n-1)(b_{5}+b_{6})-b_{7}r)\right)
g(X,U)\eta (V)\eta (W)  \nonumber \\
&&-k(n-1)\left( (b_{1}+b_{3})(ka_{0}+k(n-1)a_{4}+a_{7}r)\right.  \nonumber \\
&&\left. +(a_{2}+a_{4})(kb_{0}+k(n-1)(b_{4}+b_{6})+b_{7}r)\right) g(X,V)\eta
(U)\eta (W)  \nonumber \\
&&-\left( (b_{1}+b_{3})(-ka_{0}+k(n-1)(a_{1}+a_{2})-a_{7}r)\right.  \nonumber
\\
&&\left. +(a_{1}+a_{5})(kb_{0}+k(n-1)(b_{4}+b_{6})+b_{7}r)\right) S(X,V)\eta
(U)\eta (W)  \nonumber \\
&&-\left( (b_{2}+b_{3})(-ka_{0}+k(n-1)(a_{1}+a_{2})-a_{7}r)\right.  \nonumber
\\
&&\left. +(a_{1}+a_{5})(kb_{0}+k(n-1)(b_{5}+b_{6})+b_{7}r)\right) S(X,U)\eta
(V)\eta (W)  \nonumber \\
&&-\left( (b_{1}+b_{2})(-ka_{0}+k(n-1)(a_{1}+a_{2})-a_{7}r)\right.  \nonumber
\\
&&\left. +k(n-1)(b_{4}+b_{5})(a_{1}+a_{5})\right) S(X,W)\eta (U)\eta (V)
\nonumber \\
&&-k(n-1)\left( k(n-1)(b_{4}+b_{5})(a_{2}+a_{4})\right.  \nonumber \\
&&\left. +(b_{1}+b_{2})(ka_{0}+k(n-1)a_{4}+a_{7}r)\right) g(X,W)\eta (U)\eta
(V)  \nonumber \\
&=&L_{S^{\ell }}(\varepsilon b_{0}S^{\ell }(R(U,V)W,X)+\varepsilon
b_{4}S^{\ell +1}(X,U)g(V,W)+\varepsilon b_{5}S^{\ell +1}(X,V)g(U,W)
\nonumber \\
&&+\varepsilon b_{6}S^{\ell +1}(X,W)g(U,V)-\varepsilon k(n-1)b_{6}S^{\ell
}(X,W)g(U,V)  \nonumber \\
&&+k^{\ell }(n-1)^{\ell }(-kb_{0}+k(n-1)(b_{5}+b_{6})-b_{7}r)g(X,U)\eta
(V)\eta (W)  \nonumber \\
&&+k^{\ell }(n-1)^{\ell }(kb_{0}+k(n-1)(b_{4}+b_{6})+b_{7}r)g(X,V)\eta
(U)\eta (W)  \nonumber \\
&&-(-kb_{0}+k(n-1)(b_{2}+b_{3}+b_{5}+b_{6})-b_{7}r)S^{\ell }(X,U)\eta
(V)\eta (W)  \nonumber \\
&&-(kb_{0}+k(n-1)(b_{1}+b_{3}+b_{4}+b_{6})+b_{7}r)S^{\ell }(X,V)\eta (U)\eta
(W)  \nonumber \\
&&-k(n-1)(b_{1}+b_{2}+b_{4}+b_{5})S(X,W)\eta (U)\eta (V)  \nonumber \\
&&+k^{\ell +1}(n-1)^{\ell +1}(b_{4}+b_{5})g(X,W)\eta (U)\eta (V)  \nonumber
\\
&&-\varepsilon (-kb_{0}+k(n-1)b_{5})S^{\ell }(X,V)g(U,W)  \nonumber \\
&&-\varepsilon (kb_{0}+k(n-1)b_{4})S^{\ell }(X,U)g(V,W)  \nonumber \\
&&+k^{\ell }(n-1)^{\ell }(b_{1}+b_{3})S(X,V)\eta (U)\eta (W)  \nonumber \\
&&+k^{\ell }(n-1)^{\ell }(b_{1}+b_{2})S(X,W)\eta (U)\eta (V)  \nonumber \\
&&+k^{\ell }(n-1)^{\ell }(b_{2}+b_{3})S(X,U)\eta (V)\eta (W)).
\label{eq-T-T-2-i-ii}
\end{eqnarray}
\end{th}

\noindent {\bf Proof.} Let $M$ be an $n$-dimensional $({\cal T}_{\!a},{\cal T%
}_{\!b},S^{\ell })$-pseudosymmetric $(N(k),\xi )$-semi-Riemannian manifold.
Then
\begin{equation}
{\cal T}_{\!a}(Z,X)\cdot {\cal T}_{\!b}(U,V)W=L_{S^{\ell }}Q(S^{\ell },{\cal %
T}_{\!b})(U,V,W;Z,X).  \label{eq-txxi-1-ii}
\end{equation}%
Taking $Z=\xi $ in (\ref{eq-txxi-1-ii}), we get
\[
{\cal T}_{\!a}(\xi ,X)\cdot {\cal T}_{\!b}(U,V)W=L_{S^{\ell }}Q(S^{\ell },%
{\cal T}_{\!b})(U,V,W;\xi ,X),
\]%
which gives
\begin{eqnarray*}
&&{\cal T}_{\!a}(\xi ,X){\cal T}_{\!b}(U,V)W-{\cal T}_{\!b}({\cal T}%
_{\!a}(\xi ,X)U,V)W \\
&&\quad -\ {\cal T}_{\!b}(U,{\cal T}_{\!a}(\xi ,X)V)W-{\cal T}_{\!b}(U,V)%
{\cal T}_{\!a}(\xi ,X)W \\
&=&L_{S^{\ell }}((\xi \wedge _{S^{\ell }}X){\cal T}_{\!b}(U,V)W-{\cal T}%
_{\!b}((\xi \wedge _{S^{\ell }}X)U,V)W \\
&&-{\cal T}_{\!b}(U,(\xi \wedge _{S^{\ell }}X)V)W-{\cal T}_{\!b}(U,V)(\xi
\wedge _{S^{\ell }}X)W),
\end{eqnarray*}%
that is,
\begin{eqnarray}
&&{\cal T}_{\!a}(\xi ,X){\cal T}_{\!b}(U,V)W-{\cal T}_{\!b}({\cal T}%
_{\!a}(\xi ,X)U,V)W  \nonumber \\
&&\quad -\ {\cal T}_{\!b}(U,{\cal T}_{\!a}(\xi ,X)V)W-{\cal T}_{\!b}(U,V)%
{\cal T}_{\!a}(\xi ,X)W  \nonumber \\
&=&L_{S^{\ell }}(S^{\ell }(X,{\cal T}_{\!b}(U,V)W)\xi -S^{\ell }(\xi ,{\cal T%
}_{\!b}(U,V)W)X  \nonumber \\
&&-S^{\ell }(X,U){\cal T}_{\!b}(\xi ,V)W+S^{\ell }(\xi ,U){\cal T}%
_{\!b}(X,V)W  \nonumber \\
&&-S^{\ell }(X,V){\cal T}_{\!b}(U,\xi )W+S^{\ell }(\xi ,V){\cal T}%
_{\!b}(U,X)W  \nonumber \\
&&-S^{\ell }(X,W){\cal T}_{\!b}(U,V)\xi +S^{\ell }(\xi ,W){\cal T}%
_{\!b}(U,V)X).  \label{eq-T-T-i-ii}
\end{eqnarray}%
Taking the inner product of (\ref{eq-T-T-i-ii}) with $\xi $, we get
\begin{eqnarray}
&&{\cal T}_{\!a}(\xi ,X,{\cal T}_{\!b}(U,V)W,\xi )-{\cal T}_{\!b}({\cal T}%
_{\!a}(\xi ,X)U,V,W,\xi )  \nonumber \\
&&\quad -\ {\cal T}_{\!b}(U,{\cal T}_{\!a}(\xi ,X)V,W,\xi )-{\cal T}%
_{\!b}(U,V,{\cal T}_{\!a}(\xi ,X)W,\xi )  \nonumber \\
&=&L_{S^{\ell }}(\varepsilon S^{\ell }(X,{\cal T}_{\!b}(U,V)W)-\varepsilon
\eta (X)S^{\ell }(\xi ,{\cal T}_{\!b}(U,V)W)  \nonumber \\
&&-S^{\ell }(X,U){\cal T}_{\!b}(\xi ,V,W,\xi )+S^{\ell }(\xi ,U){\cal T}%
_{\!b}(X,V,W,\xi )  \nonumber \\
&&-S^{\ell }(X,V){\cal T}_{\!b}(U,\xi ,W,\xi )+S^{\ell }(\xi ,V){\cal T}%
_{\!b}(U,X,W,\xi )  \nonumber \\
&&-S^{\ell }(X,W){\cal T}_{\!b}(U,V,\xi ,\xi )+S^{\ell }(\xi ,W){\cal T}%
_{\!b}(U,V,X,\xi )).  \label{eq-T-T-1-i-ii}
\end{eqnarray}%
By using (\ref{eq-X-Y-xi}),\ldots ,(\ref{eq-X-xi-xi}) in (\ref{eq-T-T-1-i-ii}%
), we get (\ref{eq-T-T-2-i-ii}). $\blacksquare $

\begin{cor}
\label{th-T-T-11 copy(1)} Let $M$ be an $n$-dimensional $({\cal T}_{\!a},%
{\cal T}_{\!b},S)$-pseudosymmetric $(N(k),\xi )$-semi-Riemannian manifold.
Then
\begin{eqnarray}
&&\varepsilon b_{0}(ka_{0}+\varepsilon
k(n-1)a_{4}+a_{7}r)R(U,V,W,X)+\varepsilon a_{1}b_{0}S(X,R(U,V)W)  \nonumber
\\
&&-2k(n-1)a_{3}(kb_{0}+k(n-1)b_{4}+b_{7}r)\eta (X)\eta (U)g(V,W)  \nonumber
\\
&&-2k(n-1)a_{3}(-kb_{0}+k(n-1)b_{5}-b_{7}r)\eta (X)\eta (V)g(U,W)  \nonumber
\\
&&+\varepsilon a_{1}b_{4}S^{2}(X,U)g(V,W)+\varepsilon
a_{1}b_{5}S^{2}(X,V)g(U,W)  \nonumber \\
&&+\varepsilon a_{1}b_{6}S^{2}(X,W)g(U,V)-a_{5}(b_{1}+b_{3})S^{2}(X,V)\eta
(U)\eta (W)  \nonumber \\
&&-a_{5}(b_{1}+b_{2})S^{2}(X,W)\eta (U)\eta
(V)-a_{5}(b_{2}+b_{3})S^{2}(X,U)\eta (V)\eta (W)  \nonumber \\
&&-2a_{6}b_{1}S^{2}(V,W)\eta (X)\eta (U)-2a_{6}b_{2}S^{2}(U,W)\eta (X)\eta
(V)  \nonumber \\
&&-2a_{6}b_{3}S^{2}(U,V)\eta (X)\eta (W)-2k^{2}(n-1)a_{3}b_{6}g(U,V)\eta
(X)\eta (W)  \nonumber \\
&&-2\left( k(n-1)a_{3}b_{1}+a_{6}(kb_{0}+k(n-1)b_{4}+b_{7}r)\right) \eta
(X)\eta (U)S(V,W)  \nonumber \\
&&-2\left( k(n-1)a_{3}b_{2}+a_{6}(-kb_{0}+k(n-1)b_{5}-b_{7}r)\right) \eta
(X)\eta (V)S(U,W)  \nonumber \\
&&-2k(n-1)(a_{3}b_{3}+a_{6}b_{6})S(U,V)\eta (X)\eta (W)  \nonumber \\
&&+\varepsilon \left(
b_{4}(ka_{0}+k(n-1)a_{4}+a_{7}r)-a_{1}(kb_{0}+k(n-1)b_{4})\right)
S(X,U)g(V,W)  \nonumber \\
&&+\varepsilon \left(
b_{5}(ka_{0}+k(n-1)a_{4}+a_{7}r)-a_{1}(-kb_{0}+k(n-1)b_{5})\right)
S(X,V)g(U,W)  \nonumber \\
&&+\varepsilon b_{6}(ka_{0}+k(n-1)(a_{4}-a_{1})+a_{7}r)S(X,W)g(U,V)
\nonumber \\
&&-\varepsilon (kb_{0}+k(n-1)b_{4})(ka_{0}+k(n-1)a_{4}+a_{7}r)g(X,U)g(V,W)
\nonumber \\
&&-\varepsilon (-kb_{0}+k(n-1)b_{5})(ka_{0}+k(n-1)a_{4}+a_{7}r)g(U,W)g(X,V)
\nonumber \\
&&-\varepsilon k(n-1)b_{6}(ka_{0}+k(n-1)a_{4}+a_{7}r)g(X,W)g(U,V)  \nonumber
\\
&&-k(n-1)\left( (b_{2}+b_{3})(ka_{0}+k(n-1)a_{4}+a_{7}r)\right.  \nonumber \\
&&\left. +(a_{2}+a_{4})(-kb_{0}+k(n-1)(b_{5}+b_{6})-b_{7}r)\right)
g(X,U)\eta (V)\eta (W)  \nonumber \\
&&-k(n-1)\left( (b_{1}+b_{3})(ka_{0}+k(n-1)a_{4}+a_{7}r)\right.  \nonumber \\
&&\left. +(a_{2}+a_{4})(kb_{0}+k(n-1)(b_{4}+b_{6})+b_{7}r)\right) g(X,V)\eta
(U)\eta (W)  \nonumber \\
&&-\left( (b_{1}+b_{3})(-ka_{0}+k(n-1)(a_{1}+a_{2})-a_{7}r)\right.  \nonumber
\\
&&\left. +(a_{1}+a_{5})(kb_{0}+k(n-1)(b_{4}+b_{6})+b_{7}r)\right) S(X,V)\eta
(U)\eta (W)  \nonumber \\
&&-\left( (b_{2}+b_{3})(-ka_{0}+k(n-1)(a_{1}+a_{2})-a_{7}r)\right.  \nonumber
\\
&&\left. +(a_{1}+a_{5})(kb_{0}+k(n-1)(b_{5}+b_{6})+b_{7}r)\right) S(X,U)\eta
(V)\eta (W)  \nonumber \\
&&-\left( (b_{1}+b_{2})(-ka_{0}+k(n-1)(a_{1}+a_{2})-a_{7}r)\right.  \nonumber
\\
&&\left. +k(n-1)(b_{4}+b_{5})(a_{1}+a_{5})\right) S(X,W)\eta (U)\eta (V)
\nonumber \\
&&-k(n-1)\left( k(n-1)(b_{4}+b_{5})(a_{2}+a_{4})\right.  \nonumber \\
&&\left. +(b_{1}+b_{2})(ka_{0}+k(n-1)a_{4}+a_{7}r)\right) g(X,W)\eta (U)\eta
(V)  \nonumber \\
&=&\,L_{S}(\varepsilon b_{0}S(R(U,V)W,X)+\varepsilon
b_{4}S^{2}(X,U)g(V,W)+\varepsilon b_{5}S^{2}(X,V)g(U,W)  \nonumber \\
&&+\varepsilon b_{6}S^{2}(X,W)g(U,V)-\varepsilon k(n-1)b_{6}S(X,W)g(U,V)
\nonumber \\
&&+k(n-1)(-kb_{0}+k(n-1)(b_{5}+b_{6})-b_{7}r)g(X,U)\eta (V)\eta (W)
\nonumber \\
&&+k(n-1)(kb_{0}+k(n-1)(b_{4}+b_{6})b_{7}r)g(X,V)\eta (U)\eta (W)  \nonumber
\\
&&+k^{2}(n-1)^{2}(b_{4}+b_{5})g(X,W)\eta (U)\eta (V)  \nonumber \\
&&-(-kb_{0}+k(n-1)(b_{5}+b_{6})-b_{7}r)S(X,U)\eta (V)\eta (W)  \nonumber \\
&&-(kb_{0}+k(n-1)(b_{4}+b_{6})+b_{7}r)S(X,V)\eta (U)\eta (W)  \nonumber \\
&&-k(n-1)(b_{4}+b_{5})S(X,W)\eta (U)\eta (V)-\varepsilon
(-kb_{0}+k(n-1)b_{5})S(X,V)g(U,W)  \nonumber \\
&&-\varepsilon (kb_{0}+k(n-1)b_{4})S(X,U)g(V,W)).  \label{eq-T-T-2-i-i}
\end{eqnarray}
\end{cor}

\begin{th}
Let $M$ be an $n$-dimensional $({\cal T}_{\!a},{\cal T}_{\!b},S^{\ell })$%
-pseudosymmetric $(N(k),\xi )$-semi-Riemannian manifold. Then \vspace{-3cm}
\begin{eqnarray*}
&&(a_{1}b_{5}-a_{5}b_{1}-a_{5}b_{3})S^{2}(X,V)\eta (W) \\
&&+(a_{1}b_{6}-a_{5}b_{1}-a_{5}b_{2})S^{2}(X,W)\eta (V) \\
&&-2a_{6}b_{1}S^{2}(V,W)\eta (X) \\
&&+(b_{5}(ka_{0}+k(n-1)(a_{4}-a_{1})+a_{7}r) \\
&&-(b_{1}+b_{3})(-ka_{0}+k(n-1)(a_{1}+a_{2})-a_{7}r) \\
&&-(a_{1}+a_{5})(kb_{0}+k(n-1)(b_{4}+b_{6})+b_{7}r))S(X,V)\eta (W) \\
&&+(b_{6}(ka_{0}+k(n-1)(a_{4}-a_{1})+a_{7}r) \\
&&-(b_{1}+b_{2})(-ka_{0}+k(n-1)(a_{1}+a_{2})-a_{7}r) \\
&&-k(n-1)(a_{1}+a_{5})(b_{4}+b_{5}))S(X,W)\eta (V) \\
&&-2(k(n-1)a_{3}b_{1}+a_{6}(kb_{0}+k(n-1)b_{4}+b_{7}r))S(V,W)\eta (X) \\
&&-k(n-1)((b_{1}+b_{3}+b_{5})(ka_{0}+k(n-1)a_{4}+a_{7}r) \\
&&+(a_{2}+a_{4})(kb_{0}+k(n-1)(b_{4}+b_{6})+b_{7}r))g(X,V)\eta (W) \\
&&-k(n-1)((b_{1}+b_{2}+b_{6})(ka_{0}+k(n-1)a_{4}+a_{7}r) \\
&&+k(n-1)(a_{2}+a_{4})(b_{4}+b_{5}))g(X,W)\eta (V) \\
&&-2k(n-1)a_{3}(kb_{0}+k(n-1)b_{4}+b_{7}r)g(V,W)\eta (X) \\
&&-\varepsilon k(n-1)((kb_{0}+b_{7}r)(a_{1}+a_{2}+a_{4}+a_{5}) \\
&&+k(n-1)(b_{2}+b_{3}+b_{5}+b_{6})\times \\
&&(a_{1}+a_{2}+2a_{3}+a_{4}+a_{5}+2a_{6}))\eta (X)\eta (V)\eta (W) \\
&&=\, L_{S^{\ell
}}(-(kb_{0}+k(n-1)(b_{1}+b_{3}+b_{4}+b_{5}+b_{6})+b_{7}r)S^{\ell }(X,V)\eta
(W) \\
&&+k^{\ell }(n-1)^{\ell }(kb_{0}+k(n-1)(b_{4}+b_{6})+b_{7}r)g(X,V)\eta (W) \\
&&+k^{\ell +1}(n-1)^{\ell +1}(b_{4}+b_{5})g(X,W)\eta (V) \\
&&-k(n-1)(b_{1}+b_{2}+b_{4}+b_{5}+b_{6})S^{\ell }(X,W)\eta (V) \\
&&+b_{5}S^{\ell +1}(X,V)\eta (W)+b_{6}S^{\ell +1}(X,W)\eta (V)) \\
&&+k^{\ell }(n-1)^{\ell }(b_{1}+b_{3})S(X,V)\eta (W) \\
&&+k^{\ell }(n-1)^{\ell }(b_{1}+b_{2})S(X,W)\eta (V)).
\end{eqnarray*}
\end{th}

\begin{cor}
Let $M$ be an $n$-dimensional $({\cal T}_{\!a},{\cal T}_{\!b},S)$%
-pseudosymmetric $(N(k),\xi )$-semi-Riemannian manifold. Then
\begin{eqnarray*}
&&(a_{1}b_{5}-a_{5}b_{1}-a_{5}b_{3})S^{2}(X,V)\eta (W) \\
&&+(a_{1}b_{6}-a_{5}b_{1}-a_{5}b_{2})S^{2}(X,W)\eta (V) -
2a_{6}b_{1}S^{2}(V,W)\eta (X) \\
&&+(b_{5}(ka_{0}+k(n-1)(a_{4}-a_{1})+a_{7}r) \\
&&-(b_{1}+b_{3})(-ka_{0}+k(n-1)(a_{1}+a_{2})-a_{7}r) \\
&&-(a_{1}+a_{5})(kb_{0}+k(n-1)(b_{4}+b_{6})+b_{7}r))S(X,V)\eta (W) \\
&&+(b_{6}(ka_{0}+k(n-1)(a_{4}-a_{1})+a_{7}r) \\
&&-(b_{1}+b_{2})(-ka_{0}+k(n-1)(a_{1}+a_{2})-a_{7}r) \\
&&-k(n-1)(a_{1}+a_{5})(b_{4}+b_{5}))S(X,W)\eta (V) \\
&&-2(k(n-1)a_{3}b_{1}+a_{6}(kb_{0}+k(n-1)b_{4}+b_{7}r))S(V,W)\eta (X) \\
&&-k(n-1)((b_{1}+b_{3}+b_{5})(ka_{0}+k(n-1)a_{4}+a_{7}r) \\
&&+(a_{2}+a_{4})(kb_{0}+k(n-1)(b_{4}+b_{6})+b_{7}r))g(X,V)\eta (W) \\
&&-k(n-1)((b_{1}+b_{2}+b_{6})(ka_{0}+k(n-1)a_{4}+a_{7}r) \\
&&+k(n-1)(a_{2}+a_{4})(b_{4}+b_{5}))g(X,W)\eta (V) \\
&&-2k(n-1)a_{3}(kb_{0}+k(n-1)b_{4}+b_{7}r)g(V,W)\eta (X) \\
&&-\varepsilon k(n-1)((kb_{0}+b_{7}r)(a_{1}+a_{2}+a_{4}+a_{5}) \\
&&+k(n-1)(b_{2}+b_{3}+b_{5}+b_{6})\times \\
&&(a_{1}+a_{2}+2a_{3}+a_{4}+a_{5}+2a_{6}))\eta (X)\eta (V)\eta (W)
\end{eqnarray*}
\begin{eqnarray}
&=&L_{S}(-(kb_{0}+k(n-1)(b_{4}+b_{5}+b_{6})+b_{7}r)S(X,V)\eta (W)  \nonumber
\\
&&+k(n-1)(kb_{0}+k(n-1)(b_{4}+b_{6})+b_{7}r)g(X,V)\eta (W)  \nonumber \\
&&+k^{2}(n-1)^{2}(b_{4}+b_{5})g(X,W)\eta (V)  \nonumber \\
&&-k(n-1)(b_{4}+b_{5}+b_{6})S(X,W)\eta (V)  \nonumber \\
&&+b_{5}S^{2}(X,V)\eta (W)+b_{6}S^{2}(X,W)\eta (V)).  \label{eq-T-T-2-i-i-1}
\end{eqnarray}
\end{cor}

\begin{th}
Let $M$ be an $n$-dimensional $(R,{\cal T}_{\!a},S^{\ell })$-pseudosymmetric
$(N(k),\xi )$-semi-Riemannian manifold. Then
\begin{eqnarray}
&&k(a_{1}+a_{3}+a_{5})\left( S(X,V)-\,k(n-1)g(X,V)\right) \eta (W)  \nonumber
\\
&+&k(a_{1}+a_{2}+a_{6})\left( S(X,W)-\,k(n-1)g(X,W)\right) \eta (V)
\nonumber \\
&=&L_{S^{\ell
}}(-(ka_{0}+k(n-1)(a_{1}+a_{3}+a_{4}+a_{5}+a_{6})+a_{7}r)S^{\ell }(X,V)\eta
(W)  \nonumber \\
&&+k^{\ell }(n-1)^{\ell }(ka_{0}+k(n-1)(a_{4}+a_{6})+a_{7}r)g(X,V)\eta (W)
\nonumber \\
&&+k^{\ell +1}(n-1)^{\ell +1}(a_{4}+a_{5})g(X,W)\eta (V)  \nonumber \\
&&-k(n-1)(a_{1}+a_{2}+a_{4}+a_{5}+a_{6})S^{\ell }(X,W)\eta (V)  \nonumber \\
&&+a_{5}S^{\ell +1}(X,V)\eta (W)+a_{6}S^{\ell +1}(X,W)\eta (V)  \nonumber \\
&&+k^{\ell }(n-1)^{\ell }(a_{1}+a_{3})S(X,V)\eta (W)  \nonumber \\
&&+k^{\ell }(n-1)^{\ell }(a_{1}+a_{2})S(X,W)\eta (V)).  \label{eq-TRGP-i}
\end{eqnarray}
\end{th}

\begin{cor}
Let $M$ be an $n$-dimensional $(R,{\cal T}_{\!a},S)$-pseudosymmetric $%
(N(k),\xi )$-semi-Riemannian manifold. Then
\begin{eqnarray}
&&k(a_{1}+a_{3}+a_{5})\left( S(X,V)-\,k(n-1)g(X,V)\right) \eta (W)  \nonumber
\\
&+&k(a_{1}+a_{2}+a_{6})\left( S(X,W)-\,k(n-1)g(X,W)\right) \eta (V)
\nonumber \\
&=&L_{S}(-(ka_{0}+k(n-1)(a_{4}+a_{5}+a_{6})+a_{7}r)S(X,V)\eta (W)  \nonumber
\\
&&+k(n-1)(ka_{0}+k(n-1)(a_{4}+a_{6})+a_{7}r)g(X,V)\eta (W)  \nonumber \\
&&+k^{2}(n-1)^{2}(a_{4}+a_{5})g(X,W)\eta (V)  \nonumber \\
&&-k(n-1)(a_{4}+a_{5}+a_{6})S(X,W)\eta (V)  \nonumber \\
&&+a_{5}S^{2}(X,V)\eta (W)+a_{6}S^{2}(X,W)\eta (V)).  \label{eq-TRGP}
\end{eqnarray}
\end{cor}

\begin{cor}
Let $M$ be an $n$-dimensional $(R,R,S^{\ell })$-pseudosymmetric $\left(
N(k),\xi \right) $-semi-Riemannian manifold. If $M$ is not semisymmetric,
then
\[
S^{\ell }=k^{\ell }(n-1)^{\ell }g
\]%
and $L_{S^{\ell }}=\dfrac{1}{k^{\ell -1}(n-1)^{\ell }}$. Consequently, we
have the following\/{\rm :}~
\[
\begin{tabular}{|l|l|l|}
\hline
${\boldmath M}$ & $L_{S^{\ell }}=$ & ${\boldmath S}^{\ell }{=}$ \\ \hline
$N(k)$-contact metric & $\dfrac{1}{k^{\ell -1}(n-1)^{\ell }}$ & $k^{\ell
}(n-1)^{\ell }g$ \\ \hline
Sasakian & $\dfrac{1}{(n-1)^{\ell }}$ & $(n-1)^{\ell }g$ \\ \hline
Kenmotsu & $\dfrac{1}{(-1)^{\ell -1}(n-1)^{\ell }}$ & $(-1)^{\ell
}(n-1)^{\ell }g$ \\ \hline
$(\varepsilon )$-Sasakian & $\dfrac{1}{(\varepsilon )^{\ell -1}(n-1)^{\ell }}
$ & $(\varepsilon )^{\ell }(n-1)^{\ell }g$ \\ \hline
para-Sasakian & $\dfrac{1}{(-1)^{\ell -1}(n-1)^{\ell }}$ & $(-1)^{\ell
}(n-1)^{\ell }g$ \\ \hline
$(\varepsilon )$-para-Sasakian & $\dfrac{1}{(-\varepsilon )^{\ell
-1}(n-1)^{\ell }}$ & $(-\varepsilon )^{\ell }(n-1)^{\ell }g$ \\ \hline
\end{tabular}%
\]
\end{cor}

\noindent {\bf Proof.} Let $M$ be an $n$-dimensional $(R,R,S^{\ell })$%
-pseudosymmetric $\left( N(k),\xi \right) $-semi-Riemannian manifold, that
is
\begin{equation}
R\cdot R=L_{S^{\ell }}Q(S^{\ell },R)  \label{eq-TRGP-1-i}
\end{equation}%
holds on $M$. By putting the value for $R$ in (\ref{eq-TRGP-i}), we get
\begin{equation}
-kL_{S^{\ell }}\left( S^{\ell }(X,V)-k^{\ell }(n-1)^{\ell }g(X,V)\right)
\eta (W)=0.  \label{eq-TRGP-2-i}
\end{equation}%
Putting $W=\xi $ in (\ref{eq-TRGP-2-i}), we get
\begin{equation}
-kL_{S^{\ell }}\left( S^{\ell }(X,V)-k^{\ell }(n-1)^{\ell }g(X,V)\right) =0.
\label{eq-TRGP-3-i}
\end{equation}%
Since $M$ is not semisymmetric $L_{S^{\ell }}\not=0$. Therefore from (\ref%
{eq-TRGP-3-i}), we have
\[
S^{\ell }(X,V)=k^{\ell }(n-1)^{\ell }g(X,V).
\]%
So putting $S^{\ell }=k^{\ell }(n-1)^{\ell }g$ in (\ref{eq-TRGP-1-i}), we
get
\[
R\cdot R=k^{\ell }(n-1)^{\ell }L_{S^{\ell }}Q(g,R),
\]%
which is the condition of pseudosymmetric manifold. By comparison with the
result of pseudosymmetric manifold (Corollary \ref{cor-pseudo-1}), we get $%
L_{S^{\ell }}=\dfrac{1}{k^{\ell -1}(n-1)^{\ell }}$. This proves the result. $%
\blacksquare $

\begin{cor}
Let $M$ be an $n$-dimensional Ricci-generalized pseudosymmetric $\left(
N(k),\xi \right) $-semi-Riemannian manifold. If $M$ is not semisymmetric,
then $M$ is an Einstein manifold with scalar curvature $kn(n-1)$ and $L_{S}=%
\dfrac{1}{n-1}$. Consequently, we have the following\/{\rm :}~
\[
\begin{tabular}{|l|l|l|}
\hline
${\boldmath M}$ & $L_{S}=$ & ${\boldmath S=}$ \\ \hline
$N(k)$-contact metric & $\dfrac{1}{n-1}$ & $k(n-1)g$ \\ \hline
Sasakian & $\dfrac{1}{n-1}$ & $(n-1)g$ \\ \hline
Kenmotsu {\rm \cite{Ozgur-06}} & $\dfrac{1}{n-1}$ & $-(n-1)g$ \\ \hline
$(\varepsilon )$-Sasakian & $\dfrac{1}{n-1}$ & $\varepsilon (n-1)g$ \\ \hline
para-Sasakian & $\dfrac{1}{n-1}$ & $-(n-1)g$ \\ \hline
$(\varepsilon )$-para-Sasakian & $\dfrac{1}{n-1}$ & $-\,\varepsilon (n-1)g$
\\ \hline
\end{tabular}%
\
\]
\end{cor}


\begin{cor}
Let $M$ be an $n$-dimensional $(R,{\cal C}_{\ast },S^{\ell })$%
-pseudosymmetric $\left( N(k),\xi \right) $-semi-Riemannian manifold. If $M$
is not quasi-conformal semisymmetric, then
\begin{eqnarray*}
S^{\ell +1} &=&\left( \left( \frac{r}{n(n-1)}-k\right) \frac{a_{0}}{a_{1}}%
+\left( \frac{2r}{n}-k(n-1)\right) \right) \left( S^{\ell }-k^{\ell
}(n-1)^{\ell }g\right) \\
&&+k^{\ell }(n-1)^{\ell }S.
\end{eqnarray*}%
Consequently, we have the following\/{\rm :}~
\[
\begin{tabular}{|l|l|}
\hline
${\boldmath M}$ & ${\boldmath S^{\ell +1}=}$ \\ \hline
$N(k)$-contact metric & $\left( \left( \dfrac{r}{n(n-1)}-k\right) \dfrac{%
a_{0}}{a_{1}}+\left( \dfrac{2r}{n}-k(n-1)\right) \right) \left( S^{\ell
}-k^{\ell }(n-1)^{\ell }g\right) $ \\
& $+k^{\ell }(n-1)^{\ell }S$ \\ \hline
Sasakian & $\left( \left( \dfrac{r}{n(n-1)}-1\right) \dfrac{a_{0}}{a_{1}}%
+\left( \dfrac{2r}{n}-1(n-1)\right) \right) \left( S^{\ell }-(n-1)^{\ell
}g\right) $ \\
& $+(n-1)^{\ell }S$ \\ \hline
Kenmotsu & $\left( \left( \dfrac{r}{n(n-1)}+1\right) \dfrac{a_{0}}{a_{1}}%
+\left( \frac{2r}{n}+(n-1)\right) \right) \left( S^{\ell }-(-1)^{\ell
}(n-1)^{\ell }g\right) $ \\
& $+(-1)^{\ell }(n-1)^{\ell }S$ \\ \hline
$(\varepsilon )$-Sasakian & $\left( \left( \dfrac{r}{n(n-1)}-\varepsilon
\right) \dfrac{a_{0}}{a_{1}}+\left( \frac{2r}{n}-\varepsilon (n-1)\right)
\right) \left( S^{\ell }-(\varepsilon )^{\ell }(n-1)^{\ell }g\right) $ \\
& $+(\varepsilon )^{\ell }(n-1)^{\ell }S$ \\ \hline
para-Sasakian & $\left( \left( \dfrac{r}{n(n-1)}+1\right) \dfrac{a_{0}}{a_{1}%
}+\left( \frac{2r}{n}+(n-1)\right) \right) \left( S^{\ell }-(-1)^{\ell
}(n-1)^{\ell }g\right) $ \\
& $+(-1)^{\ell }(n-1)^{\ell }S$ \\ \hline
$(\varepsilon )$-para-Sasakian & $\left( \left( \dfrac{r}{n(n-1)}%
+\varepsilon \right) \dfrac{a_{0}}{a_{1}}+\left( \frac{2r}{n}+\varepsilon
(n-1)\right) \right) \left( S^{\ell }-(-\varepsilon )^{\ell }(n-1)^{\ell
}g\right) $ \\
& $+(-\varepsilon )^{\ell }(n-1)^{\ell }S$ \\ \hline
\end{tabular}
\]
\end{cor}

\begin{cor}
Let $M$ be an $n$-dimensional $(R,{\cal C}_{\ast },S)$-pseudosymmetric $%
\left( N(k),\xi \right) $-semi-Riemannian manifold. If $M$ is not
quasi-conformal-semisymmetric, then
\begin{eqnarray*}
S^{2} &=&\left( \left( \frac{r}{n(n-1)}-k\right) \frac{a_{0}}{a_{1}}+\frac{2r%
}{n}\right) S \\
&&-k(n-1)\left( \left( \frac{r}{n(n-1)}-k\right) \frac{a_{0}}{a_{1}}+\left(
\frac{2r}{n}-k(n-1)\right) \right) g.
\end{eqnarray*}%
Consequently, we have the following\/{\rm :}~
\[
\begin{tabular}{|l|l|}
\hline
${\boldmath M}$ & ${\boldmath S}^{2}{=}$ \\ \hline
$N(k)$-contact metric & $\left( \left( \dfrac{r}{n(n-1)}-k\right) \dfrac{%
a_{0}}{a_{1}}+\dfrac{2r}{n}\right) S$ \\
& $-k(n-1)\left( \left( \dfrac{r}{n(n-1)}-k\right) \dfrac{a_{0}}{a_{1}}%
+\left( \dfrac{2r}{n}-k(n-1)\right) \right) g$ \\ \hline
Sasakian & $\left( \left( \dfrac{r}{n(n-1)}-1\right) \dfrac{a_{0}}{a_{1}}+%
\dfrac{2r}{n}\right) S$ \\
& $-(n-1)\left( \left( \dfrac{r}{n(n-1)}-1\right) \dfrac{a_{0}}{a_{1}}%
+\left( \dfrac{2r}{n}-(n-1)\right) \right) g$ \\ \hline
Kenmotsu & $\left( \left( \dfrac{r}{n(n-1)}+1\right) \dfrac{a_{0}}{a_{1}}+%
\dfrac{2r}{n}\right) S$ \\
& $+(n-1)\left( \left( \dfrac{r}{n(n-1)}+1\right) \dfrac{a_{0}}{a_{1}}%
+\left( \dfrac{2r}{n}+(n-1)\right) \right) g$ \\ \hline
$(\varepsilon )$-Sasakian & $\left( \left( \dfrac{r}{n(n-1)}-\varepsilon
\right) \dfrac{a_{0}}{a_{1}}+\dfrac{2r}{n}\right) S$ \\
& $-\varepsilon (n-1)\left( \left( \dfrac{r}{n(n-1)}-\varepsilon \right)
\dfrac{a_{0}}{a_{1}}+\left( \dfrac{2r}{n}-\varepsilon (n-1)\right) \right) g$
\\ \hline
para-Sasakian & $\left( \left( \dfrac{r}{n(n-1)}+1\right) \dfrac{a_{0}}{a_{1}%
}+\dfrac{2r}{n}\right) S$ \\
& $+(n-1)\left( \left( \dfrac{r}{n(n-1)}+1\right) \dfrac{a_{0}}{a_{1}}%
+\left( \dfrac{2r}{n}+(n-1)\right) \right) g$ \\ \hline
$(\varepsilon )$-para-Sasakian & $\left( \left( \dfrac{r}{n(n-1)}%
+\varepsilon \right) \dfrac{a_{0}}{a_{1}}+\dfrac{2r}{n}\right) S$ \\
& $+\varepsilon (n-1)\left( \left( \dfrac{r}{n(n-1)}+\varepsilon \right)
\dfrac{a_{0}}{a_{1}}+\left( \dfrac{2r}{n}+\varepsilon (n-1)\right) \right) g$
\\ \hline
\end{tabular}%
\]
\end{cor}

\begin{cor}
Let $M$ be an $n$-dimensional $(R,{\cal C},S^{\ell })$-pseudosymmetric $%
\left( N(k),\xi \right) $-semi-Riemannian manifold. If $M$ is not
Weyl-semisymmetric, then
\begin{eqnarray*}
S^{\ell +1} &=&\left( \frac{r}{n-1}-k\right) \left( S^{\ell }-k^{\ell
}(n-1)^{\ell }g\right) +k^{\ell }(n-1)^{\ell }S.
\end{eqnarray*}%
Consequently, we have the following\/{\rm :}~
\[
\begin{tabular}{|l|l|}
\hline
${\boldmath M}$ & ${\boldmath S^{\ell +1}=}$ \\ \hline
$N(k)$-contact metric & $\left( \dfrac{r}{n-1}-k\right) \left( S^{\ell
}-k^{\ell }(n-1)^{\ell }g\right) $ \\
& $+k^{\ell }(n-1)^{\ell }S$ \\ \hline
Sasakian & $\left( \dfrac{r}{n-1}-1\right) \left( S^{\ell }-(n-1)^{\ell
}g\right) $ \\
& $+(n-1)^{\ell }S$ \\ \hline
Kenmotsu & $\left( \dfrac{r}{n-1}+1\right) \left( S^{\ell }-(-1)^{\ell
}(n-1)^{\ell }g\right) $ \\
& $+(-1)^{\ell }(n-1)^{\ell }S$ \\ \hline
$(\varepsilon )$-Sasakian & $\left( \dfrac{r}{n-1}-\varepsilon \right)
\left( S^{\ell }-(\varepsilon )^{\ell }(n-1)^{\ell }g\right) $ \\
& $+(\varepsilon )^{\ell }(n-1)^{\ell }S$ \\ \hline
para-Sasakian & $\left( \dfrac{r}{n-1}+1\right) \left( S^{\ell }-(-1)^{\ell
}(n-1)^{\ell }g\right) $ \\
& $+(-1)^{\ell }(n-1)^{\ell }S$ \\ \hline
$(\varepsilon )$-para-Sasakian & $\left( \dfrac{r}{n-1}+\varepsilon \right)
\left( S^{\ell }-(-\varepsilon )^{\ell }(n-1)^{\ell }g\right) $ \\
& $+(-\varepsilon )^{\ell }(n-1)^{\ell }S$ \\ \hline
\end{tabular}%
\
\]
\end{cor}

\begin{cor}
Let $M$ be an $n$-dimensional $(R,{\cal C},S)$-pseudosymmetric $\left(
N(k),\xi \right) $-semi-Riemannian manifold. If $M$ is not
Weyl-semisymmetric, then
\begin{eqnarray*}
S^{2} &=&\left( k(n-2)+\frac{r}{n-1}\right) S +k(n-1)\left( k-\frac{r}{n-1}%
\right) g.
\end{eqnarray*}%
Consequently, we have the following\/{\rm :}~
\[
\begin{tabular}{|l|l|}
\hline
${\boldmath M}$ & ${\boldmath S}^{2}{=}$ \\ \hline
$N(k)$-contact metric & $\left( k(n-2)+\frac{r}{n-1}\right) S$ \\
& $+k(n-1)\left( k-\frac{r}{n-1}\right) g$ \\ \hline
Sasakian & $\left( (n-2)+\frac{r}{n-1}\right) S$ \\
& $+(n-1)\left( 1-\frac{r}{n-1}\right) g$ \\ \hline
Kenmotsu & $\left( -(n-2)+\frac{r}{n-1}\right) S$ \\
& $+(n-1)\left( 1+\frac{r}{n-1}\right) g$ \\ \hline
$(\varepsilon )$-Sasakian & $\left( \varepsilon (n-2)+\frac{r}{n-1}\right) S$
\\
& $+\varepsilon (n-1)\left( \varepsilon -\frac{r}{n-1}\right) g$ \\ \hline
para-Sasakian & $\left( -(n-2)+\frac{r}{n-1}\right) S$ \\
& $+(n-1)\left( 1+\frac{r}{n-1}\right) g$ \\ \hline
$(\varepsilon )$-para-Sasakian & $\left( -\varepsilon (n-2)+\frac{r}{n-1}%
\right) S$ \\
& $+\varepsilon (n-1)\left( \varepsilon +\frac{r}{n-1}\right) g$ \\ \hline
\end{tabular}%
\
\]
\end{cor}

\begin{cor}
Let $M$ be an $n$-dimensional $(R,{\cal L},S^{\ell })$-pseudosymmetric $%
\left( N(k),\xi \right) $-semi-Riemannian manifold. If $M$ is not
conharmonic semisymmetric, then
\[
S^{\ell +1}=-kS^{\ell }+k^{\ell }(n-1)^{\ell }S+k^{\ell +1}(n-1)^{\ell }g.
\]%
Consequently, we have the following\/{\rm :}~
\[
\begin{tabular}{|l|l|}
\hline
${\boldmath M}$ & ${\boldmath S^{\ell +1}=}$ \\ \hline
$N(k)$-contact metric & $-kS^{\ell }+k^{\ell }(n-1)^{\ell }S+k^{\ell
+1}(n-1)^{\ell }g$ \\ \hline
Sasakian & $-S^{\ell }+(n-1)^{\ell }S+(n-1)^{\ell }g$ \\ \hline
Kenmotsu & $S^{\ell }+(-1)^{\ell }(n-1)^{\ell }S+(-1)^{\ell +1}(n-1)^{\ell
}g $ \\ \hline
$(\varepsilon )$-Sasakian & $-\varepsilon S^{\ell }+(\varepsilon )^{\ell
}(n-1)^{\ell }S+(\varepsilon )^{\ell +1}(n-1)^{\ell }g$ \\ \hline
para-Sasakian & $S^{\ell }+(-1)^{\ell }(n-1)^{\ell }S+(-1)^{\ell
+1}(n-1)^{\ell }g$ \\ \hline
$(\varepsilon )$-para-Sasakian & $\varepsilon S^{\ell }+(-\varepsilon
)^{\ell }(n-1)^{\ell }S+(-\varepsilon )^{\ell +1}(n-1)^{\ell }g$ \\ \hline
\end{tabular}%
\
\]
\end{cor}

\begin{cor}
Let $M$ be an $n$-dimensional $(R,{\cal L},S)$-pseudosymmetric $\left(
N(k),\xi \right) $-semi-Riemannian manifold. If $M$ is not conharmonic
semisymmetric, then
\[
S^{2}=k(n-2)S+k^{2}(n-1)g.
\]%
Consequently, we have the following\/{\rm :}~
\[
\begin{tabular}{|l|l|}
\hline
${\boldmath M}$ & ${\boldmath S}^{2}{=}$ \\ \hline
$N(k)$-contact metric & $k(n-2)S+k^{2}(n-1)g$ \\ \hline
Sasakian & $(n-2)S+(n-1)g$ \\ \hline
Kenmotsu & $-(n-2)S+(n-1)g$ \\ \hline
$(\varepsilon )$-Sasakian & $\varepsilon (n-2)S+(n-1)g$ \\ \hline
para-Sasakian & $-(n-2)S+(n-1)g$ \\ \hline
$(\varepsilon )$-para-Sasakian & $-\varepsilon (n-2)S+(n-1)g$ \\ \hline
\end{tabular}%
\
\]
\end{cor}

\begin{cor}
Let $M$ be an $n$-dimensional $(R,{\cal V},S^{\ell })$-pseudosymmetric $%
\left( N(k),\xi \right) $-semi-Riemannian manifold. If $M$ is not
concircular semisymmetric, then $M$ either satisfies
\[
S^{\ell }=k^{\ell }(n-1)^{\ell }g
\]%
or scalar curvature is $kn(n-1)$ and $L_{S^{\ell }}=\dfrac{1}{k^{\ell
-1}(n-1)^{\ell }}$. Consequently, we have the following\/{\rm :}~
\[
\begin{tabular}{|l|l|l|}
\hline
${\boldmath M}$ & $L_{S^{\ell }}=$ & {\bf Result} \\ \hline
$N(k)$-contact metric & $\dfrac{1}{k^{\ell -1}(n-1)^{\ell }}$ & $S^{\ell
}=k^{\ell }(n-1)^{\ell }g\quad {\rm or\quad }r=kn(n-1)$ \\ \hline
Sasakian & $\dfrac{1}{(n-1)^{\ell }}$ & $S^{\ell }=(n-1)^{\ell }g\quad {\rm %
or\quad }r=n(n-1)$ \\ \hline
Kenmotsu & $\dfrac{1}{(-1)^{\ell -1}(n-1)^{\ell }}$ & $S^{\ell }=(-1)^{\ell
}(n-1)^{\ell }g\quad {\rm or\quad }r=-n(n-1)$ \\ \hline
$(\varepsilon )$-Sasakian & $\dfrac{1}{(\varepsilon )^{\ell -1}(n-1)^{\ell }}
$ & $S^{\ell }=(\varepsilon )^{\ell }(n-1)^{\ell }g\quad {\rm or\quad }%
r=\varepsilon n(n-1)$ \\ \hline
para-Sasakian & $\dfrac{1}{(-1)^{\ell -1}(n-1)^{\ell }}$ & $S^{\ell
}=(-1)^{\ell }(n-1)^{\ell }g\quad {\rm or\quad }r=-n(n-1)$ \\ \hline
$(\varepsilon )$-para-Sasakian & $\dfrac{1}{(-\varepsilon )^{\ell
-1}(n-1)^{\ell }}$ & $S^{\ell }=(-\varepsilon )^{\ell }(n-1)^{\ell }g\quad
{\rm or\quad }r=-\varepsilon n(n-1)$ \\ \hline
\end{tabular}%
\]
\end{cor}

\begin{cor}
Let $M$ be an $n$-dimensional $(R,{\cal V},S)$-pseudosymmetric $\left(
N(k),\xi \right) $-semi-Riemannian manifold. If $M$ is not concircularly
semisymmetric, then $M$ is either an Einstein manifold or scalar curvature
is $kn(n-1)$ and $L_{S}=\dfrac{1}{n-1}$. Consequently, we have the
following\/{\rm :}~
\[
\begin{tabular}{|l|l|l|}
\hline
${\boldmath M}$ & $L_{S}=$ & ${\boldmath S=}$ \\ \hline
$N(k)$-contact metric & $\frac{1}{n-1}$ & $S=k(n-1)g\quad {\rm or\quad }%
r=kn(n-1)$ \\ \hline
Sasakian & $\frac{1}{n-1}$ & $S=(n-1)g\quad {\rm or\quad }r=n(n-1)$ \\ \hline
Kenmotsu {\rm \cite{Ozgur-06}} & $\frac{1}{n-1}$ & $S=-(n-1)g\quad {\rm %
or\quad }r=-n(n-1)$ \\ \hline
$(\varepsilon )$-Sasakian & $\frac{1}{n-1}$ & $S=\varepsilon (n-1)g\quad
{\rm or\quad }r=\varepsilon n(n-1)$ \\ \hline
para-Sasakian & $\frac{1}{n-1}$ & $S=-(n-1)g\quad {\rm or\quad }r=-n(n-1)$
\\ \hline
$(\varepsilon )$-para-Sasakian & $\frac{1}{n-1}$ & $S=-\,\varepsilon
(n-1)g\quad {\rm or\quad }r=-\varepsilon n(n-1)$ \\ \hline
\end{tabular}%
\]
\end{cor}

\begin{cor}
Let $M$ be an $n$-dimensional $(R,{\cal P}_{\ast },S^{\ell })$%
-pseudosymmetric $\left( N(k),\xi \right) $-semi-Riemannian manifold such
that $a_{0}+(n-1)a_{1}\not=0$. If $M$ is not pseudo-projective
semisymmetric, then
\begin{eqnarray*}
&&\left( \left( k-\frac{r}{n(n-1)}\right) a_{0}+\left( k(n-1)-\frac{r}{n}%
\right) a_{1}\right) S^{\ell } \\
&=&k^{\ell }(n-1)^{\ell }\left( \left( k-\frac{r}{n(n-1)}\right)
a_{0}-\left( \frac{r}{n}\right) a_{1}\right) g \\
&&+k^{\ell }(n-1)^{\ell }a_{1}S.
\end{eqnarray*}%
Consequently, we have the following\/{\rm :}~
\[
\begin{tabular}{|l|l|}
\hline
${\boldmath M}$ & ${\boldmath S^{\ell }=}$ \\ \hline
$N(k)$-contact metric & $\dfrac{k^{\ell }(n-1)^{\ell }a_{1}}{{\left( \left(
k-\dfrac{r}{n(n-1)}\right) a_{0}+\left( k(n-1)-\dfrac{r}{n}\right)
a_{1}\right) }}S$ \\
& $+\dfrac{k^{\ell }(n-1)^{\ell }\left( \left( k-\dfrac{r}{n(n-1)}\right)
a_{0}-\left( \dfrac{r}{n}\right) a_{1}\right) }{{\left( \left( k-\dfrac{r}{%
n(n-1)}\right) a_{0}+\left( k(n-1)-\dfrac{r}{n}\right) a_{1}\right) }}g$ \\
\hline
Sasakian & $\dfrac{(n-1)^{\ell }a_{1}}{{\left( \left( 1-\dfrac{r}{n(n-1)}%
\right) a_{0}+\left( (n-1)-\dfrac{r}{n}\right) a_{1}\right) }}S$ \\
& $+\dfrac{(n-1)^{\ell }\left( \left( 1-\dfrac{r}{n(n-1)}\right)
a_{0}-\left( \frac{r}{n}\right) a_{1}\right) }{{\left( \left( 1-\dfrac{r}{%
n(n-1)}\right) a_{0}+\left( (n-1)-\dfrac{r}{n}\right) a_{1}\right) }}g$ \\
\hline
Kenmotsu & $\dfrac{(-1)^{\ell }(n-1)^{\ell }a_{1}}{{\left( \left( -1-\dfrac{r%
}{n(n-1)}\right) a_{0}+\left( -(n-1)-\dfrac{r}{n}\right) a_{1}\right) }}S$
\\
& $+\dfrac{(-1)^{\ell }(n-1)^{\ell }\left( \left( -1-\frac{r}{n(n-1)}\right)
a_{0}-\left( \dfrac{r}{n}\right) a_{1}\right) }{{\left( \left( -1-\dfrac{r}{%
n(n-1)}\right) a_{0}+\left( -(n-1)-\dfrac{r}{n}\right) a_{1}\right) }}g$ \\
\hline
$(\varepsilon )$-Sasakian & $\dfrac{(\varepsilon )^{\ell }(n-1)^{\ell }a_{1}%
}{{\left( \left( \varepsilon -\dfrac{r}{n(n-1)}\right) a_{0}+\left(
\varepsilon (n-1)-\dfrac{r}{n}\right) a_{1}\right) }}S$ \\
\hline
\end{tabular}%
\
\]
\end{cor}

\medskip

\begin{cor}
Let $M$ be an $n$-dimensional $(R,{\cal P}_{\ast },S)$-pseudosymmetric $%
\left( N(k),\xi \right) $-semi-Riemannian manifold such that $%
a_{0}+(n-1)a_{1}\not=0$. If $M$ is not pseudo-projective semisymmetric, then
either $M$ is an Einstein manifold or $r=\dfrac{kn(n-1)a_{0}}{%
a_{0}+(n-1)a_{1}}$ and $L_{S}=\dfrac{1}{n-1}$. Consequently, we have the
following\/{\rm :}~
\[
\begin{tabular}{|l|l|l|}
\hline
${\boldmath M}$ & $L_{S}=$ & ${\boldmath S=}$ \\ \hline
$N(k)$-contact metric & $\dfrac{1}{n-1}$ & $S=k(n-1)g\quad {\rm or\quad }r=%
\dfrac{kn(n-1)a_{0}}{a_{0}+(n-1)a_{1}}$ \\ \hline
Sasakian & $\dfrac{1}{n-1}$ & $S=(n-1)g\quad {\rm or\quad }r=\dfrac{%
n(n-1)a_{0}}{a_{0}+(n-1)a_{1}}$ \\ \hline
Kenmotsu & $\dfrac{1}{n-1}$ & $S=-(n-1)g\quad {\rm or\quad }r=-\,\dfrac{%
n(n-1)a_{0}}{a_{0}+(n-1)a_{1}}$ \\ \hline
$(\varepsilon )$-Sasakian & $\dfrac{1}{n-1}$ & $S=\varepsilon (n-1)g\quad
{\rm or\quad }r=\dfrac{\varepsilon n(n-1)a_{0}}{a_{0}+(n-1)a_{1}}$ \\ \hline
para-Sasakian & $\dfrac{1}{n-1}$ & $S=-(n-1)g\quad {\rm or\quad }r=-\dfrac{%
n(n-1)a_{0}}{a_{0}+(n-1)a_{1}}$ \\ \hline
$(\varepsilon )$-para-Sasakian & $\dfrac{1}{n-1}$ & $S=-\,\varepsilon
(n-1)g\quad {\rm or\quad }r=-\dfrac{\varepsilon n(n-1)a_{0}}{a_{0}+(n-1)a_{1}%
}$ \\ \hline
\end{tabular}%
\
\]
\end{cor}

\begin{cor}
Let $M$ be an $n$-dimensional $(R,{\cal P},S^{\ell })$-pseudosymmetric $%
\left( N(k),\xi \right) $-semi-Riemannian manifold. If $M$ is not projective
semisymmetric, then $M$ is an Einstein manifold with scalar curvature $%
kn(n-1)$. Consequently, we have the following\/{\rm :}~
\[
\begin{tabular}{|l|l|}
\hline
${\boldmath M}$ & ${\boldmath S=}$ \\ \hline
$N(k)$-contact metric & $k(n-1)g$ \\ \hline
Sasakian & $(n-1)g$ \\ \hline
Kenmotsu & $-\,(n-1)g$ \\ \hline
$(\varepsilon )$-Sasakian & $\varepsilon (n-1)g$ \\ \hline
para-Sasakian & $-\,(n-1)g$ \\ \hline
$(\varepsilon )$-para-Sasakian & $-\,\varepsilon (n-1)g$ \\ \hline
\end{tabular}%
\]
\end{cor}

\begin{cor}
Let $M$ be an $n$-dimensional $(R,{\cal P},S)$-pseudosymmetric $\left(
N(k),\xi \right) $-semi-Riemannian manifold. If $M$ is not projective
semisymmetric, then $M$ is an Einstein manifold with scalar curvature $%
kn(n-1)$ and $L_{S}=\dfrac{1}{n-1}$. Consequently, we have the following\/%
{\rm :}~
\[
\begin{tabular}{|l|l|l|}
\hline
${\boldmath M}$ & $L_{S}=$ & ${\boldmath S=}$ \\ \hline
$N(k)$-contact metric & $\frac{1}{n-1}$ & $k(n-1)g$ \\ \hline
Sasakian & $\frac{1}{n-1}$ & $(n-1)g$ \\ \hline
Kenmotsu {\rm \cite{Ozgur-06}} & $\frac{1}{n-1}$ & $-(n-1)g$ \\ \hline
$(\varepsilon )$-Sasakian & $\frac{1}{n-1}$ & $\varepsilon (n-1)g$ \\ \hline
para-Sasakian & $\frac{1}{n-1}$ & $-(n-1)g$ \\ \hline
$(\varepsilon )$-para-Sasakian & $\frac{1}{n-1}$ & $-\,\varepsilon (n-1)g$
\\ \hline
\end{tabular}%
\]
\end{cor}

\begin{cor}
Let $M$ be an $n$-dimensional $(R,{\cal M},S^{\ell })$-pseudosymmetric $%
\left( N(k),\xi \right) $-semi-Riemannian manifold. If $M$ is not ${\cal M}$%
-semisymmetric, then
\[
S^{\ell +1}=k(n-1)S^{\ell }+k^{\ell }(n-1)^{\ell }S-k^{\ell +1}(n-1)^{\ell
+1}g.
\]%
Consequently, we have the following\/{\rm :}~
\[
\begin{tabular}{|l|l|}
\hline
${\boldmath M}$ & ${\boldmath S^{\ell +1}=}$ \\ \hline
$N(k)$-contact metric & $k(n-1)S^{\ell }+k^{\ell }(n-1)^{\ell }S-k^{\ell
+1}(n-1)^{\ell +1}g$ \\ \hline
Sasakian & $(n-1)S^{\ell }+(n-1)^{\ell }S-(n-1)^{\ell +1}g$ \\ \hline
Kenmotsu & $-(n-1)S^{\ell }+(-1)^{\ell }(n-1)^{\ell }S-(-1)^{\ell
+1}(n-1)^{\ell +1}g$ \\ \hline
$(\varepsilon )$-Sasakian & $\varepsilon (n-1)S^{\ell }+(\varepsilon )^{\ell
}(n-1)^{\ell }S-(\varepsilon )^{\ell +1}(n-1)^{\ell +1}g$ \\ \hline
para-Sasakian & $-(n-1)S^{\ell }+(-1)^{\ell }(n-1)^{\ell }S-(-1)^{\ell
+1}(n-1)^{\ell +1}g$ \\ \hline
$(\varepsilon )$-para-Sasakian & $-\varepsilon (n-1)S^{\ell }+(-\varepsilon
)^{\ell }(n-1)^{\ell }S-(-\varepsilon )^{\ell +1}(n-1)^{\ell +1}g$ \\ \hline
\end{tabular}%
\]
\end{cor}

\begin{cor}
Let $M$ be an $n$-dimensional $(R,{\cal M},S)$-pseudosymmetric $\left(
N(k),\xi \right) $-semi-Riemannian manifold. If $M$ is not ${\cal M}$%
-semisymmetric, then
\[
S^{2}=2(n-1)kS-k^{2}(n-1)^{2}g.
\]%
Consequently, we have the following\/{\rm :}~
\[
\begin{tabular}{|l|l|}
\hline
${\boldmath M}$ & ${\boldmath S}^{2}{=}$ \\ \hline
$N(k)$-contact metric & $2(n-1)kS-k^{2}(n-1)^{2}g$ \\ \hline
Sasakian & $2(n-1)S-(n-1)^{2}g$ \\ \hline
Kenmotsu & $-2(n-1)S-(n-1)^{2}g$ \\ \hline
$(\varepsilon )$-Sasakian & $2(n-1)\varepsilon S-(n-1)^{2}g$ \\ \hline
para-Sasakian & $-2(n-1)S-(n-1)^{2}g$ \\ \hline
$(\varepsilon )$-para-Sasakian & $-2(n-1)\varepsilon S-(n-1)^{2}g$ \\ \hline
\end{tabular}%
\]
\end{cor}

\begin{cor}
Let $M$ be an $n$-dimensional $(R,{\cal W}_{0},S^{\ell })$-pseudosymmetric $%
\left( N(k),\xi \right) $-semi-Riemannian manifold. If $M$ is not ${\cal W}%
_{0}$-semisymmetric, then
\[
S^{\ell +1}=k(n-1)S^{\ell }+k^{\ell }(n-1)^{\ell }S-k^{\ell +1}(n-1)^{\ell
+1}g.
\]%
Consequently, we have the following\/{\rm :}~
\[
\begin{tabular}{|l|l|}
\hline
${\boldmath M}$ & ${\boldmath S^{\ell +1}=}$ \\ \hline
$N(k)$-contact metric & $k(n-1)S^{\ell }+k^{\ell }(n-1)^{\ell }S-k^{\ell
+1}(n-1)^{\ell +1}g$ \\ \hline
Sasakian & $(n-1)S^{\ell }+(n-1)^{\ell }S-(n-1)^{\ell +1}g$ \\ \hline
Kenmotsu & $-(n-1)S^{\ell }+(-1)^{\ell }(n-1)^{\ell }S-(-1)^{\ell
+1}(n-1)^{\ell +1}g$ \\ \hline
$(\varepsilon )$-Sasakian & $\varepsilon (n-1)S^{\ell }+(\varepsilon )^{\ell
}(n-1)^{\ell }S-(\varepsilon )^{\ell +1}(n-1)^{\ell +1}g$ \\ \hline
para-Sasakian & $-(n-1)S^{\ell }+(-1)^{\ell }(n-1)^{\ell }S-(-1)^{\ell
+1}(n-1)^{\ell +1}g$ \\ \hline
$(\varepsilon )$-para-Sasakian & $-\varepsilon (n-1)S^{\ell }+(-\varepsilon
)^{\ell }(n-1)^{\ell }S-(-\varepsilon )^{\ell +1}(n-1)^{\ell +1}g$ \\ \hline
\end{tabular}%
\]
\end{cor}

\begin{cor}
Let $M$ be an $n$-dimensional $(R,{\cal W}_{0},S)$-pseudosymmetric $\left(
N(k),\xi \right) $-semi-Riemannian manifold. If $M$ is not ${\cal W}_{0}$%
-semisymmetric, then
\[
S^{2}=2(n-1)kS-k^{2}(n-1)^{2}g.
\]%
Consequently, we have the following\/{\rm :}~
\[
\begin{tabular}{|l|l|}
\hline
${\boldmath M}$ & ${\boldmath S}^{2}{=}$ \\ \hline
$N(k)$-contact metric & $2(n-1)kS-k^{2}(n-1)^{2}g$ \\ \hline
Sasakian & $2(n-1)S-(n-1)^{2}g$ \\ \hline
Kenmotsu & $-2(n-1)S-(n-1)^{2}g$ \\ \hline
$(\varepsilon )$-Sasakian & $2(n-1)\varepsilon S-(n-1)^{2}g$ \\ \hline
para-Sasakian & $-2(n-1)S-(n-1)^{2}g$ \\ \hline
$(\varepsilon )$-para-Sasakian & $-2(n-1)\varepsilon S-(n-1)^{2}g$ \\ \hline
\end{tabular}%
\]
\end{cor}

\begin{cor}
Let $M$ be an $n$-dimensional $(R,{\cal W}_{0}^{\ast },S^{\ell })$%
-pseudosymmetric $\left( N(k),\xi \right) $-semi-Riemannian manifold. If $M$
is not ${\cal W}_{0}^{\ast }$-semisymmetric, then
\[
S^{\ell +1}=-k(n-1)S^{\ell }+k^{\ell }(n-1)^{\ell }S+k^{\ell +1}(n-1)^{\ell
+1}g.
\]%
Consequently, we have the following\/{\rm :}~
\[
\begin{tabular}{|l|l|}
\hline
${\boldmath M}$ & ${\boldmath S^{\ell +1}=}$ \\ \hline
$N(k)$-contact metric & $-k(n-1)S^{\ell }+k^{\ell }(n-1)^{\ell }S+k^{\ell
+1}(n-1)^{\ell +1}g$ \\ \hline
Sasakian & $-(n-1)S^{\ell }+(n-1)^{\ell }S+(n-1)^{\ell +1}g$ \\ \hline
Kenmotsu & $(n-1)S^{\ell }+(-1)^{\ell }(n-1)^{\ell }S+(-1)^{\ell
+1}(n-1)^{\ell +1}g$ \\ \hline
$(\varepsilon )$-Sasakian & $-\varepsilon (n-1)S^{\ell }+(\varepsilon
)^{\ell }(n-1)^{\ell }S+(\varepsilon )^{\ell +1}(n-1)^{\ell +1}g$ \\ \hline
para-Sasakian & $(n-1)S^{\ell }+(-1)^{\ell }(n-1)^{\ell }S+(-1)^{\ell
+1}(n-1)^{\ell +1}g$ \\ \hline
$(\varepsilon )$-para-Sasakian & $\varepsilon (n-1)S^{\ell }+(-\varepsilon
)^{\ell }(n-1)^{\ell }S+(-\varepsilon )^{\ell +1}(n-1)^{\ell +1}g$ \\ \hline
\end{tabular}%
\]
\end{cor}

\begin{cor}
Let $M$ be an $n$-dimensional $(R,{\cal W}_{0}^{\ast },S)$-pseudosymmetric $%
\left( N(k),\xi \right) $-semi-Riemannian manifold. If $M$ is not ${\cal W}%
_{0}^{\ast }$-semisymmetric, then
\[
S^{2}=k^{2}(n-1)^{2}g.
\]%
Consequently, we have the following\/{\rm :}~
\[
\begin{tabular}{|l|l|}
\hline
${\boldmath M}$ & ${\boldmath S}^{2}{=}$ \\ \hline
$N(k)$-contact metric & $k^{2}(n-1)^{2}g$ \\ \hline
Sasakian & $(n-1)^{2}g$ \\ \hline
Kenmotsu & $(n-1)^{2}g$ \\ \hline
$(\varepsilon )$-Sasakian & $(n-1)^{2}g$ \\ \hline
para-Sasakian & $(n-1)^{2}g$ \\ \hline
$(\varepsilon )$-para-Sasakian & $(n-1)^{2}g$ \\ \hline
\end{tabular}%
\]
\end{cor}

\begin{cor}
Let $M$ be an $n$-dimensional $(R,{\cal W}_{1},S^{\ell })$-pseudosymmetric $%
\left( N(k),\xi \right) $-semi-Riemannian manifold. If $M$ is not ${\cal W}%
_{1}$-semisymmetric, then
\[
2S^{\ell }=k^{\ell -1}(n-1)^{\ell -1}S+k^{\ell }(n-1)^{\ell }g.
\]%
Consequently, we have the following\/{\rm :}~
\[
\begin{tabular}{|l|l|}
\hline
${\boldmath M}$ & ${\boldmath2S^{\ell }=}$ \\ \hline
$N(k)$-contact metric & $k^{\ell -1}(n-1)^{\ell -1}S+k^{\ell }(n-1)^{\ell }g$
\\ \hline
Sasakian & $(n-1)^{\ell -1}S+(n-1)^{\ell }g$ \\ \hline
Kenmotsu & $(-1)^{\ell -1}(n-1)^{\ell -1}S+(-1)^{\ell }(n-1)^{\ell }g$ \\
\hline
$(\varepsilon )$-Sasakian & $(\varepsilon )^{\ell -1}(n-1)^{\ell
-1}S+(\varepsilon )^{\ell }(n-1)^{\ell }g$ \\ \hline
para-Sasakian & $(-1)^{\ell -1}(n-1)^{\ell -1}S+(-1)^{\ell }(n-1)^{\ell }g$
\\ \hline
$(\varepsilon )$-para-Sasakian & $(-\varepsilon )^{\ell -1}(n-1)^{\ell
-1}S+(-\varepsilon )^{\ell }(n-1)^{\ell }g$ \\ \hline
\end{tabular}%
\]
\end{cor}

\begin{cor}
Let $M$ be an $n$-dimensional $(R,{\cal W}_{1},S)$-pseudosymmetric $\left(
N(k),\xi \right) $-semi-Riemannian manifold. If $M$ is not ${\cal W}_{1}$%
-semisymmetric, then $M$ is an Einstein manifold with scalar curvature $%
kn(n-1)$ and $L_{S}=\dfrac{1}{n-1}$. Consequently, we have the following\/%
{\rm :}~
\[
\begin{tabular}{|l|l|l|}
\hline
${\boldmath M}$ & $L_{S}=$ & ${\boldmath S=}$ \\ \hline
$N(k)$-contact metric & $\dfrac{1}{n-1}$ & $k(n-1)g$ \\ \hline
Sasakian & $\dfrac{1}{n-1}$ & $(n-1)g$ \\ \hline
Kenmotsu & $\dfrac{1}{n-1}$ & $-\,(n-1)g$ \\ \hline
$(\varepsilon )$-Sasakian & $\dfrac{1}{n-1}$ & $\varepsilon (n-1)g$ \\ \hline
para-Sasakian & $\dfrac{1}{n-1}$ & $-\,(n-1)g$ \\ \hline
$(\varepsilon )$-para-Sasakian & $\dfrac{1}{n-1}$ & $-\,\varepsilon (n-1)g$
\\ \hline
\end{tabular}%
\]
\end{cor}

\begin{cor}
Let $M$ be an $n$-dimensional $(R,{\cal W}_{1}^{\ast },S^{\ell })$%
-pseudosymmetric $\left( N(k),\xi \right) $-semi-Riemannian manifold. If $M$
is not ${\cal W}_{1}^{\ast }$-semisymmetric, then $M$ is an Einstein
manifold with scalar curvature $kn(n-1)$. Consequently, we have the
following\/{\rm :}~
\[
\begin{tabular}{|l|l|}
\hline
${\boldmath M}$ & ${\boldmath S=}$ \\ \hline
$N(k)$-contact metric & $k(n-1)g$ \\ \hline
Sasakian & $(n-1)g$ \\ \hline
Kenmotsu & $-(n-1)g$ \\ \hline
$(\varepsilon )$-Sasakian & $\varepsilon (n-1)g$ \\ \hline
para-Sasakian & $-(n-1)g$ \\ \hline
$(\varepsilon )$-para-Sasakian & $-\,\varepsilon (n-1)g$ \\ \hline
\end{tabular}%
\]
\end{cor}

\begin{cor}
Let $M$ be an $n$-dimensional $(R,{\cal W}_{1}^{\ast },S)$-pseudosymmetric $%
\left( N(k),\xi \right) $-semi-Riemannian manifold. If $M$ is not ${\cal W}%
_{1}^{\ast }$-semisymmetric, then $M$ is an Einstein manifold with scalar
curvature $kn(n-1)$ and $L_{S}=\dfrac{1}{n-1}$. Consequently, we have the
following\/{\rm :}~
\[
\begin{tabular}{|l|l|l|}
\hline
${\boldmath M}$ & $L_{S}=$ & ${\boldmath S=}$ \\ \hline
$N(k)$-contact metric & $\dfrac{1}{n-1}$ & $k(n-1)g$ \\ \hline
Sasakian & $\dfrac{1}{n-1}$ & $(n-1)g$ \\ \hline
Kenmotsu & $\dfrac{1}{n-1}$ & $-(n-1)g$ \\ \hline
$(\varepsilon )$-Sasakian & $\dfrac{1}{n-1}$ & $\varepsilon (n-1)g$ \\ \hline
para-Sasakian & $\dfrac{1}{n-1}$ & $-(n-1)g$ \\ \hline
$(\varepsilon )$-para-Sasakian & $\dfrac{1}{n-1}$ & $-\,\varepsilon (n-1)g$
\\ \hline
\end{tabular}%
\]
\end{cor}

\begin{cor}
Let $M$ be an $n$-dimensional $(R,{\cal W}_{2},S^{\ell })$-pseudosymmetric $%
\left( N(k),\xi \right) $-semi-Riemannian manifold. If $M$ is not ${\cal W}%
_{2}$-semisymmetric, then
\[
L_{S^{\ell }}S^{\ell +1}=k(n-1)L_{S^{\ell }}S^{\ell }+kS-k^{2}(n-1)g.
\]%
Consequently, we have the following\/{\rm :}~
\[
\begin{tabular}{|l|l|}
\hline
${\boldmath M}$ & ${\boldmath L_{S^{\ell }}S^{\ell +1}=}$ \\ \hline
$N(k)$-contact metric & $k(n-1)L_{S^{\ell }}S^{\ell }+kS-k^{2}(n-1)g$ \\
\hline
Sasakian & $(n-1)L_{S^{\ell }}S^{\ell }+S-(n-1)g$ \\ \hline
Kenmotsu & $-(n-1)L_{S^{\ell }}S^{\ell }-S-(n-1)g$ \\ \hline
$(\varepsilon )$-Sasakian & $\varepsilon (n-1)L_{S^{\ell }}S^{\ell
}+\varepsilon S-(n-1)g$ \\ \hline
para-Sasakian & $-(n-1)L_{S^{\ell }}S^{\ell }-S-(n-1)g$ \\ \hline
$(\varepsilon )$-para-Sasakian & $-\varepsilon (n-1)L_{S^{\ell }}S^{\ell
}-\varepsilon S-(n-1)g$ \\ \hline
\end{tabular}%
\]
\end{cor}

\begin{cor}
Let $M$ be an $n$-dimensional $(R,{\cal W}_{2},S)$-pseudosymmetric $\left(
N(k),\xi \right) $-semi-Riemannian manifold. If $M$ is not ${\cal W}_{2}$%
-semisymmetric, then
\[
L_{S}S^{2}=k\left( (n-1)L_{S}+1\right) S-k^{2}(n-1)g.
\]%
Consequently, we have the following\/{\rm :}~
\[
\begin{tabular}{|l|l|}
\hline
${\boldmath M}$ & ${\boldmath L_{S}S^{2}=}$ \\ \hline
$N(k)$-contact metric & $k\left( (n-1)L_{S}+1\right) S-k^{2}(n-1)g$ \\ \hline
Sasakian & $\left( (n-1)L_{S}+1\right) S-(n-1)g$ \\ \hline
Kenmotsu & $-\,\left( (n-1)L_{S}+1\right) S-(n-1)g$ \\ \hline
$(\varepsilon )$-Sasakian & $\varepsilon \left( (n-1)L_{S}+1\right) S-(n-1)g$
\\ \hline
para-Sasakian & $-\,\left( (n-1)L_{S}+1\right) S-(n-1)g$ \\ \hline
$(\varepsilon )$-para-Sasakian & $-\,\varepsilon \left( (n-1)L_{S}+1\right)
S-(n-1)g$ \\ \hline
\end{tabular}%
\]
\end{cor}

\begin{cor}
Let $M$ be an $n$-dimensional $(R,{\cal W}_{3},S^{\ell })$-pseudosymmetric $%
\left( N(k),\xi \right) $-semi-Riemannian manifold. If $M$ is not ${\cal W}%
_{3}$-semisymmetric, then
\[
S^{\ell }=k^{\ell }(n-1)^{\ell }g
\]%
and $L_{S^{\ell }}=\dfrac{1}{k^{\ell -1}(n-1)^{\ell }}$. Consequently, we
have the following\/{\rm :}~
\[
\begin{tabular}{|l|l|l|}
\hline
${\boldmath M}$ & $L_{S^{\ell }}=$ & ${\boldmath S}^{\ell }{=}$ \\ \hline
$N(k)$-contact metric & $\dfrac{1}{k^{\ell -1}(n-1)^{\ell }}$ & $k^{\ell
}(n-1)^{\ell }g$ \\ \hline
Sasakian & $\dfrac{1}{(n-1)^{\ell }}$ & $(n-1)^{\ell }g$ \\ \hline
Kenmotsu & $\dfrac{1}{(-1)^{\ell -1}(n-1)^{\ell }}$ & $(-1)^{\ell
}(n-1)^{\ell }g$ \\ \hline
$(\varepsilon )$-Sasakian & $\dfrac{1}{(\varepsilon )^{\ell -1}(n-1)^{\ell }}
$ & $(\varepsilon )^{\ell }(n-1)^{\ell }g$ \\ \hline
para-Sasakian & $\dfrac{1}{(-1)^{\ell -1}(n-1)^{\ell }}$ & $(-1)^{\ell
}(n-1)^{\ell }g$ \\ \hline
$(\varepsilon )$-para-Sasakian & $\dfrac{1}{(-\varepsilon )^{\ell
-1}(n-1)^{\ell }}$ & $(-\varepsilon )^{\ell }(n-1)^{\ell }g$ \\ \hline
\end{tabular}%
\]
\end{cor}

\begin{cor}
Let $M$ be an $n$-dimensional $(R,{\cal W}_{3},S)$-pseudosymmetric $\left(
N(k),\xi \right) $-semi-Riemannian manifold. If $M$ is not ${\cal W}_{3}$%
-semisymmetric, then $M$ is an Einstein manifold with scalar curvature $%
kn(n-1)$ and $L_{S}=\dfrac{1}{n-1}$. Consequently, we have the following\/%
{\rm :}~
\[
\begin{tabular}{|l|l|l|}
\hline
${\boldmath M}$ & $L_{S}=$ & ${\boldmath S=}$ \\ \hline
$N(k)$-contact metric & $\dfrac{1}{n-1}$ & $k(n-1)g$ \\ \hline
Sasakian & $\dfrac{1}{n-1}$ & $(n-1)g$ \\ \hline
Kenmotsu & $\dfrac{1}{n-1}$ & $-\,(n-1)g$ \\ \hline
$(\varepsilon )$-Sasakian & $\dfrac{1}{n-1}$ & $\varepsilon (n-1)g$ \\ \hline
para-Sasakian & $\dfrac{1}{n-1}$ & $-\,(n-1)g$ \\ \hline
$(\varepsilon )$-para-Sasakian & $\dfrac{1}{n-1}$ & $-\,\varepsilon (n-1)g$
\\ \hline
\end{tabular}%
\]
\end{cor}

\begin{cor}
Let $M$ be an $n$-dimensional $(R,{\cal W}_{4},S^{\ell })$-pseudosymmetric $%
\left( N(k),\xi \right) $-semi-Riemannian manifold. If $M$ is not ${\cal W}%
_{4}$-semisymmetric, then
\[
L_{S^{\ell }}S^{\ell +1}=k(n-1)L_{S^{\ell }}S^{\ell }+kS-k^{2}(n-1)g.
\]%
Consequently, we have the following\/{\rm :}~
\[
\begin{tabular}{|l|l|}
\hline
${\boldmath M}$ & ${\boldmath L_{S^{\ell }}S^{\ell +1}=}$ \\ \hline
$N(k)$-contact metric & $k(n-1)L_{S^{\ell }}S^{\ell }+kS-k^{2}(n-1)g$ \\
\hline
Sasakian & $(n-1)L_{S^{\ell }}S^{\ell }+S-(n-1)g$ \\ \hline
Kenmotsu & $-(n-1)L_{S^{\ell }}S^{\ell }-S-(n-1)g$ \\ \hline
$(\varepsilon )$-Sasakian & $\varepsilon (n-1)L_{S^{\ell }}S^{\ell
}+\varepsilon S-(n-1)g$ \\ \hline
para-Sasakian & $-(n-1)L_{S^{\ell }}S^{\ell }-S-(n-1)g$ \\ \hline
$(\varepsilon )$-para-Sasakian & $-\varepsilon (n-1)L_{S^{\ell }}S^{\ell
}-\varepsilon S-(n-1)g$ \\ \hline
\end{tabular}%
\]
\end{cor}

\begin{cor}
Let $M$ be an $n$-dimensional $(R,{\cal W}_{4},S)$-pseudosymmetric $\left(
N(k),\xi \right) $-semi-Riemannian manifold. If $M$ is not ${\cal W}_{4}$%
-semisymmetric, then
\[
L_{S}S^{2}=k\left( (n-1)L_{S}+1\right) S-k^{2}(n-1)g.
\]%
Consequently, we have the following\/{\rm :}~
\[
\begin{tabular}{|l|l|}
\hline
${\boldmath M}$ & ${\boldmath L_{S}S^{2}=}$ \\ \hline
$N(k)$-contact metric & $k\left( (n-1)L_{S}+1\right) S-k^{2}(n-1)g$ \\ \hline
Sasakian & $\left( (n-1)L_{S}+1\right) S-(n-1)g$ \\ \hline
Kenmotsu & $-\,\left( (n-1)L_{S}+1\right) S-(n-1)g$ \\ \hline
$(\varepsilon )$-Sasakian & $\varepsilon \left( (n-1)L_{S}+1\right) S-(n-1)g$
\\ \hline
para-Sasakian & $-\,\left( (n-1)L_{S}+1\right) S-(n-1)g$ \\ \hline
$(\varepsilon )$-para-Sasakian & $-\,\varepsilon \left( (n-1)L_{S}+1\right)
S-(n-1)g$ \\ \hline
\end{tabular}%
\]
\end{cor}

\begin{cor}
Let $M$ be an $n$-dimensional $(R,{\cal W}_{5},S^{\ell })$-pseudosymmetric $%
\left( N(k),\xi \right) $-semi-Riemannian manifold. If $M$ is not ${\cal W}%
_{5}$-semisymmetric, then either
\[
S=k(n-1)g
\]%
or $L_{S^{\ell }}=\dfrac{1}{k^{\ell -1}(n-1)^{\ell }}$. Consequently, we
have the following\/{\rm :}~
\[
\begin{tabular}{|l|l|l|}
\hline
${\boldmath M}$ & $L_{S^{\ell }}=$ & ${\boldmath S}=$ \\ \hline
$N(k)$-contact metric & $\dfrac{1}{k^{\ell -1}(n-1)^{\ell }}$ & $k(n-1)g$ \\
\hline
Sasakian & $\dfrac{1}{(n-1)^{\ell }}$ & $(n-1)g$ \\ \hline
Kenmotsu & $\dfrac{1}{(-1)^{\ell -1}(n-1)^{\ell }}$ & $-(n-1)g$ \\ \hline
$(\varepsilon )$-Sasakian & $\dfrac{1}{(-1)^{\ell -1}(n-1)^{\ell }}$ & $%
\varepsilon (n-1)g$ \\ \hline
para-Sasakian & $\dfrac{1}{(-1)^{\ell -1}(n-1)^{\ell }}$ & $-(n-1)g$ \\
\hline
$(\varepsilon )$-para-Sasakian & $\dfrac{1}{(-\varepsilon )^{\ell
-1}(n-1)^{\ell }}$ & $-\varepsilon (n-1)g$ \\ \hline
\end{tabular}%
\]
\end{cor}

\begin{cor}
Let $M$ be an $n$-dimensional $(R,{\cal W}_{5},S)$-pseudosymmetric $\left(
N(k),\xi \right) $-semi-Riemannian manifold. If $M$ is not ${\cal W}_{5}$%
-semisymmetric, then $M$ is an Einstein manifold with scalar curvature $%
kn(n-1)$ and $L_{S}=\dfrac{1}{n-1}$. Consequently, we have the following\/%
{\rm :}~
\[
\begin{tabular}{|l|l|l|}
\hline
${\boldmath M}$ & $L_{S}=$ & ${\boldmath S=}$ \\ \hline
$N(k)$-contact metric & $\frac{1}{n-1}$ & $k(n-1)g$ \\ \hline
Sasakian & $\frac{1}{n-1}$ & $(n-1)g$ \\ \hline
Kenmotsu & $\frac{1}{n-1}$ & $-\,(n-1)g$ \\ \hline
$(\varepsilon )$-Sasakian & $\frac{1}{n-1}$ & $\varepsilon (n-1)g$ \\ \hline
para-Sasakian & $\frac{1}{n-1}$ & $-\,(n-1)g$ \\ \hline
$(\varepsilon )$-para-Sasakian & $\frac{1}{n-1}$ & $-\,\varepsilon (n-1)g$
\\ \hline
\end{tabular}%
\]
\end{cor}

\begin{cor}
Let $M$ be an $n$-dimensional $(R,{\cal W}_{6},S^{\ell })$-pseudosymmetric $%
\left( N(k),\xi \right) $-semi-Riemannian manifold. If $M$ is not ${\cal W}%
_{6}$-semisymmetric, then
\[
(n-1)L_{S^{\ell }}S^{\ell }=(1-k^{\ell -1}(n-1)^{\ell }L_{S^{\ell
}})S+(2k^{\ell }(n-1)^{\ell +1}L_{S^{\ell }}-k(n-1))g.
\]%
Consequently, we have the following\/{\rm :}~
\[
\begin{tabular}{|l|l|}
\hline
${\boldmath M}$ & ${\boldmath(n-1)L_{S^{\ell }}S^{\ell }=}$ \\ \hline
$N(k)$-contact metric & $(1-k^{\ell -1}(n-1)^{\ell }L_{S^{\ell
}})S+(2k^{\ell }(n-1)^{\ell +1}L_{S^{\ell }}-k(n-1))g$ \\ \hline
Sasakian & $(1-k^{\ell -1}(n-1)^{\ell }L_{S^{\ell }})S+(2k^{\ell
}(n-1)^{\ell +1}L_{S^{\ell }}-k(n-1))g$ \\ \hline
Kenmotsu & $(1-k^{\ell -1}(n-1)^{\ell }L_{S^{\ell }})S+(2k^{\ell
}(n-1)^{\ell +1}L_{S^{\ell }}-k(n-1))g$ \\ \hline
$(\varepsilon )$-Sasakian & $(1-k^{\ell -1}(n-1)^{\ell }L_{S^{\ell
}})S+(2k^{\ell }(n-1)^{\ell +1}L_{S^{\ell }}-k(n-1))g$ \\ \hline
para-Sasakian & $(1-k^{\ell -1}(n-1)^{\ell }L_{S^{\ell }})S+(2k^{\ell
}(n-1)^{\ell +1}L_{S^{\ell }}-k(n-1))g$ \\ \hline
$(\varepsilon )$-para-Sasakian & $(1-k^{\ell -1}(n-1)^{\ell }L_{S^{\ell
}})S+(2k^{\ell }(n-1)^{\ell +1}L_{S^{\ell }}-k(n-1))g$ \\ \hline
\end{tabular}%
\]
\end{cor}

\begin{cor}
Let $M$ be an $n$-dimensional $(R,{\cal W}_{6},S)$-pseudosymmetric $\left(
N(k),\xi \right) $-semi-Riemannian manifold. If $M$ is not ${\cal W}_{6}$%
-semisymmetric, then $M$ is an Einstein manifold with scalar curvature $%
kn(n-1)$ and $L_{S}=\dfrac{1}{2(n-1)}$. Consequently, we have the following\/%
{\rm :}~
\[
\begin{tabular}{|l|l|l|}
\hline
${\boldmath M}$ & $L_{S}=$ & ${\boldmath S=}$ \\ \hline
$N(k)$-contact metric & $\frac{1}{2(n-1)}$ & $k(n-1)g$ \\ \hline
Sasakian & $\frac{1}{2(n-1)}$ & $(n-1)g$ \\ \hline
Kenmotsu & $\frac{1}{2(n-1)}$ & $-\,(n-1)g$ \\ \hline
$(\varepsilon )$-Sasakian & $\frac{1}{2(n-1)}$ & $\varepsilon (n-1)g$ \\
\hline
para-Sasakian & $\frac{1}{2(n-1)}$ & $-\,(n-1)g$ \\ \hline
$(\varepsilon )$-para-Sasakian & $\frac{1}{2(n-1)}$ & $-\,\varepsilon (n-1)g$
\\ \hline
\end{tabular}%
\]
\end{cor}

\begin{cor}
Let $M$ be an $n$-dimensional $(R,{\cal W}_{7},S^{\ell })$-pseudosymmetric $%
\left( N(k),\xi \right) $-semi-Riemannian manifold. If $M$ is not ${\cal W}%
_{7}$-semisymmetric, then
\[
(n-1)L_{S^{\ell }}S^{\ell }=(1-k^{\ell -1}(n-1)^{\ell }L_{S^{\ell
}})S+(2k^{\ell }(n-1)^{\ell +1}L_{S^{\ell }}-k(n-1))g.
\]%
Consequently, we have the following\/{\rm :}~
\[
\begin{tabular}{|l|l|}
\hline
${\boldmath M}$ & ${\boldmath(n-1)L_{S^{\ell }}S^{\ell }=}$ \\ \hline
$N(k)$-contact metric & $(1-k^{\ell -1}(n-1)^{\ell }L_{S^{\ell
}})S+(2k^{\ell }(n-1)^{\ell +1}L_{S^{\ell }}-k(n-1))g$ \\ \hline
Sasakian & $(1-(n-1)^{\ell }L_{S^{\ell }})S+(2(n-1)^{\ell +1}L_{S^{\ell
}}-(n-1))g$ \\ \hline
Kenmotsu & $(1-(-1)^{\ell -1}(n-1)^{\ell }L_{S^{\ell }})S+(2(-1)^{\ell
}(n-1)^{\ell +1}L_{S^{\ell }}+(n-1))g$ \\ \hline
$(\varepsilon )$-Sasakian & $(1-(\varepsilon )^{\ell -1}(n-1)^{\ell
}L_{S^{\ell }})S+(2(\varepsilon )^{\ell }(n-1)^{\ell +1}L_{S^{\ell
}}-\varepsilon (n-1))g$ \\ \hline
para-Sasakian & $(1-(-1)^{\ell -1}(n-1)^{\ell }L_{S^{\ell }})S+(2(-1)^{\ell
}(n-1)^{\ell +1}L_{S^{\ell }}+(n-1))g$ \\ \hline
$(\varepsilon )$-para-Sasakian & $(1-(-\varepsilon )^{\ell -1}(n-1)^{\ell
}L_{S^{\ell }})S+(2(-\varepsilon )^{\ell }(n-1)^{\ell +1}L_{S^{\ell
}}+\varepsilon (n-1))g$ \\ \hline
\end{tabular}%
\]
\end{cor}

\begin{cor}
Let $M$ be an $n$-dimensional $(R,{\cal W}_{7},S)$-pseudosymmetric $\left(
N(k),\xi \right) $-semi-Riemannian manifold. If $M$ is not ${\cal W}_{7}$%
-semisymmetric, then $M$ is an Einstein manifold with scalar curvature $%
kn(n-1)$ and $L_{S}=\dfrac{1}{2(n-1)}$. Consequently, we have the following\/%
{\rm :}~
\[
\begin{tabular}{|l|l|l|}
\hline
${\boldmath M}$ & $L_{S}=$ & ${\boldmath S=}$ \\ \hline
$N(k)$-contact metric & $\frac{1}{2(n-1)}$ & $k(n-1)g$ \\ \hline
Sasakian & $\frac{1}{2(n-1)}$ & $(n-1)g$ \\ \hline
Kenmotsu & $\frac{1}{2(n-1)}$ & $-\,(n-1)g$ \\ \hline
$(\varepsilon )$-Sasakian & $\frac{1}{2(n-1)}$ & $\varepsilon (n-1)g$ \\
\hline
para-Sasakian & $\frac{1}{2(n-1)}$ & $-\,(n-1)g$ \\ \hline
$(\varepsilon )$-para-Sasakian & $\frac{1}{2(n-1)}$ & $-\,\varepsilon (n-1)g$
\\ \hline
\end{tabular}%
\]
\end{cor}

\begin{cor}
Let $M$ be an $n$-dimensional $(R,{\cal W}_{8},S^{\ell })$-pseudosymmetric $%
\left( N(k),\xi \right) $-semi-Riemannian manifold. If $M$ is not ${\cal W}%
_{8}$-semisymmetric, then
\[
S^{\ell }=k^{\ell }(n-1)^{\ell }g
\]%
and $L_{S^{\ell }}=\dfrac{1}{k^{\ell -1}(n-1)^{\ell }}$. Consequently, we
have the following\/{\rm :}~
\[
\begin{tabular}{|l|l|l|}
\hline
${\boldmath M}$ & $L_{S^{\ell }}=$ & ${\boldmath S}^{\ell }=$ \\ \hline
$N(k)$-contact metric & $\dfrac{1}{k^{\ell -1}(n-1)^{\ell }}$ & $k^{\ell
}(n-1)^{\ell }g$ \\ \hline
Sasakian & $\dfrac{1}{(n-1)^{\ell }}$ & $(n-1)^{\ell }g$ \\ \hline
Kenmotsu & $\dfrac{1}{(-1)^{\ell -1}(n-1)^{\ell }}$ & $(-1)^{\ell
}(n-1)^{\ell }g$ \\ \hline
$(\varepsilon )$-Sasakian & $\dfrac{1}{(\varepsilon )^{\ell -1}(n-1)^{\ell }}
$ & $(\varepsilon )^{\ell }(n-1)^{\ell }g$ \\ \hline
para-Sasakian & $\dfrac{1}{(-1)^{\ell -1}(n-1)^{\ell }}$ & $(-1)^{\ell
}(n-1)^{\ell }g$ \\ \hline
$(\varepsilon )$-para-Sasakian & $\dfrac{1}{(-\varepsilon )^{\ell
-1}(n-1)^{\ell }}$ & $(-\varepsilon )^{\ell }(n-1)^{\ell }g$ \\ \hline
\end{tabular}%
\]
\end{cor}

\begin{cor}
Let $M$ be an $n$-dimensional $(R,{\cal W}_{8},S)$-pseudosymmetric $\left(
N(k),\xi \right) $-semi-Riemannian manifold. If $M$ is not ${\cal W}_{8}$%
-semisymmetric, then $M$ is an Einstein manifold with scalar curvature $%
kn(n-1)$ and $L_{S}=\dfrac{1}{n-1}$. Consequently, we have the following\/%
{\rm :}~
\[
\begin{tabular}{|l|l|l|}
\hline
${\boldmath M}$ & $L_{S}=$ & ${\boldmath S=}$ \\ \hline
$N(k)$-contact metric & $\dfrac{1}{n-1}$ & $k(n-1)g$ \\ \hline
Sasakian & $\dfrac{1}{n-1}$ & $(n-1)g$ \\ \hline
Kenmotsu & $\dfrac{1}{n-1}$ & $-\,(n-1)g$ \\ \hline
$(\varepsilon )$-Sasakian & $\dfrac{1}{n-1}$ & $\varepsilon (n-1)g$ \\ \hline
para-Sasakian & $\dfrac{1}{n-1}$ & $-\,(n-1)g$ \\ \hline
$(\varepsilon )$-para-Sasakian & $\dfrac{1}{n-1}$ & $-\,\varepsilon (n-1)g$
\\ \hline
\end{tabular}%
\]
\end{cor}

\begin{cor}
Let $M$ be an $n$-dimensional $(R,{\cal W}_{9},S^{\ell })$-pseudosymmetric $%
\left( N(k),\xi \right) $-semi-Riemannian manifold. If $M$ is not ${\cal W}%
_{9}$-semisymmetric, then
\[
S^{\ell }=k^{\ell }(n-1)^{\ell }g
\]%
and $L_{S^{\ell }}=\dfrac{1}{k^{\ell -1}(n-1)^{\ell }}$. Consequently, we
have the following\/{\rm :}~
\[
\begin{tabular}{|l|l|l|}
\hline
${\boldmath M}$ & $L_{S^{\ell }}=$ & ${\boldmath S}^{\ell }{=}$ \\ \hline
$N(k)$-contact metric & $\dfrac{1}{k^{\ell -1}(n-1)^{\ell }}$ & $k^{\ell
}(n-1)^{\ell }g$ \\ \hline
Sasakian & $\dfrac{1}{(n-1)^{\ell }}$ & $(n-1)^{\ell }g$ \\ \hline
Kenmotsu & $\dfrac{1}{(-1)^{\ell -1}(n-1)^{\ell }}$ & $(-1)^{\ell
}(n-1)^{\ell }g$ \\ \hline
$(\varepsilon )$-Sasakian & $\dfrac{1}{(\varepsilon )^{\ell -1}(n-1)^{\ell }}
$ & $(\varepsilon )^{\ell }(n-1)^{\ell }g$ \\ \hline
para-Sasakian & $\dfrac{1}{(-1)^{\ell -1}(n-1)^{\ell }}$ & $(-1)^{\ell
}(n-1)^{\ell }g$ \\ \hline
$(\varepsilon )$-para-Sasakian & $\dfrac{1}{(-\varepsilon )^{\ell
-1}(n-1)^{\ell }}$ & $(-\varepsilon )^{\ell }(n-1)^{\ell }g$ \\ \hline
\end{tabular}%
\]
\end{cor}

\begin{cor}
Let $M$ be an $n$-dimensional $(R,{\cal W}_{9},S)$-pseudosymmetric $\left(
N(k),\xi \right) $-semi-Riemannian manifold. If $M$ is not ${\cal W}_{9}$%
-semisymmetric, then $M$ is an Einstein manifold with scalar curvature $%
kn(n-1)$ and $L_{S}=\dfrac{1}{n-1}$. Consequently, we have the following\/%
{\rm :}~
\[
\begin{tabular}{|l|l|l|}
\hline
${\boldmath M}$ & $L_{S}=$ & ${\boldmath S=}$ \\ \hline
$N(k)$-contact metric & $\dfrac{1}{n-1}$ & $k(n-1)g$ \\ \hline
Sasakian & $\dfrac{1}{n-1}$ & $(n-1)g$ \\ \hline
Kenmotsu & $\dfrac{1}{n-1}$ & $-\,(n-1)g$ \\ \hline
$(\varepsilon )$-Sasakian & $\dfrac{1}{n-1}$ & $\varepsilon (n-1)g$ \\ \hline
para-Sasakian & $\dfrac{1}{n-1}$ & $-\,(n-1)g$ \\ \hline
$(\varepsilon )$-para-Sasakian & $\dfrac{1}{n-1}$ & $-\,\varepsilon (n-1)g$
\\ \hline
\end{tabular}%
\]
\end{cor}

\section{$({\cal T}_{\!a},S_{{\cal T}_{b}})$-pseudosymmetry\label{sect-TSP}}

In this section, we determine the results for an $n$-dimensional $\left(
N(k),\xi \right) $-semi-Riemannian manifold satisfy ${\cal T}_{\!a}\cdot S_{%
{\cal T}_{b}}=LQ(g,S_{{\cal T}_{\!b}})$.

\begin{defn-new}
A semi-Riemannian manifold is said to be $({\cal T}_{\!a},S_{{\cal T}_{b}})$%
-pseudosymmetric if
\begin{equation}
{\cal T}_{\!a}\cdot S_{{\cal T}_{b}}=LQ(g,S_{{\cal T}_{\!b}}),
\label{T.S=LQ}
\end{equation}
where $L$ is some smooth function defined on $M$. In particular, it is said
to be $(R\cdot S_{{\cal T}_{a}})$-pseudosymmetric if it satisfies
\begin{equation}
R\cdot S_{{\cal T}_{a}}=LQ(g,S_{{\cal T}_{a}}),  \label{T.S=LQ1}
\end{equation}%
holds on the set ${\cal U}=\left\{ x\in M:\left( S_{{\cal T}_{a}}-\dfrac{%
tr(S_{{\cal T}_{a}})}{n}g\right) _{x}\not=0\right\} $, where $L$ is some
function defined on ${\cal U}$.
\end{defn-new}

\begin{rem-new}
If in {\rm (\ref{T.S=LQ1})}, $S_{{\cal T}_{a}}$ is replaced by $S$ then it
is said to be Ricci-pseudosymmetric.
\end{rem-new}

\begin{th}
Let $M$ be an $n$-dimensional $({\cal T}_{\!a},S_{{\cal T}_{b}})$%
-pseudosymmetric $(N(k),\xi )$-semi-Riemannian manifold. Then
\begin{eqnarray*}
&&\varepsilon a_{5}(b_{0}+nb_{1}+b_{2}+b_{3}+b_{5}+b_{6})S^{2}(Y,U) \\
&&+\ \left\{ \varepsilon (b_{0}+nb_{1}+b_{2}+b_{3}+b_{5}+b_{6})\times \right.
\\
&&(-ka_{0}+k(n-1)a_{1}+k(n-1)a_{2}-a_{7}r) \\
&&\qquad \left. +\ \varepsilon (a_{1}+a_{5})(b_{4}r+(n-1)b_{7}r)\right\}
S(Y,U) \\
&&+\ \left\{ \varepsilon k(n-1)(a_{2}+a_{4})(b_{4}r+(n-1)b_{7}r)\right. \\
&&\qquad +\ \varepsilon k(n-1)(b_{0}+nb_{1}+b_{2}+b_{3}+b_{5}+b_{6})\times \\
&&\left. (ka_{0}+k(n-1)a_{4}+a_{7}r)\right\} g(Y,U) \\
&&+\ k(n-1)(a_{1}+a_{2}+2a_{3}+a_{4}+a_{5}+2a_{6})\times \\
&&\left\{ b_{4}r+(n-1)b_{7}r\right. \\
&&\qquad \left. +\ k(n-1)(b_{0}+nb_{1}+b_{2}+b_{3}+b_{5}+b_{6})\right\} \eta
(Y)\eta (U) \\
&=&L(b_{0}+nb_{1}+b_{2}+b_{3}+b_{5}+b_{6})(\varepsilon
k(n-1)g(Y,U)-\varepsilon S(Y,U)).
\end{eqnarray*}%
In particular, if $M$ is an $n$-dimensional $({\cal T}_{\!a},S_{{\cal T}%
_{a}})$-pseudosymmetric $(N(k),\xi )$-semi-Riemannian manifold, then
\begin{eqnarray*}
&&\varepsilon a_{5}(a_{0}+na_{1}+a_{2}+a_{3}+a_{5}+a_{6})S^{2}(Y,U) \\
&&+\ \left\{ \varepsilon (a_{0}+na_{1}+a_{2}+a_{3}+a_{5}+a_{6})\times \right.
\\
&&(-ka_{0}+k(n-1)a_{1}+k(n-1)a_{2}-a_{7}r) \\
&&\qquad \left. +\ \varepsilon (a_{1}+a_{5})(a_{4}r+(n-1)a_{7}r)\right\}
S(Y,U) \\
&&+\ \left\{ \varepsilon k(n-1)(a_{2}+a_{4})(a_{4}r+(n-1)a_{7}r)\right. \\
&&\qquad +\ \varepsilon k(n-1)(a_{0}+na_{1}+a_{2}+a_{3}+a_{5}+a_{6})\times \\
&&\left. (ka_{0}+k(n-1)a_{4}+a_{7}r)\right\} g(Y,U) \\
&&+\ k(n-1)(a_{1}+a_{2}+2a_{3}+a_{4}+a_{5}+2a_{6}) \\
&&\left\{ a_{4}r+(n-1)a_{7}r\right. \\
&&\qquad \left. +\ \varepsilon
k(n-1)(a_{0}+na_{1}+a_{2}+a_{3}+a_{5}+a_{6})\right\} \eta (Y)\eta (U) \\
&=&L(a_{0}+na_{1}+a_{2}+a_{3}+a_{5}+a_{6})(\varepsilon
k(n-1)g(Y,U)-\varepsilon S(Y,U)).
\end{eqnarray*}
\end{th}

\noindent {\bf Proof.} Let $M$ be an $n$-dimensional $({\cal T}_{\!a},S_{%
{\cal T}_{b}})$-pseudosymmetric $(N(k),\xi )$-semi-Riemannian manifold. Then
\begin{equation}
{\cal T}_{\!a}(X,Y)\cdot S_{{\cal T}_{\!b}}(U,V)=LQ(g,S_{{\cal T}%
_{\!b}})(U,V;X,Y).  \label{eq-Ric-pseudo}
\end{equation}%
Taking $X=\xi =V$ in (\ref{eq-Ric-pseudo}), we have%
\[
{\cal T}_{\!a}(\xi ,Y)\cdot S_{{\cal T}_{\!b}}(U,\xi )=LQ(g,S_{{\cal T}%
_{\!b}})(U,\xi ;\xi ,Y),
\]%
which gives
\begin{eqnarray}
&&S_{{\cal T}_{\!b}}({\cal T}_{\!a}(\xi ,Y)U,\xi )+S_{{\cal T}_{\!b}}(U,%
{\cal T}_{\!a}(\xi ,Y)\xi )  \nonumber \\
&=&L\left( S_{{\cal T}_{\!b}}((\xi \wedge Y)U,\xi )+S_{{\cal T}%
_{\!b}}(U,(\xi \wedge Y)\xi )\right) .  \label{eq-T-S-11}
\end{eqnarray}%
Using (\ref{eq-cond}), (\ref{eq-ricci}), (\ref{eq-xi-X-xi}), (\ref{eq-xi-Y-Z}%
), (\ref{eq-ric-T1}) and (\ref{eq-ric-T2}) in (\ref{eq-T-S-11}), we get the
result. $\blacksquare $

\begin{th}
\label{GCT-rsss} Let $M$ be an $n$-dimensional $({\cal T}_{\!a},S)$%
-pseudosymmetric $(N(k),\xi )$-semi-Riemannian manifold. Then
\begin{eqnarray*}
&&\varepsilon a_{5}\,S^{2}(Y,U)-E\,S(Y,U)-Fg(Y,U)-G\eta (Y)\eta \left(
U\right) \\
&=&L(\varepsilon k(n-1)g(Y,U)-\varepsilon S(Y,U)),
\end{eqnarray*}%
where
\[
E=\varepsilon \,(ka_{0}+a_{7}r-k(n-1)a_{1}-k(n-1)a_{2}),
\]%
\[
F=-\,\varepsilon k(n-1)(ka_{0}+k(n-1)a_{4}+a_{7}r),
\]%
\[
G=-\,k^{2}(n-1)^{2}(a_{1}+a_{2}+2a_{3}+a_{4}+a_{5}+2a_{6}).
\]
\end{th}

In view of Theorem~\ref{GCT-rsss}, we have the following

\begin{cor}
Let $M$ be an $n$-dimensional Ricci-pseudosymmetric $\left( N(k),\xi \right)
$-semi-Riemannian manifold. Then we have the following table\/{\rm :}~%
\[
\begin{tabular}{|l|l|l|}
\hline
${\boldmath M}$ & ${\boldmath}L=$ & ${\boldmath S=}$ \\ \hline
$N(k)$-contact metric & $k$ & $k(n-1)g$ \\ \hline
Sasakian & $1$ & $(n-1)g$ \\ \hline
Kenmotsu & $-\,1$ & $-\,(n-1)g$ \\ \hline
$(\varepsilon )$-Sasakian & $\varepsilon $ & $\varepsilon (n-1)g$ \\ \hline
para-Sasakian & $-\,1$ & $-\,(n-1)g$ \\ \hline
$(\varepsilon )$-para-Sasakian & $-\,\varepsilon $ & $-\,\varepsilon (n-1)g$
\\ \hline
\end{tabular}%
\]
\end{cor}

\begin{cor}
Let $M$ be an $n$-dimensional $({\cal C}_{\ast },S)$-pseudosymmetric $\left(
N(k),\xi \right) $-semi-Riemannian manifold. Then we have the following
table\/{\rm :}~%
\[
\begin{tabular}{|l|l|}
\hline
${\boldmath M}$ & ${\boldmath S^{2}=}$ \\ \hline
$N(k)$-contact metric & $-\,\left( \left( k-\dfrac{r}{n(n-1)}\right) \dfrac{%
a_{0}}{a_{1}}-\dfrac{2r}{n}-\dfrac{L}{a_{1}}\right) S$ \\
& $+k(n-1)\left( \left( k-\dfrac{r}{n(n-1)}\right) \dfrac{a_{0}}{a_{1}}%
+k(n-1)-\dfrac{2r}{n}-\dfrac{L}{a_{1}}\right) g$ \\ \hline
Sasakian & $-\,\left( \left( 1-\dfrac{r}{n(n-1)}\right) \dfrac{a_{0}}{a_{1}}-%
\dfrac{2r}{n}-\dfrac{L}{a_{1}}\right) S$ \\
& $+(n-1)\left( \left( 1-\dfrac{r}{n(n-1)}\right) \dfrac{a_{0}}{a_{1}}+(n-1)-%
\dfrac{2r}{n}-\dfrac{L}{a_{1}}\right) g$ \\ \hline
Kenmotsu & $\left( \left( 1+\dfrac{r}{n(n-1)}\right) \dfrac{a_{0}}{a_{1}}+%
\dfrac{2r}{n}+\dfrac{L}{a_{1}}\right) S$ \\
& $+(n-1)\left( \left( 1+\dfrac{r}{n(n-1)}\right) \dfrac{a_{0}}{a_{1}}+(n-1)+%
\dfrac{2r}{n}+\dfrac{L}{a_{1}}\right) g$ \\ \hline
$(\varepsilon )$-Sasakian & $-\,\varepsilon \left( \left( 1-\dfrac{%
\varepsilon r}{n(n-1)}\right) \dfrac{a_{0}}{a_{1}}-\dfrac{2\varepsilon r}{n}-%
\dfrac{\varepsilon L}{a_{1}}\right) S$ \\
& $+\varepsilon (n-1)\left( \left( \varepsilon -\dfrac{r}{n(n-1)}\right)
\dfrac{a_{0}}{a_{1}}+\varepsilon (n-1)-\dfrac{2r}{n}-\dfrac{L}{a_{1}}\right)
g$ \\ \hline
para-Sasakian & $\left( \left( 1+\dfrac{r}{n(n-1)}\right) \dfrac{a_{0}}{a_{1}%
}+\dfrac{2r}{n}+\dfrac{L}{a_{1}}\right) S$ \\
& $+(n-1)\left( \left( 1+\dfrac{r}{n(n-1)}\right) \dfrac{a_{0}}{a_{1}}+(n-1)+%
\dfrac{2r}{n}+\dfrac{L}{a_{1}}\right) g$ \\ \hline
$(\varepsilon )$-para-Sasakian & $\varepsilon \left( \left( 1+\dfrac{%
\varepsilon r}{n(n-1)}\right) \dfrac{a_{0}}{a_{1}}+\dfrac{2\varepsilon r}{n}+%
\dfrac{\varepsilon L}{a_{1}}\right) S$ \\
& $+\varepsilon (n-1)\left( \left( \varepsilon +\dfrac{r}{n(n-1)}\right)
\dfrac{a_{0}}{a_{1}}+\varepsilon (n-1)+\dfrac{2r}{n}+\dfrac{L}{a_{1}}\right)
g$ \\ \hline
\end{tabular}%
\]
\end{cor}

\begin{cor}
Let $M$ be an $n$-dimensional $({\cal C},S)$-pseudosymmetric $\left(
N(k),\xi \right) $-semi-Riemannian manifold. Then we have the following
table\/{\rm :}~
\[
\begin{tabular}{|l|l|}
\hline
${\boldmath M}$ & ${\boldmath S^{2}=}$ \\ \hline
$N(k)$-contact metric & $\left( \dfrac{r}{n-1}+\left( k-L\right)
(n-2)\right) S-k(r-\left( L(n-2)+k\right) (n-1))g$ \\ \hline
Sasakian & $\left( \dfrac{r}{n-1}+\left( 1-L\right) (n-2)\right) S-(r-\left(
L(n-2)+1\right) (n-1))g$ \\ \hline
Kenmotsu & $\left( \dfrac{r}{n-1}-\left( 1+L\right) (n-2)\right) S+(r-\left(
L(n-2)-1\right) (n-1))g$ \\ \hline
$(\varepsilon )$-Sasakian & $\left( \dfrac{r}{n-1}+\left( \varepsilon
-L\right) (n-2)\right) S-\varepsilon (r-\left( L(n-2)+\varepsilon \right)
(n-1))g$ \\ \hline
para-Sasakian & $\left( \dfrac{r}{n-1}-\left( 1+L\right) (n-2)\right)
S+1(r-\left( L(n-2)-1\right) (n-1))g$ \\ \hline
$(\varepsilon )$-para-Sasakian & $\left( \dfrac{r}{n-1}-\left( \varepsilon
+L\right) (n-2)\right) S+\varepsilon (r-\left( L(n-2)-\varepsilon \right)
(n-1))g$ \\ \hline
\end{tabular}%
\]
\end{cor}

\begin{cor}
Let $M$ be an $n$-dimensional $({\cal L},S)$-pseudosymmetric $\left(
N(k),\xi \right) $-semi-Riemannian manifold. Then we have the following
table\/{\rm :}~%
\[
\begin{tabular}{|l|l|}
\hline
${\boldmath M}$ & ${\boldmath S^{2}=}$ \\ \hline
$N(k)$-contact metric & $(n-2)\left( k-L\right) S+k(n-1)(k+(n-2)L)g$ \\
\hline
Sasakian & $(n-2)\left( 1-L\right) S+(n-1)(1+(n-2)L)g$ \\ \hline
Kenmotsu & $-\,(n-2)\left( 1+L\right) S-(n-1)(-1+(n-2)L)g$ \\ \hline
$(\varepsilon )$-Sasakian & $(n-2)\left( \varepsilon -L\right) S+\varepsilon
(n-1)(\varepsilon +(n-2)L)g$ \\ \hline
para-Sasakian & $-\,(n-2)\left( 1+L\right) S-(n-1)(-1+(n-2)L)g$ \\ \hline
$(\varepsilon )$-para-Sasakian & $-\,(n-2)\left( \varepsilon +L\right)
S-\varepsilon (n-1)(-\varepsilon +(n-2)L)g$ \\ \hline
\end{tabular}%
\]
\end{cor}

\begin{cor}
Let $M$ be an $n$-dimensional $({\cal V},S)$-pseudosymmetric $\left(
N(k),\xi \right) $-semi-Riemannian manifold. Then we have the following
table\/{\rm :}~%
\[
\begin{tabular}{|l|l|l|}
\hline
${\boldmath M}$ & ${\boldmath L=}$ & ${\boldmath S=}$ \\ \hline
$N(k)$-contact metric & $k-\dfrac{r}{n(n-1)}$ & $k(n-1)g$ \\ \hline
Sasakian & $1-\dfrac{r}{n(n-1)}$ & $(n-1)g$ \\ \hline
Kenmotsu & $-\,1-\dfrac{r}{n(n-1)}$ & $-\,(n-1)g$ \\ \hline
$(\varepsilon )$-Sasakian & $\varepsilon -\dfrac{r}{n(n-1)}$ & $\varepsilon
(n-1)g$ \\ \hline
para-Sasakian & $-\,1-\dfrac{r}{n(n-1)}$ & $-\,(n-1)g$ \\ \hline
$(\varepsilon )$-para-Sasakian & $-\,\varepsilon -\dfrac{r}{n(n-1)}$ & $%
-\,\varepsilon (n-1)g$ \\ \hline
\end{tabular}%
\]
\end{cor}

\begin{cor}
Let $M$ be an $n$-dimensional $({\cal P}_{\ast },S)$-pseudosymmetric $\left(
N(k),\xi \right) $-semi-Riemannian manifold. Then we have the following
table\/{\rm :}~%
\[
\begin{tabular}{|l|l|l|}
\hline
${\boldmath M}$ & ${\boldmath L=}$ & ${\boldmath S=}$ \\ \hline
$N(k)$-contact metric & $\left( k-\dfrac{r}{n(n-1)}\right) a_{0}-\dfrac{r}{n}%
a_{1}$ & $k(n-1)g$ \\ \hline
Sasakian & $\left( 1-\dfrac{r}{n(n-1)}\right) a_{0}-\dfrac{r}{n}a_{1}$ & $%
(n-1)g$ \\ \hline
Kenmotsu & $\left( -1-\dfrac{r}{n(n-1)}\right) a_{0}-\dfrac{r}{n}a_{1}$ & $%
-\,(n-1)g$ \\ \hline
$(\varepsilon )$-Sasakian & $\left( \varepsilon -\dfrac{r}{n(n-1)}\right)
a_{0}-\dfrac{r}{n}a_{1}$ & $\varepsilon (n-1)g$ \\ \hline
para-Sasakian & $\left( -1-\dfrac{r}{n(n-1)}\right) a_{0}-\dfrac{r}{n}a_{1}$
& $-\,(n-1)g$ \\ \hline
$(\varepsilon )$-para-Sasakian & $\left( -\varepsilon -\dfrac{r}{n(n-1)}%
\right) a_{0}-\dfrac{r}{n}a_{1}$ & $-\,\varepsilon (n-1)g$ \\ \hline
\end{tabular}%
\]
\end{cor}

\begin{cor}
Let $M$ be an $n$-dimensional $({\cal P},S)$-pseudosymmetric $\left(
N(k),\xi \right) $-semi-Riemannian manifold. Then we have the following
table\/{\rm :}~%
\[
\begin{tabular}{|l|l|l|}
\hline
${\boldmath M}$ & ${\boldmath L=}$ & ${\boldmath S=}$ \\ \hline
$N(k)$-contact metric & $k$ & $k(n-1)g$ \\ \hline
Sasakian & $1$ & $(n-1)g$ \\ \hline
Kenmotsu & $-\,1$ & $-\,(n-1)g$ \\ \hline
$(\varepsilon )$-Sasakian & $\varepsilon $ & $\varepsilon (n-1)g$ \\ \hline
para-Sasakian & $-\,1$ & $-\,(n-1)g$ \\ \hline
$(\varepsilon )$-para-Sasakian & $-\,\varepsilon $ & $-\,\varepsilon (n-1)g$
\\ \hline
\end{tabular}%
\]
\end{cor}

\begin{cor}
Let $M$ be an $n$-dimensional $({\cal M},S)$-pseudosymmetric $\left(
N(k),\xi \right) $-semi-Riemannian manifold. Then we have the following
table\/{\rm :}~%
\[
\begin{tabular}{|l|l|}
\hline
${\boldmath M}$ & ${\boldmath S^{2}=}$ \\ \hline
$N(k)$-contact metric & $2(n-1)\left( k-L\right) S-k(n-1)^{2}\left(
k-2L\right) g$ \\ \hline
Sasakian & $2(n-1)\left( 1-L\right) S-(n-1)^{2}\left( 1-2L\right) g$ \\
\hline
Kenmotsu & $-\,2(n-1)\left( 1+L\right) S-(n-1)^{2}\left( 1+2L\right) g$ \\
\hline
$(\varepsilon )$-Sasakian & $2(n-1)\left( \varepsilon -L\right)
S-\varepsilon (n-1)^{2}\left( \varepsilon -2L\right) g$ \\ \hline
para-Sasakian & $-\,2(n-1)\left( 1+L\right) S-(n-1)^{2}\left( 1+2L\right) g$
\\ \hline
$(\varepsilon )$-para-Sasakian & $-\,2(n-1)\left( \varepsilon +L\right)
S-\varepsilon (n-1)^{2}\left( \varepsilon +2L\right) g$ \\ \hline
\end{tabular}%
\]
\end{cor}

\begin{cor}
Let $M$ be an $n$-dimensional $({\cal W}_{0},S)$-pseudosymmetric $\left(
N(k),\xi \right) $-semi-Riemannian manifold. Then we have the following
table\/{\rm :}~%
\[
\begin{tabular}{|l|l|}
\hline
${\boldmath M}$ & ${\boldmath S^{2}=}$ \\ \hline
$N(k)$-contact metric & $(n-1)(2k-L)S+k(n-1)^{2}\left( L-k\right) g$ \\
\hline
Sasakian & $(n-1)(2-L)S+(n-1)^{2}\left( L-1\right) g$ \\ \hline
Kenmotsu & $-\,(n-1)(2+L)S-(n-1)^{2}\left( L+1\right) g$ \\ \hline
$(\varepsilon )$-Sasakian & $(n-1)(2\varepsilon -L)S+\varepsilon
(n-1)^{2}\left( L-\varepsilon \right) g$ \\ \hline
para-Sasakian & $-\,(n-1)(2+L)S-(n-1)^{2}\left( L+1\right) g$ \\ \hline
$(\varepsilon )$-para-Sasakian & $-\,(n-1)(2\varepsilon +L)S-\varepsilon
(n-1)^{2}\left( L+\varepsilon \right) g$ \\ \hline
\end{tabular}%
\]
\end{cor}

\begin{cor}
Let $M$ be an $n$-dimensional $({\cal W}_{0}^{\ast },S)$-pseudosymmetric $%
\left( N(k),\xi \right) $-semi-Riemannian manifold. Then we have the
following table\/{\rm :}~%
\[
\begin{tabular}{|l|l|}
\hline
${\boldmath M}$ & ${\boldmath S^{2}=}$ \\ \hline
$N(k)$-contact metric & $L(n-1)S+k(n-1)^{2}(k-L)g$ \\ \hline
Sasakian & $L(n-1)S+(n-1)^{2}(1-L)g$ \\ \hline
Kenmotsu & $L(n-1)S+(n-1)^{2}(1+L)g$ \\ \hline
$(\varepsilon )$-Sasakian & $L(n-1)S+\varepsilon (n-1)^{2}(\varepsilon -L)g$
\\ \hline
para-Sasakian & $L(n-1)S+(n-1)^{2}(1+L)g$ \\ \hline
$(\varepsilon )$-para-Sasakian & $L(n-1)S+\varepsilon (n-1)^{2}(\varepsilon
+L)g$ \\ \hline
\end{tabular}%
\]
\end{cor}

\begin{cor}
Let $M$ be an $n$-dimensional $({\cal W}_{1},S)$-pseudosymmetric $\left(
N(k),\xi \right) $-semi-Riemannian manifold. Then we have the following
table\/{\rm :}~%
\[
\begin{tabular}{|l|l|l|}
\hline
${\boldmath M}$ & ${\boldmath L=}$ & ${\boldmath S=}$ \\ \hline
$N(k)$-contact metric & $k$ & $k(n-1)g$ \\ \hline
Sasakian & $1$ & $(n-1)g$ \\ \hline
Kenmotsu & $-\,1$ & $-\,(n-1)g$ \\ \hline
$(\varepsilon )$-Sasakian & $\varepsilon $ & $\varepsilon (n-1)g$ \\ \hline
para-Sasakian & $-\,1$ & $-\,(n-1)g$ \\ \hline
$(\varepsilon )$-para-Sasakian & $-\,\varepsilon $ & $-\,\varepsilon (n-1)g$
\\ \hline
\end{tabular}%
\]
\end{cor}

\begin{cor}
Let $M$ be an $n$-dimensional $({\cal W}_{1}^{\ast },S)$-pseudosymmetric $%
\left( N(k),\xi \right) $-semi-Riemannian manifold. Then we have the
following table\/{\rm :}~%
\[
\begin{tabular}{|l|l|l|}
\hline
${\boldmath M}$ & ${\boldmath L=}$ & ${\boldmath S=}$ \\ \hline
$N(k)$-contact metric & $k$ & $k(n-1)g$ \\ \hline
Sasakian & $1$ & $(n-1)g$ \\ \hline
Kenmotsu & $-\,1$ & $-\,(n-1)g$ \\ \hline
$(\varepsilon )$-Sasakian & $\varepsilon $ & $\varepsilon (n-1)g$ \\ \hline
para-Sasakian & $-\,1$ & $-\,(n-1)g$ \\ \hline
$(\varepsilon )$-para-Sasakian & $-\,\varepsilon $ & $-\,\varepsilon (n-1)g$
\\ \hline
\end{tabular}%
\]
\end{cor}

\begin{cor}
Let $M$ be an $n$-dimensional $({\cal W}_{2},S)$-pseudosymmetric $\left(
N(k),\xi \right) $-semi-Riemannian manifold. Then we have the following
table\/{\rm :}~%
\[
\begin{tabular}{|l|l|}
\hline
${\boldmath M}$ & ${\boldmath S^{2}=}$ \\ \hline
$N(k)$-contact metric & $(n-1)\left( k-L\right) S+k(n-1)^{2}Lg$ \\ \hline
Sasakian & $(n-1)\left( 1-L\right) S+(n-1)^{2}Lg$ \\ \hline
Kenmotsu & $-\,(n-1)\left( 1+L\right) S-(n-1)^{2}Lg$ \\ \hline
$(\varepsilon )$-Sasakian & $(n-1)\left( \varepsilon -L\right) S+\varepsilon
(n-1)^{2}Lg$ \\ \hline
para-Sasakian & $-\,(n-1)\left( 1+L\right) S-(n-1)^{2}Lg$ \\ \hline
$(\varepsilon )$-para-Sasakian & $-\,(n-1)\left( \varepsilon +L\right)
S-\varepsilon (n-1)^{2}Lg$ \\ \hline
\end{tabular}%
\]
\end{cor}

\begin{cor}
Let $M$ be an $n$-dimensional $({\cal W}_{3},S)$-pseudosymmetric $\left(
N(k),\xi \right) $-semi-Riemannian manifold. Then we have the following
table\/{\rm :}~%
\[
\begin{tabular}{|l|l|l|}
\hline
${\boldmath M}$ & ${\boldmath L=}$ & ${\boldmath S=}$ \\ \hline
$N(k)$-contact metric & $2k$ & $k(n-1)g$ \\ \hline
Sasakian & $2$ & $(n-1)g$ \\ \hline
Kenmotsu & $-\,2$ & $-\,(n-1)g$ \\ \hline
$(\varepsilon )$-Sasakian & $2\varepsilon $ & $\varepsilon (n-1)g$ \\ \hline
para-Sasakian & $-\,2$ & $-\,(n-1)g$ \\ \hline
$(\varepsilon )$-para-Sasakian & $-\,2\varepsilon $ & $-\,\varepsilon (n-1)g$
\\ \hline
\end{tabular}%
\]
\end{cor}

\begin{cor}
Let $M$ be an $n$-dimensional $({\cal W}_{4},S)$-pseudosymmetric $\left(
N(k),\xi \right) $-semi-Riemannian manifold. Then we have the following
table\/{\rm :}~%
\[
\begin{tabular}{|l|l|}
\hline
${\boldmath M}$ & ${\boldmath S^{2}=}$ \\ \hline
$N(k)$-contact metric & $(n-1)\left( k-L\right) S+k(n-1)^{2}\left(
L-k\right) g+\varepsilon k^{2}(n-1)^{2}\eta \otimes \eta $ \\ \hline
Sasakian & $(n-1)\left( 1-L\right) S+(n-1)^{2}\left( L-1\right)
g+\varepsilon (n-1)^{2}\eta \otimes \eta $ \\ \hline
Kenmotsu & $-\,(n-1)\left( 1+L\right) S-(n-1)^{2}\left( L+1\right)
g+\varepsilon (n-1)^{2}\eta \otimes \eta $ \\ \hline
$(\varepsilon )$-Sasakian & $(n-1)\left( \varepsilon -L\right) S+\varepsilon
(n-1)^{2}\left( L-\varepsilon \right) g+\varepsilon (n-1)^{2}\eta \otimes
\eta $ \\ \hline
para-Sasakian & $-\,(n-1)\left( 1+L\right) S-(n-1)^{2}\left( L+1\right)
g+\varepsilon (n-1)^{2}\eta \otimes \eta $ \\ \hline
$(\varepsilon )$-para-Sasakian & $-\,(n-1)\left( \varepsilon +L\right)
S-\varepsilon (n-1)^{2}\left( L+\varepsilon \right) g+\varepsilon
(n-1)^{2}\eta \otimes \eta $ \\ \hline
\end{tabular}%
\]
\end{cor}

\begin{cor}
Let $M$ be an $n$-dimensional $({\cal W}_{5},S)$-pseudosymmetric $\left(
N(k),\xi \right) $-semi-Riemannian manifold. Then we have the following
table\/{\rm :}~%
\[
\begin{tabular}{|l|l|}
\hline
${\boldmath M}$ & ${\boldmath S^{2}=}$ \\ \hline
$N(k)$-contact metric & $(n-1)\left( 2k-L\right) S+k(n-1)^{2}\left(
L-k\right) g$ \\ \hline
Sasakian & $(n-1)\left( 2-L\right) S+(n-1)^{2}\left( L-1\right) g$ \\ \hline
Kenmotsu & $-\,(n-1)\left( 2+L\right) S-(n-1)^{2}\left( L+1\right) g$ \\
\hline
$(\varepsilon )$-Sasakian & $(n-1)\left( 2\varepsilon -L\right)
S+\varepsilon (n-1)^{2}\left( L-\varepsilon \right) g$ \\ \hline
para-Sasakian & $-\,(n-1)\left( 2+L\right) S-(n-1)^{2}\left( L+1\right) g$
\\ \hline
$(\varepsilon )$-para-Sasakian & $-\,(n-1)\left( 2\varepsilon +L\right)
S+k(n-1)^{2}\left( L+\varepsilon \right) g$ \\ \hline
\end{tabular}%
\]
\end{cor}

\begin{cor}
Let $M$ be an $n$-dimensional $({\cal W}_{6},S)$-pseudosymmetric $\left(
N(k),\xi \right) $-semi-Riemannian manifold. Then we have the following
table\/{\rm :}~%
\[
\begin{tabular}{|l|l|}
\hline
${\boldmath M}$ & {\bf Result} \\ \hline
$N(k)$-contact metric & $\left( 2k-L\right) S+k(n-1)\left( L-k\right)
g+k^{2}(n-1)\eta \otimes \eta =0$ \\ \hline
Sasakian & $\left( 2-L\right) S+(n-1)\left( L-1\right) g+(n-1)\eta \otimes
\eta =0$ \\ \hline
Kenmotsu & $-\,\left( 2+L\right) S-(n-1)\left( L+1\right) g+(n-1)\eta
\otimes \eta =0$ \\ \hline
$(\varepsilon )$-Sasakian & $\left( 2\varepsilon -L\right) S+\varepsilon
(n-1)\left( L-\varepsilon \right) g+\varepsilon (n-1)\eta \otimes \eta =0$
\\ \hline
para-Sasakian & $-\,\left( 2+L\right) S-(n-1)\left( L+1\right) g+(n-1)\eta
\otimes \eta =0$ \\ \hline
$(\varepsilon )$-para-Sasakian & $-\,\left( 2\varepsilon +L\right)
S-\varepsilon (n-1)\left( L+\varepsilon \right) g+\varepsilon (n-1)\eta
\otimes \eta =0$ \\ \hline
\end{tabular}%
\]
\end{cor}

\begin{cor}
Let $M$ be an $n$-dimensional $({\cal W}_{7},S)$-pseudosymmetric $\left(
N(k),\xi \right) $-semi-Riemannian manifold. Then we have the following
table\/{\rm :}~%
\[
\begin{tabular}{|l|l|l|}
\hline
${\boldmath M}$ & ${\boldmath L=}$ & ${\boldmath S=}$ \\ \hline
$N(k)$-contact metric & $2k$ & $k(n-1)g$ \\ \hline
Sasakian & $2$ & $(n-1)g$ \\ \hline
Kenmotsu & $-\,2$ & $-\,(n-1)g$ \\ \hline
$(\varepsilon )$-Sasakian & $2\varepsilon $ & $\varepsilon (n-1)g$ \\ \hline
para-Sasakian & $-\,2$ & $-\,(n-1)g$ \\ \hline
$(\varepsilon )$-para-Sasakian & $-\,2\varepsilon $ & $-\,\varepsilon (n-1)g$
\\ \hline
\end{tabular}%
\]
\end{cor}

\begin{cor}
Let $M$ be an $n$-dimensional $({\cal W}_{8},S)$-pseudosymmetric $\left(
N(k),\xi \right) $-semi-Riemannian manifold. Then we have the following
table\/{\rm :}~%
\[
\begin{tabular}{|l|l|}
\hline
${\boldmath M}$ & {\bf Result} \\ \hline
$N(k)$-contact metric & $\left( 2k-L\right) S+k(n-1)\left( L-k\right)
g-k^{2}(n-1)\eta \otimes \eta =0$ \\ \hline
Sasakian & $\left( 2-L\right) S+(n-1)\left( L-1\right) g-(n-1)\eta \otimes
\eta =0$ \\ \hline
Kenmotsu & $-\,\left( 2+L\right) S-(n-1)\left( L+1\right) g-(n-1)\eta
\otimes \eta =0$ \\ \hline
$(\varepsilon )$-Sasakian & $\left( 2\varepsilon -L\right) S+\varepsilon
(n-1)\left( L-\varepsilon \right) g-\varepsilon (n-1)\eta \otimes \eta =0$
\\ \hline
para-Sasakian & $-\,\left( 2+L\right) S-(n-1)\left( L+1\right) g-(n-1)\eta
\otimes \eta =0$ \\ \hline
$(\varepsilon )$-para-Sasakian & $-\,\left( 2\varepsilon +L\right)
S-\varepsilon (n-1)\left( L+\varepsilon \right) g-\varepsilon (n-1)\eta
\otimes \eta =0$ \\ \hline
\end{tabular}%
\]
\end{cor}

\begin{cor}
Let $M$ be an $n$-dimensional $({\cal W}_{9},S)$-pseudosymmetric $\left(
N(k),\xi \right) $-semi-Riemannian manifold. Then we have the following
table\/{\rm :}~%
\[
\begin{tabular}{|l|l|}
\hline
${\boldmath M}$ & {\bf Result} \\ \hline
$N(k)$-contact metric & $\left( L-k\right) S-k(n-1)Lg+k^{2}(n-1)\eta \otimes
\eta =0$ \\ \hline
Sasakian & $\left( L-1\right) S-(n-1)Lg+(n-1)\eta \otimes \eta =0$ \\ \hline
Kenmotsu & $\left( L+1\right) S+(n-1)Lg+(n-1)\eta \otimes \eta =0$ \\ \hline
$(\varepsilon )$-Sasakian & $\left( L-\varepsilon \right) S-\varepsilon
(n-1)Lg+\varepsilon (n-1)\eta \otimes \eta =0$ \\ \hline
para-Sasakian & $\left( L+1\right) S+(n-1)Lg+(n-1)\eta \otimes \eta =0$ \\
\hline
$(\varepsilon )$-para-Sasakian & $\left( L+\varepsilon \right) S+\varepsilon
(n-1)Lg+\varepsilon (n-1)\eta \otimes \eta =0$ \\ \hline
\end{tabular}%
\]
\end{cor}

\begin{th}
\label{th-T-ric-pseudo} Let $M$ be an $n$-dimensional $(R,S_{{\cal T}_{a}})$%
-pseudosymmetric $(N(k),\xi )$-semi-Riemannian manifold such that
\[
a_{0}+na_{1}+a_{2}+a_{3}+a_{5}+a_{6}\not=0.
\]%
Then $M$ is either Einstein manifold, that is,
\[
S=k(n-1)g
\]%
or $L=k$ holds on $M$. Consequently, we have the following table\/{\rm :}
\[
\begin{tabular}{|l|l|l|l|}
\hline
{\bf Manifold} & {\bf Condition} & ${\boldmath S=}$ & ${\boldmath L=}$ \\
\hline
$N(k)$-contact metric & $R\cdot S_{{\cal T}_{a}}=LQ(g,S_{{\cal T}_{a}})$ & $%
k(n-1)g$ & $k$ \\ \hline
Sasakian & $R\cdot S_{{\cal T}_{a}}=LQ(g,S_{{\cal T}_{a}})$ & $(n-1)g$ & $1$
\\ \hline
Kenmotsu & $R\cdot S_{{\cal T}_{a}}=LQ(g,S_{{\cal T}_{a}})$ & $-\,(n-1)g$ & $%
-\,1$ \\ \hline
$(\varepsilon )$-Sasakian & $R\cdot S_{{\cal T}_{a}}=LQ(g,S_{{\cal T}_{a}})$
& $\varepsilon (n-1)g$ & $\varepsilon $ \\ \hline
para-Sasakian & $R\cdot S_{{\cal T}_{a}}=LQ(g,S_{{\cal T}_{a}})$ & $%
-\,(n-1)g $ & $-\,1$ \\ \hline
$(\varepsilon )$-para-Sasakian & $R\cdot S_{{\cal T}_{a}}=LQ(g,S_{{\cal T}%
_{a}})$ & $-\,\varepsilon (n-1)g$ & $-\,\varepsilon $ \\ \hline
\end{tabular}%
\ \
\]
\end{th}

\begin{rem-new}
The conclusions of Theorem {\rm \ref{th-T-ric-pseudo}} remain true if $S_{%
{\cal T}_{a}}$ is replaced by $S$.
\end{rem-new}

\begin{cor}
{\rm (\cite{Ozgur-06}, \cite{Hong-Ozgur-Tripathi-06})} If an $n$-dimensional
Kenmotsu manifold $M$ is Ricci-pseudosymmetric then either $M$ is an
Einstein manifold with the scalar curvature $r=n(1-n)$ or $L=-1$ holds on $M$%
.
\end{cor}

\section{$({\cal T}_{\!a},S_{{\cal T}_{\!b}},S^{\ell })$-pseudosymmetry\label%
{sect-TSSP}}

In this section, we determine the result for an $n$-dimensional $\left(
N(k),\xi \right) $-semi-Riemannian manifold satisfy ${\cal T}_{\!a}\cdot S_{%
{\cal T}_{b}}=LQ(S^{\ell },S_{{\cal T}_{\!b}})$.

\begin{defn-new}
A semi-Riemannian manifold $M$ is called $({\cal T}_{\!a},S_{{\cal T}%
_{\!b}},S^{\ell })$-pseudosymmetric if
\[
{\cal T}_{\!a}\cdot S_{{\cal T}_{b}}=LQ(S^{\ell },S_{{\cal T}_{\!b}}),
\]%
where $L$ is some smooth function defined on $M$. In particular, $M$ is said
to be $(R,S_{{\cal T}_{\!a}},S^{\ell })$-pseudosymmetric if
\[
R\cdot S_{{\cal T}_{a}}=LQ(S^{\ell },S_{{\cal T}_{\!a}}).
\]
\end{defn-new}

\begin{th}
Let $M$ be an $n$-dimensional $({\cal T}_{\!a},S_{{\cal T}_{\!b}},S^{\ell })$%
-pseudosymmetric $(N(k),\xi )$-semi-Riemannian manifold. Then
\begin{eqnarray*}
&&\varepsilon a_{5}(b_{0}+nb_{1}+b_{2}+b_{3}+b_{5}+b_{6})S^{2}(Y,U) \\
&&+\ \left\{ \varepsilon (b_{0}+nb_{1}+b_{2}+b_{3}+b_{5}+b_{6})\times \right.
\\
&&(-ka_{0}+k(n-1)a_{1}+k(n-1)a_{2}-a_{7}r) \\
&&\qquad \left. +\ \varepsilon (a_{1}+a_{5})(b_{4}r+(n-1)b_{7}r)\right\}
S(Y,U) \\
&&+\ \left\{ \varepsilon k(n-1)(a_{2}+a_{4})(b_{4}r+(n-1)b_{7}r)\right. \\
&&\qquad +\ \varepsilon k(n-1)(b_{0}+nb_{1}+b_{2}+b_{3}+b_{5}+b_{6})\times \\
&&\left. (ka_{0}+k(n-1)a_{4}+a_{7}r)\right\} g(Y,U) \\
&&+\ k(n-1)(a_{1}+a_{2}+2a_{3}+a_{4}+a_{5}+2a_{6})\times \\
&&\left\{ b_{4}r+(n-1)b_{7}r\right. \\
&&\qquad \left. +\ k(n-1)(b_{0}+nb_{1}+b_{2}+b_{3}+b_{5}+b_{6})\right\} \eta
(Y)\eta (U) \\
&=&L\varepsilon ((b_{0}+nb_{1}+b_{2}+b_{3}+b_{5}+b_{6})\times \\
&&(k(n-1)S^{\ell }(Y,U)-k^{\ell }(n-1)^{\ell }S(Y,U)) \\
&&+\,(b_{4}+(n-1)b_{7})r(k^{\ell }(n-1)^{\ell }g(Y,U)-S^{\ell }(Y,U))).
\end{eqnarray*}%
In particular, if $M$ be an $n$-dimensional $({\cal T}_{\!a},S_{{\cal T}%
_{\!a}},S^{\ell })$-pseudosymmetric $(N(k),\xi )$-semi-Riemannian manifold.
Then
\begin{eqnarray*}
&&\varepsilon a_{5}(a_{0}+na_{1}+a_{2}+a_{3}+a_{5}+a_{6})S^{2}(Y,U) \\
&&+\ \left\{ \varepsilon (a_{0}+na_{1}+a_{2}+a_{3}+a_{5}+a_{6})\times \right.
\\
&&(-ka_{0}+k(n-1)a_{1}+k(n-1)a_{2}-a_{7}r) \\
&&\qquad \left. +\ \varepsilon (a_{1}+a_{5})(a_{4}r+(n-1)a_{7}r)\right\}
S(Y,U) \\
&&+\ \left\{ \varepsilon k(n-1)(a_{2}+a_{4})(a_{4}r+(n-1)a_{7}r)\right. \\
&&\qquad +\ \varepsilon k(n-1)(a_{0}+na_{1}+a_{2}+a_{3}+a_{5}+a_{6})\times \\
&&\left. (ka_{0}+k(n-1)a_{4}+a_{7}r)\right\} g(Y,U) \\
&&+\ k(n-1)(a_{1}+a_{2}+2a_{3}+a_{4}+a_{5}+2a_{6}) \\
&&\left\{ a_{4}r+(n-1)a_{7}r\right. \\
&&\qquad \left. +\ \varepsilon
k(n-1)(a_{0}+na_{1}+a_{2}+a_{3}+a_{5}+a_{6})\right\} \eta (Y)\eta (U) \\
&=&L\varepsilon ((a_{0}+na_{1}+a_{2}+a_{3}+a_{5}+a_{6})\times \\
&&(k(n-1)S^{\ell }(Y,U)-k^{\ell }(n-1)^{\ell }S(Y,U)) \\
&&+\,(a_{4}+(n-1)a_{7})r(k^{\ell }(n-1)^{\ell }g(Y,U)-S^{\ell }(Y,U))).
\end{eqnarray*}
\end{th}

\noindent {\bf Proof.} Let $M$ be an $n$-dimensional $({\cal T}_{\!a},S_{%
{\cal T}_{\!b}},S^{\ell })$-pseudosymmetric $\left( N(k),\xi \right) $%
-semi-Riemannian manifold. Then
\begin{equation}
{\cal T}_{\!a}(X,Y)\cdot S_{{\cal T}_{\!b}}(U,V)=LQ(S^{\ell },S_{{\cal T}%
_{\!b}})(U,V;X,Y).  \label{eq-Ric-pseudo-1}
\end{equation}%
Taking $X=\xi =V$ in (\ref{eq-Ric-pseudo-1}), we have
\[
{\cal T}_{\!a}(\xi ,Y)\cdot S_{{\cal T}_{\!b}}(U,\xi )=LQ(S^{\ell },S_{{\cal %
T}_{\!b}})(U,\xi ;\xi ,Y),
\]%
which gives
\begin{eqnarray}
&&S_{{\cal T}_{\!b}}({\cal T}_{\!a}(\xi ,Y)U,\xi )+S_{{\cal T}_{\!b}}(U,%
{\cal T}_{\!a}(\xi ,Y)\xi )  \nonumber \\
&=&L\left( S_{{\cal T}_{\!b}}((\xi \wedge _{S^{\ell }}Y)U,\xi )+S_{{\cal T}%
_{\!b}}(U,(\xi \wedge _{S^{\ell }}Y)\xi )\right) .  \label{eq-T-S-111}
\end{eqnarray}%
Using (\ref{eq-cond}), (\ref{eq-Sp-QX-xi}), (\ref{eq-xi-X-xi}), (\ref%
{eq-xi-Y-Z}), (\ref{eq-ric-T1}) and (\ref{eq-ric-T2}) in (\ref{eq-T-S-111}),
we get the result. $\blacksquare $

For $\ell =1$, we have the following result.

\begin{cor}
Let $M$ be an $n$-dimensional $({\cal T}_{\!a},S_{{\cal T}_{\!b}},S)$%
-pseudosymmetric $(N(k),\xi )$-semi-Riemannian manifold. Then
\begin{eqnarray*}
&&\varepsilon a_{5}(b_{0}+nb_{1}+b_{2}+b_{3}+b_{5}+b_{6})S^{2}(Y,U) \\
&&+\ \left\{ \varepsilon (b_{0}+nb_{1}+b_{2}+b_{3}+b_{5}+b_{6})\times \right.
\\
&&(-ka_{0}+k(n-1)a_{1}+k(n-1)a_{2}-a_{7}r) \\
&&\qquad \left. +\ \varepsilon (a_{1}+a_{5})(b_{4}r+(n-1)b_{7}r)\right\}
S(Y,U) \\
&&+\ \left\{ \varepsilon k(n-1)(a_{2}+a_{4})(b_{4}r+(n-1)b_{7}r)\right. \\
&&\qquad +\ \varepsilon k(n-1)(b_{0}+nb_{1}+b_{2}+b_{3}+b_{5}+b_{6})\times \\
&&\left. (ka_{0}+k(n-1)a_{4}+a_{7}r)\right\} g(Y,U) \\
&&+\ k(n-1)(a_{1}+a_{2}+2a_{3}+a_{4}+a_{5}+2a_{6})\times \\
&&\left\{ b_{4}r+(n-1)b_{7}r\right. \\
&&\qquad \left. +\ k(n-1)(b_{0}+nb_{1}+b_{2}+b_{3}+b_{5}+b_{6})\right\} \eta
(Y)\eta (U) \\
&=&L\varepsilon (b_{4}+(n-1)b_{7})r(k(n-1)g(Y,U)-S(Y,U)).
\end{eqnarray*}%
In particular, if $M$ be an $n$-dimensional $({\cal T}_{\!a},S_{{\cal T}%
_{\!a}},S)$-pseudosymmetric $(N(k),\xi )$-semi-Riemannian manifold. Then
\begin{eqnarray*}
&&\varepsilon a_{5}(a_{0}+na_{1}+a_{2}+a_{3}+a_{5}+a_{6})S^{2}(Y,U) \\
&&+\ \left\{ \varepsilon (a_{0}+na_{1}+a_{2}+a_{3}+a_{5}+a_{6})\times \right.
\\
&&(-ka_{0}+k(n-1)a_{1}+k(n-1)a_{2}-a_{7}r) \\
&&\qquad \left. +\ \varepsilon (a_{1}+a_{5})(a_{4}r+(n-1)a_{7}r)\right\}
S(Y,U) \\
&&+\ \left\{ \varepsilon k(n-1)(a_{2}+a_{4})(a_{4}r+(n-1)a_{7}r)\right. \\
&&\qquad +\ \varepsilon k(n-1)(a_{0}+na_{1}+a_{2}+a_{3}+a_{5}+a_{6})\times \\
&&\left. (ka_{0}+k(n-1)a_{4}+a_{7}r)\right\} g(Y,U) \\
&&+\ k(n-1)(a_{1}+a_{2}+2a_{3}+a_{4}+a_{5}+2a_{6}) \\
&&\left\{ a_{4}r+(n-1)a_{7}r\right. \\
&&\qquad \left. +\ \varepsilon
k(n-1)(a_{0}+na_{1}+a_{2}+a_{3}+a_{5}+a_{6})\right\} \eta (Y)\eta (U) \\
&=&L\varepsilon (a_{4}+(n-1)a_{7})r(k(n-1)g(Y,U)-S(Y,U)).
\end{eqnarray*}
\end{cor}

\begin{th}
\label{GCT-rssss} Let $M$ be an $n$-dimensional $\left( {\cal T}%
_{\!a},S,S^{\ell }\right) $-pseudosymmetric $\left( N(k),\xi \right)$%
-semi-Riemannian manifold. Then
\begin{eqnarray*}
&&\varepsilon a_{5}\,S^{2}(Y,U)-E\,S(Y,U)-Fg(Y,U)-G\eta (Y)\eta \left(
U\right) \\
&=&\varepsilon L(k(n-1)S^{\ell }(Y,U)-k^{\ell }(n-1)^{\ell }S(Y,U)),
\end{eqnarray*}%
where
\[
E=\varepsilon \,(ka_{0}+a_{7}r-k(n-1)a_{1}-k(n-1)a_{2}),
\]%
\[
F=-\,\varepsilon k(n-1)(ka_{0}+k(n-1)a_{4}+a_{7}r),
\]%
\[
G=-\,k^{2}(n-1)^{2}(a_{1}+a_{2}+2a_{3}+a_{4}+a_{5}+2a_{6}).
\]
\end{th}

In view of Theorem~\ref{GCT-rssss}, we have the following

\begin{cor}
Let $M$ be an $n$-dimensional $\left( R,S,S^{\ell }\right) $-pseudosymmetric
$\left( N(k),\xi \right) $-semi-Riemannian manifold. Then we have the
following table\/{\rm :}~%
\[
\begin{tabular}{|l|l|}
\hline
${\boldmath M}$ & ${\boldmath LS}^{\ell }{=}$ \\ \hline
$N(k)$-contact metric & $\left( Lk^{\ell -1}(n-1)^{\ell -1}-\dfrac{1}{n-1}%
\right) S+kg$ \\ \hline
Sasakian & $\left( L(n-1)^{\ell -1}-\dfrac{1}{n-1}\right) S+g$ \\ \hline
Kenmotsu & $\left( L(-1)^{\ell -1}(n-1)^{\ell -1}-\dfrac{1}{n-1}\right) S-g$
\\ \hline
$(\varepsilon )$-Sasakian & $\left( L(\varepsilon )^{\ell -1}(n-1)^{\ell -1}-%
\dfrac{1}{n-1}\right) S+\varepsilon g$ \\ \hline
para-Sasakian & $\left( L(-1)^{\ell -1}(n-1)^{\ell -1}-\dfrac{1}{n-1}\right)
S-g$ \\ \hline
$(\varepsilon )$-para-Sasakian & $\left( L(-\varepsilon )^{\ell
-1}(n-1)^{\ell -1}-\dfrac{1}{n-1}\right) S-\varepsilon g$ \\ \hline
\end{tabular}%
\]
\end{cor}

\begin{cor}
Let $M$ be an $n$-dimensional $\left( {\cal C}_{\ast },S,S^{\ell }\right) $%
-pseudosymmetric $\left( N(k),\xi \right) $-semi-Riemannian manifold. Then
we have the following table\/{\rm :}~
\[
\begin{tabular}{|l|l|}
\hline
${\boldmath M}$ & ${\boldmath L(n-1)S^{\ell }=}$ \\ \hline
$N(k)$-contact metric & $-\,\dfrac{a_{1}}{k}S^{2}-\left( \left( 1-\dfrac{r}{%
kn(n-1)}\right) a_{0}-\dfrac{2r}{kn}a_{1}-k^{\ell -1}(n-1)^{\ell }L\right) S$
\\
& $+(n-1)\left( \left( k-\dfrac{r}{n(n-1)}\right) a_{0}+\left( k(n-1)-\dfrac{%
2r}{n}\right) a_{1}\right) g$ \\ \hline
Sasakian & $-\,a_{1}S^{2}-\left( \left( 1-\dfrac{r}{n(n-1)}\right) a_{0}-%
\dfrac{2r}{n}a_{1}-(n-1)^{\ell }L\right) S$ \\
& $+(n-1)\left( \left( 1-\dfrac{r}{n(n-1)}\right) a_{0}+\left( (n-1)-\dfrac{%
2r}{n}\right) a_{1}\right) g$ \\ \hline
Kenmotsu & $a_{1}S^{2}-\left( \left( 1+\dfrac{r}{n(n-1)}\right) a_{0}+\dfrac{%
2r}{n}a_{1}-(-1)^{\ell -1}(n-1)^{\ell }L\right) S$ \\
& $-\,(n-1)\left( \left( 1+\dfrac{r}{n(n-1)}\right) a_{0}+\left( (n-1)+%
\dfrac{2r}{n}\right) a_{1}\right) g$ \\ \hline
$(\varepsilon )$-Sasakian & $-\,\varepsilon a_{1}S^{2}-\left( \left( 1-%
\dfrac{\varepsilon r}{n(n-1)}\right) a_{0}-\dfrac{2\varepsilon r}{n}%
a_{1}-(\varepsilon )^{\ell -1}(n-1)^{\ell }L\right) S$ \\
& $+(n-1)\left( \left( \varepsilon -\dfrac{r}{n(n-1)}\right) a_{0}+\left(
\varepsilon (n-1)-\dfrac{2r}{n}\right) a_{1}\right) g$ \\ \hline
para-Sasakian & $a_{1}S^{2}-\left( \left( 1+\dfrac{r}{n(n-1)}\right) a_{0}+%
\dfrac{2r}{n}a_{1}-(-1)^{\ell -1}(n-1)^{\ell }L\right) S$ \\
& $-\,(n-1)\left( \left( 1+\dfrac{r}{n(n-1)}\right) a_{0}+\left( (n-1)+%
\dfrac{2r}{n}\right) a_{1}\right) g$ \\ \hline
$(\varepsilon )$-para-Sasakian & $\varepsilon a_{1}S^{2}-\left( \left( 1+%
\dfrac{\varepsilon r}{n(n-1)}\right) a_{0}+\dfrac{2\varepsilon r}{n}%
a_{1}-(-\varepsilon )^{\ell -1}(n-1)^{\ell }L\right) S$ \\
& $-\,(n-1)\left( \left( \varepsilon +\dfrac{r}{n(n-1)}\right) a_{0}+\left(
\varepsilon (n-1)+\dfrac{2r}{n}\right) a_{1}\right) g$ \\ \hline
\end{tabular}%
\]
\end{cor}

\begin{cor}
Let $M$ be an $n$-dimensional $\left( {\cal C},S,S^{\ell }\right) $%
-pseudosymmetric $\left( N(k),\xi \right) $-semi-Riemannian manifold. Then
we have the following table\/{\rm :}~
\[
\begin{tabular}{|l|l|}
\hline
${\boldmath M}$ & ${\boldmath LS^{\ell }=}$ \\ \hline
$N(k)$-contact metric & $-\,\left( \dfrac{r}{k(n-1)^{2}(n-2)}+\dfrac{1}{n-1}%
-k^{\ell -1}(n-1)^{\ell -1}L\right) S$ \\
& $+\,\left( k-\dfrac{n-1}{n-2}+\dfrac{r}{(n-1)(n-2)}\right) g+\dfrac{1}{%
k(n-1)(n-2)}S^{2}$ \\ \hline
Sasakian & $-\,\left( \dfrac{r}{(n-1)^{2}(n-2)}+\dfrac{1}{n-1}-(n-1)^{\ell
-1}L\right) S$ \\
& $+\,\left( 1-\dfrac{n-1}{n-2}+\dfrac{r}{(n-1)(n-2)}\right) g+\dfrac{1}{%
(n-1)(n-2)}S^{2}$ \\ \hline
Kenmotsu & $-\,\left( -\dfrac{r}{(n-1)^{2}(n-2)}+\dfrac{1}{n-1}-(-1)^{\ell
-1}(n-1)^{\ell -1}L\right) S$ \\
& $+\,\left( -1-\dfrac{n-1}{n-2}+\dfrac{r}{(n-1)(n-2)}\right) g-\,\dfrac{1}{%
(n-1)(n-2)}S^{2}$ \\ \hline
$(\varepsilon )$-Sasakian & $-\,\left( \dfrac{\varepsilon r}{(n-1)^{2}(n-2)}+%
\dfrac{1}{n-1}-(\varepsilon )^{\ell -1}(n-1)^{\ell -1}L\right) S$ \\
& $+\,\left( \varepsilon -\dfrac{n-1}{n-2}+\dfrac{r}{(n-1)(n-2)}\right) g+%
\dfrac{\varepsilon }{(n-1)(n-2)}S^{2}$ \\ \hline
para-Sasakian & $-\,\left( -\dfrac{r}{(n-1)^{2}(n-2)}+\dfrac{1}{n-1}%
-(-1)^{\ell -1}(n-1)^{\ell -1}L\right) S$ \\
& $+\,\left( -1-\dfrac{n-1}{n-2}+\dfrac{r}{(n-1)(n-2)}\right) g-\,\dfrac{1}{%
(n-1)(n-2)}S^{2}$ \\ \hline
$(\varepsilon )$-para-Sasakian & $-\,\left( -\dfrac{\varepsilon r}{%
(n-1)^{2}(n-2)}+\dfrac{1}{n-1}-(-\varepsilon )^{\ell -1}(n-1)^{\ell
-1}L\right) S$ \\
& $+\,\left( -\varepsilon -\dfrac{n-1}{n-2}+\dfrac{r}{(n-1)(n-2)}\right) g-%
\dfrac{\varepsilon }{(n-1)(n-2)}S^{2}$ \\ \hline
\end{tabular}%
\]
\end{cor}

\begin{cor}
Let $M$ be an $n$-dimensional $\left( {\cal L},S,S^{\ell }\right) $%
-pseudosymmetric $\left( N(k),\xi \right) $-semi-Riemannian manifold. Then
we have the following table\/{\rm :}~
\[
\begin{tabular}{|l|l|}
\hline
${\boldmath M}$ & ${\boldmath LS^{\ell }=}$ \\ \hline
$N(k)$-contact metric & $\dfrac{1}{k(n-1)(n-2)}S^{2}+\left( k^{\ell
-1}(n-1)^{\ell -1}L-\dfrac{1}{n-1}\right) S-\dfrac{k}{n-2}g$ \\ \hline
Sasakian & $\dfrac{1}{(n-1)(n-2)}S^{2}+\left( (n-1)^{\ell -1}L-\dfrac{1}{n-1}%
\right) S-\dfrac{1}{n-2}g$ \\ \hline
Kenmotsu & $-\,\dfrac{1}{(n-1)(n-2)}S^{2}+\left( (-1)^{\ell -1}(n-1)^{\ell
-1}L-\dfrac{1}{n-1}\right) S+\dfrac{1}{n-2}g$ \\ \hline
$(\varepsilon )$-Sasakian & $\dfrac{\varepsilon }{(n-1)(n-2)}S^{2}+\left(
(\varepsilon )^{\ell -1}(n-1)^{\ell -1}L-\dfrac{1}{n-1}\right) S-\dfrac{%
\varepsilon }{n-2}g$ \\ \hline
para-Sasakian & $-\,\dfrac{1}{(n-1)(n-2)}S^{2}+\left( (-1)^{\ell
-1}(n-1)^{\ell -1}L-\dfrac{1}{n-1}\right) S+\dfrac{1}{n-2}g$ \\ \hline
$(\varepsilon )$-para-Sasakian & $-\,\dfrac{\varepsilon }{(n-1)(n-2)}%
S^{2}+\left( (-\varepsilon )^{\ell -1}(n-1)^{\ell -1}L-\dfrac{1}{n-1}\right)
S+\dfrac{\varepsilon }{n-2}g$ \\ \hline
\end{tabular}%
\]
\end{cor}

\begin{cor}
Let $M$ be an $n$-dimensional $\left( {\cal V},S,S^{\ell }\right) $%
-pseudosymmetric $\left( N(k),\xi \right) $-semi-Riemannian manifold. Then
we have the following table\/{\rm :}~
\[
\begin{tabular}{|l|l|}
\hline
${\boldmath M}$ & ${\boldmath}LS^{\ell }=$ \\ \hline
$N(k)$-contact metric & $-\,\left( \dfrac{1}{n-1}-\dfrac{r}{kn(n-1)^{2}}%
-k^{\ell -1}(n-1)^{\ell -1}\right) S$ \\
& $+\left( k-\dfrac{r}{n(n-1)}\right) g$ \\ \hline
Sasakian & $-\,\left( \dfrac{1}{n-1}-\dfrac{r}{n(n-1)^{2}}-(n-1)^{\ell
-1}\right) S$ \\
& $+\left( 1-\dfrac{r}{n(n-1)}\right) g$ \\ \hline
Kenmotsu & $-\,\left( \dfrac{1}{n-1}+\dfrac{r}{n(n-1)^{2}}-(-1)^{\ell
-1}(n-1)^{\ell -1}\right) S$ \\
& $-\left( 1+\dfrac{r}{n(n-1)}\right) g$ \\ \hline
$(\varepsilon )$-Sasakian & $-\,\left( \dfrac{1}{n-1}-\dfrac{\varepsilon r}{%
n(n-1)^{2}}-(\varepsilon )^{\ell -1}(n-1)^{\ell -1}\right) S$ \\
& $+\left( \varepsilon -\dfrac{r}{n(n-1)}\right) g$ \\ \hline
para-Sasakian & $-\,\left( \dfrac{1}{n-1}+\dfrac{r}{n(n-1)^{2}}-(-1)^{\ell
-1}(n-1)^{\ell -1}\right) S$ \\
& $-\left( 1+\dfrac{r}{n(n-1)}\right) g$ \\ \hline
$(\varepsilon )$-para-Sasakian & $-\,\left( \dfrac{1}{n-1}+\dfrac{%
\varepsilon r}{n(n-1)^{2}}-(-\varepsilon )^{\ell -1}(n-1)^{\ell -1}\right) S$
\\
& $-\left( \varepsilon +\dfrac{r}{n(n-1)}\right) g$ \\ \hline
\end{tabular}%
\]
\end{cor}

\begin{cor}
Let $M$ be an $n$-dimensional $\left( {\cal P}_{\ast },S,S^{\ell }\right) $%
-pseudosymmetric $\left( N(k),\xi \right) $-semi-Riemannian manifold. Then
we have the following table\/{\rm :}~
\[
\begin{tabular}{|l|l|}
\hline
${\boldmath M}$ & ${\boldmath LS}^{\ell }{=}$ \\ \hline
$N(k)$-contact metric & $-\,\left( \left( \dfrac{1}{n-1}-\dfrac{r}{%
kn(n-1)^{2}}\right) a_{0}-\dfrac{r}{kn(n-1)}a_{1}-k^{\ell -1}(n-1)^{\ell
-1}\right) S$ \\
& $+\left( \left( k-\dfrac{r}{n(n-1)}\right) a_{0}-\dfrac{r}{n}a_{1}\right)
g $ \\ \hline
Sasakian & $-\,\left( \left( \dfrac{1}{n-1}-\dfrac{r}{n(n-1)^{2}}\right)
a_{0}-\dfrac{r}{n(n-1)}a_{1}-(n-1)^{\ell -1}\right) S$ \\
& $+\left( \left( 1-\dfrac{r}{n(n-1)}\right) a_{0}-\dfrac{r}{n}a_{1}\right) $
\\ \hline
Kenmotsu & $-\,\left( \left( \dfrac{1}{n-1}+\dfrac{r}{n(n-1)^{2}}\right)
a_{0}+\dfrac{r}{n(n-1)}a_{1}-(-1)^{\ell -1}(n-1)^{\ell -1}\right) S$ \\
& $+\left( \left( -1-\dfrac{r}{n(n-1)}\right) a_{0}-\dfrac{r}{n}a_{1}\right)
$ \\ \hline
$(\varepsilon )$-Sasakian & $-\,\left( \left( \dfrac{1}{n-1}-\dfrac{%
\varepsilon r}{n(n-1)^{2}}\right) a_{0}-\dfrac{\varepsilon r}{n(n-1)}%
a_{1}-(\varepsilon )^{\ell -1}(n-1)^{\ell -1}\right) S$ \\
& $+\left( \left( \varepsilon -\dfrac{r}{n(n-1)}\right) a_{0}-\dfrac{r}{n}%
a_{1}\right) $ \\ \hline
para-Sasakian & $-\,\left( \left( \dfrac{1}{n-1}+\dfrac{r}{n(n-1)^{2}}%
\right) a_{0}+\dfrac{r}{n(n-1)}a_{1}-(-1)^{\ell -1}(n-1)^{\ell -1}\right) S$
\\
& $+\left( \left( -1-\dfrac{r}{n(n-1)}\right) a_{0}-\dfrac{r}{n}a_{1}\right)
$ \\ \hline
$(\varepsilon )$-para-Sasakian & $-\,\left( \left( \dfrac{1}{n-1}+\dfrac{%
\varepsilon r}{n(n-1)^{2}}\right) a_{0}+\dfrac{\varepsilon r}{n(n-1)}%
a_{1}-(-\varepsilon )^{\ell -1}(n-1)^{\ell -1}\right) S$ \\
& $+\left( \left( -\varepsilon -\dfrac{r}{n(n-1)}\right) a_{0}-\dfrac{r}{n}%
a_{1}\right) $ \\ \hline
\end{tabular}%
\]
\end{cor}

\begin{cor}
Let $M$ be an $n$-dimensional $\left( {\cal P},S,S^{\ell }\right) $%
-pseudosymmetric $\left( N(k),\xi \right) $-semi-Riemannian manifold. Then
we have the following table\/{\rm :}~%
\[
\begin{tabular}{|l|l|}
\hline
${\boldmath M}$ & ${\boldmath LS}^{\ell }{=}$ \\ \hline
$N(k)$-contact metric & $\left( k^{\ell -1}(n-1)^{\ell -1}L-\dfrac{1}{n-1}%
\right) S+kg$ \\ \hline
Sasakian & $\left( (n-1)^{\ell -1}L-\dfrac{1}{n-1}\right) S+g$ \\ \hline
Kenmotsu & $\left( (-1)^{\ell -1}(n-1)^{\ell -1}L-\dfrac{1}{n-1}\right) S-g$
\\ \hline
$(\varepsilon )$-Sasakian & $\left( (\varepsilon )^{\ell -1}(n-1)^{\ell -1}L-%
\dfrac{1}{n-1}\right) S+\varepsilon g$ \\ \hline
para-Sasakian & $\left( (-1)^{\ell -1}(n-1)^{\ell -1}L-\dfrac{1}{n-1}\right)
S-g$ \\ \hline
$(\varepsilon )$-para-Sasakian & $\left( (-\varepsilon )^{\ell
-1}(n-1)^{\ell -1}L-\dfrac{1}{n-1}\right) S-\varepsilon g$ \\ \hline
\end{tabular}%
\]
\end{cor}

\begin{cor}
Let $M$ be an $n$-dimensional $\left( {\cal M},S,S^{\ell }\right) $%
-pseudosymmetric $\left( N(k),\xi \right) $-semi-Riemannian manifold. Then
we have the following table\/{\rm :}~%
\[
\begin{tabular}{|l|l|}
\hline
${\boldmath M}$ & ${\boldmath LS^{\ell }=}$ \\ \hline
$N(k)$-contact metric & $\dfrac{1}{2k(n-1)^{2}}S^{2}+\left( Lk^{\ell
-1}(n-1)^{\ell -1}-\dfrac{1}{n-1}\right) S+\dfrac{k}{2}g$ \\ \hline
Sasakian & $\dfrac{1}{2(n-1)^{2}}S^{2}+\left( L(n-1)^{\ell -1}-\dfrac{1}{n-1}%
\right) S+\dfrac{1}{2}g$ \\ \hline
Kenmotsu & $-\,\dfrac{1}{2(n-1)^{2}}S^{2}+\left( L(-1)^{\ell -1}(n-1)^{\ell
-1}-\dfrac{1}{n-1}\right) S-\dfrac{1}{2}g$ \\ \hline
$(\varepsilon )$-Sasakian & $\dfrac{\varepsilon }{2(n-1)^{2}}S^{2}+\left(
L(\varepsilon )^{\ell -1}(n-1)^{\ell -1}-\dfrac{1}{n-1}\right) S+\dfrac{%
\varepsilon }{2}g$ \\ \hline
para-Sasakian & $-\,\dfrac{1}{2(n-1)^{2}}S^{2}+\left( L(-1)^{\ell
-1}(n-1)^{\ell -1}-\dfrac{1}{n-1}\right) S-\dfrac{1}{2}g$ \\ \hline
$(\varepsilon )$-para-Sasakian & $-\,\dfrac{\varepsilon }{2(n-1)^{2}}%
S^{2}+\left( L(-\varepsilon )^{\ell -1}(n-1)^{\ell -1}-\dfrac{1}{n-1}\right)
S-\dfrac{\varepsilon }{2}g$ \\ \hline
\end{tabular}%
\]
\end{cor}

\begin{cor}
Let $M$ be an $n$-dimensional $\left( {\cal W}_{0},S,S^{\ell }\right) $%
-pseudosymmetric $\left( N(k),\xi \right) $-semi-Riemannian manifold. Then
we have the following table\/{\rm :}~%
\[
\begin{tabular}{|l|l|}
\hline
${\boldmath M}$ & ${\boldmath LS}^{\ell }{=}$ \\ \hline
$N(k)$-contact metric & $\dfrac{1}{k(n-1)^{2}}S^{2}+\left( k^{\ell
-1}(n-1)^{\ell -1}L-\dfrac{2}{n-1}\right) S+kg$ \\ \hline
Sasakian & $\dfrac{1}{(n-1)^{2}}S^{2}+\left( (n-1)^{\ell -1}L-\dfrac{2}{n-1}%
\right) S+g$ \\ \hline
Kenmotsu & $-\,\dfrac{1}{(n-1)^{2}}S^{2}+\left( (-1)^{\ell -1}(n-1)^{\ell
-1}L-\dfrac{2}{n-1}\right) S-g$ \\ \hline
$(\varepsilon )$-Sasakian & $\dfrac{\varepsilon }{(n-1)^{2}}S^{2}+\left(
(\varepsilon )^{\ell -1}(n-1)^{\ell -1}L-\dfrac{2}{n-1}\right) S+\varepsilon
g$ \\ \hline
para-Sasakian & $-\,\dfrac{1}{(n-1)^{2}}S^{2}+\left( (-1)^{\ell
-1}(n-1)^{\ell -1}L-\dfrac{2}{n-1}\right) S-g$ \\ \hline
$(\varepsilon )$-para-Sasakian & $-\,\dfrac{1}{(n-1)^{2}}S^{2}+\left(
(-\varepsilon )^{\ell -1}(n-1)^{\ell -1}L-\dfrac{2}{n-1}\right)
S-\varepsilon g$ \\ \hline
\end{tabular}%
\]
\end{cor}

\begin{cor}
Let $M$ be an $n$-dimensional $\left( {\cal W}_{0}^{\ast },S,S^{\ell
}\right) $-pseudosymmetric $\left( N(k),\xi \right) $-semi-Riemannian
manifold. Then we have the following table\/{\rm :}~%
\[
\begin{tabular}{|l|l|}
\hline
${\boldmath M}$ & ${\boldmath LS}^{\ell }{=}$ \\ \hline
$N(k)$-contact metric & $-\,\dfrac{1}{k(n-1)^{2}}S^{2}+k^{\ell
-1}(n-1)^{\ell -1}LS+kg$ \\ \hline
Sasakian & $-\,\dfrac{1}{(n-1)^{2}}S^{2}+(n-1)^{\ell -1}LS+g$ \\ \hline
Kenmotsu & $\dfrac{1}{(n-1)^{2}}S^{2}+(-1)^{\ell -1}(n-1)^{\ell -1}LS-g$ \\
\hline
$(\varepsilon )$-Sasakian & $-\,\dfrac{\varepsilon }{(n-1)^{2}}%
S^{2}+(\varepsilon )^{\ell -1}(n-1)^{\ell -1}LS+\varepsilon g$ \\ \hline
para-Sasakian & $\dfrac{1}{(n-1)^{2}}S^{2}+(-1)^{\ell -1}(n-1)^{\ell -1}LS-g$
\\ \hline
$(\varepsilon )$-para-Sasakian & $\dfrac{1}{(n-1)^{2}}S^{2}+(-\varepsilon
)^{\ell -1}(n-1)^{\ell -1}LS-\varepsilon g$ \\ \hline
\end{tabular}%
\]
\end{cor}

\begin{cor}
Let $M$ be an $n$-dimensional $\left( {\cal W}_{1},S,S^{\ell }\right) $%
-pseudosymmetric $\left( N(k),\xi \right) $-semi-Riemannian manifold. Then
we have the following table\/{\rm :}~%
\[
\begin{tabular}{|l|l|}
\hline
${\boldmath M}$ & ${\boldmath LS}^{\ell }{=}$ \\ \hline
$N(k)$-contact metric & $\left( k^{\ell -1}(n-1)^{\ell -1}L-\dfrac{1}{n-1}%
\right) S+kg$ \\ \hline
Sasakian & $\left( (n-1)^{\ell -1}L-\dfrac{1}{n-1}\right) S+g$ \\ \hline
Kenmotsu & $\left( (-1)^{\ell -1}(n-1)^{\ell -1}L-\dfrac{1}{n-1}\right) S-g$
\\ \hline
$(\varepsilon )$-Sasakian & $\left( (\varepsilon )^{\ell -1}(n-1)^{\ell -1}L-%
\dfrac{1}{n-1}\right) S+\varepsilon g$ \\ \hline
para-Sasakian & $\left( (-1)^{\ell -1}(n-1)^{\ell -1}L-\dfrac{1}{n-1}\right)
S-g$ \\ \hline
$(\varepsilon )$-para-Sasakian & $\left( (-\varepsilon )^{\ell
-1}(n-1)^{\ell -1}L-\dfrac{1}{n-1}\right) S-\varepsilon g$ \\ \hline
\end{tabular}%
\]
\end{cor}

\begin{cor}
Let $M$ be an $n$-dimensional $\left( {\cal W}_{1}^{\ast },S,S^{\ell
}\right) $-pseudosymmetric $\left( N(k),\xi \right) $-semi-Riemannian
manifold. Then we have the following table\/{\rm :}~%
\[
\begin{tabular}{|l|l|}
\hline
${\boldmath M}$ & ${\boldmath LS}^{\ell }{=}$ \\ \hline
$N(k)$-contact metric & $\left( k^{\ell -1}(n-1)^{\ell -1}L-\dfrac{1}{n-1}%
\right) S+kg$ \\ \hline
Sasakian & $\left( (n-1)^{\ell -1}L-\dfrac{1}{n-1}\right) S+g$ \\ \hline
Kenmotsu & $\left( (-1)^{\ell -1}(n-1)^{\ell -1}L-\dfrac{1}{n-1}\right) S-g$
\\ \hline
$(\varepsilon )$-Sasakian & $\left( (\varepsilon )^{\ell -1}(n-1)^{\ell -1}L-%
\dfrac{1}{n-1}\right) S+\varepsilon g$ \\ \hline
para-Sasakian & $\left( (-1)^{\ell -1}(n-1)^{\ell -1}L-\dfrac{1}{n-1}\right)
S-g$ \\ \hline
$(\varepsilon )$-para-Sasakian & $\left( (-\varepsilon )^{\ell
-1}(n-1)^{\ell -1}L-\dfrac{1}{n-1}\right) S-\varepsilon g$ \\ \hline
\end{tabular}%
\]
\end{cor}

\begin{cor}
Let $M$ be an $n$-dimensional $\left( {\cal W}_{2},S,S^{\ell }\right) $%
-pseudosymmetric $\left( N(k),\xi \right) $-semi-Riemannian manifold. Then
we have the following table\/{\rm :}~%
\[
\begin{tabular}{|l|l|}
\hline
${\boldmath M}$ & ${\boldmath LS}^{\ell }{=}$ \\ \hline
$N(k)$-contact metric & $\dfrac{1}{k(n-1)^{2}}S^{2}+\left( k^{\ell
-1}(n-1)^{\ell -1}L-\dfrac{1}{n-1}\right) S$ \\ \hline
Sasakian & $\dfrac{1}{(n-1)^{2}}S^{2}+\left( (n-1)^{\ell -1}L-\dfrac{1}{n-1}%
\right) S$ \\ \hline
Kenmotsu & $-\,\dfrac{1}{(n-1)^{2}}S^{2}+\left( (-1)^{\ell -1}(n-1)^{\ell
-1}L-\dfrac{1}{n-1}\right) S$ \\ \hline
$(\varepsilon )$-Sasakian & $\dfrac{\varepsilon }{(n-1)^{2}}S^{2}+\left(
(\varepsilon )^{\ell -1}(n-1)^{\ell -1}L-\dfrac{1}{n-1}\right) S$ \\ \hline
para-Sasakian & $-\,\dfrac{1}{(n-1)^{2}}S^{2}+\left( (-1)^{\ell
-1}(n-1)^{\ell -1}L-\dfrac{1}{n-1}\right) S$ \\ \hline
$(\varepsilon )$-para-Sasakian & $-\,\dfrac{1}{(n-1)^{2}}S^{2}+\left(
(-\varepsilon )^{\ell -1}(n-1)^{\ell -1}L-\dfrac{1}{n-1}\right) S$ \\ \hline
\end{tabular}%
\]
\end{cor}

\begin{cor}
Let $M$ be an $n$-dimensional $\left( {\cal W}_{3},S,S^{\ell }\right) $%
-pseudosymmetric $\left( N(k),\xi \right) $-semi-Riemannian manifold. Then
we have the following table\/{\rm :}~%
\[
\begin{tabular}{|l|l|}
\hline
${\boldmath M}$ & ${\boldmath LS}^{\ell }{=}$ \\ \hline
$N(k)$-contact metric & $\left( k^{\ell -1}(n-1)^{\ell -1}L-\dfrac{2}{n-1}%
\right) S+2kg$ \\ \hline
Sasakian & $\left( (n-1)^{\ell -1}L-\dfrac{2}{n-1}\right) S+2g$ \\ \hline
Kenmotsu & $\left( (-1)^{\ell -1}(n-1)^{\ell -1}L-\dfrac{2}{n-1}\right) S-2g$
\\ \hline
$(\varepsilon )$-Sasakian & $\left( (\varepsilon )^{\ell -1}(n-1)^{\ell -1}L-%
\dfrac{2}{n-1}\right) S+2\varepsilon g$ \\ \hline
para-Sasakian & $\left( (-1)^{\ell -1}(n-1)^{\ell -1}L-\dfrac{2}{n-1}\right)
S-2g$ \\ \hline
$(\varepsilon )$-para-Sasakian & $\left( (-\varepsilon )^{\ell
-1}(n-1)^{\ell -1}L-\dfrac{2}{n-1}\right) S-2\varepsilon g$ \\ \hline
\end{tabular}%
\]
\end{cor}

\begin{cor}
Let $M$ be an $n$-dimensional $\left( {\cal W}_{4},S,S^{\ell }\right) $%
-pseudosymmetric $\left( N(k),\xi \right) $-semi-Riemannian manifold. Then
we have the following table\/{\rm :}~
\[
\begin{tabular}{|l|l|}
\hline
${\boldmath M}$ & ${\boldmath LS}^{\ell }{=}$ \\ \hline
$N(k)$-contact metric & $\dfrac{1}{k(n-1)^{2}}S^{2}+\left( k^{\ell
-1}(n-1)^{\ell -1}L-\dfrac{1}{n-1}\right) S+kg-\varepsilon k\eta \otimes
\eta $ \\ \hline
Sasakian & $\dfrac{1}{(n-1)^{2}}S^{2}+\left( (n-1)^{\ell -1}L-\dfrac{1}{n-1}%
\right) S+g-\eta \otimes \eta $ \\ \hline
Kenmotsu & $-\,\dfrac{1}{(n-1)^{2}}S^{2}+\left( (-1)^{\ell -1}(n-1)^{\ell
-1}L-\dfrac{1}{n-1}\right) S-g+\eta \otimes \eta $ \\ \hline
$(\varepsilon )$-Sasakian & $\dfrac{\varepsilon }{(n-1)^{2}}S^{2}+\left(
(\varepsilon )^{\ell -1}(n-1)^{\ell -1}L-\dfrac{1}{n-1}\right) S+\varepsilon
g-\eta \otimes \eta $ \\ \hline
para-Sasakian & $-\,\dfrac{1}{(n-1)^{2}}S^{2}+\left( (-1)^{\ell
-1}(n-1)^{\ell -1}L-\dfrac{1}{n-1}\right) S-g+\eta \otimes \eta $ \\ \hline
$(\varepsilon )$-para-Sasakian & $-\,\dfrac{\varepsilon }{(n-1)^{2}}%
S^{2}+\left( (-\varepsilon )^{\ell -1}(n-1)^{\ell -1}L-\dfrac{1}{n-1}\right)
S-\varepsilon g+\eta \otimes \eta $ \\ \hline
\end{tabular}%
\]
\end{cor}

\begin{cor}
Let $M$ be an $n$-dimensional $\left( {\cal W}_{5},S,S^{\ell }\right) $%
-pseudosymmetric $\left( N(k),\xi \right) $-semi-Riemannian manifold. Then
we have the following table\/{\rm :}~%
\[
\begin{tabular}{|l|l|}
\hline
${\boldmath M}$ & ${\boldmath LS}^{\ell }{=}$ \\ \hline
$N(k)$-contact metric & $\dfrac{1}{k(n-1)^{2}}S^{2}+\left( k^{\ell
-1}(n-1)^{\ell -1}L-\dfrac{2}{n-1}\right) S+kg$ \\ \hline
Sasakian & $\dfrac{1}{(n-1)^{2}}S^{2}+\left( (n-1)^{\ell -1}L-\dfrac{2}{n-1}%
\right) S+g$ \\ \hline
Kenmotsu & $-\,\dfrac{1}{(n-1)^{2}}S^{2}+\left( (-1)^{\ell -1}(n-1)^{\ell
-1}L-\dfrac{2}{n-1}\right) S-g$ \\ \hline
$(\varepsilon )$-Sasakian & $\dfrac{\varepsilon }{(n-1)^{2}}S^{2}+\left(
(\varepsilon )^{\ell -1}(n-1)^{\ell -1}L-\dfrac{2}{n-1}\right) S+\varepsilon
g$ \\ \hline
para-Sasakian & $-\,\dfrac{1}{(n-1)^{2}}S^{2}+\left( (-1)^{\ell
-1}(n-1)^{\ell -1}L-\dfrac{2}{n-1}\right) S-g$ \\ \hline
$(\varepsilon )$-para-Sasakian & $-\,\dfrac{\varepsilon }{(n-1)^{2}}%
S^{2}+\left( (-\varepsilon )^{\ell -1}(n-1)^{\ell -1}L-\dfrac{2}{n-1}\right)
S-\varepsilon g$ \\ \hline
\end{tabular}%
\]
\end{cor}

\begin{cor}
Let $M$ be an $n$-dimensional $\left( {\cal W}_{6},S,S^{\ell }\right) $%
-pseudosymmetric $\left( N(k),\xi \right) $-semi-Riemannian manifold. Then
we have the following table\/{\rm :}~%
\[
\begin{tabular}{|l|l|}
\hline
${\boldmath M}$ & ${\boldmath LS}^{\ell }{=}$ \\ \hline
$N(k)$-contact metric & $\left( k^{\ell -1}(n-1)^{\ell -1}L-\dfrac{2}{n-1}%
\right) S+kg+k\eta \otimes \eta $ \\ \hline
Sasakian & $\left( (n-1)^{\ell -1}L-\dfrac{2}{n-1}\right) S+g+\eta \otimes
\eta $ \\ \hline
Kenmotsu & $\left( (-1)^{\ell -1}(n-1)^{\ell -1}L-\dfrac{2}{n-1}\right)
S-g-\eta \otimes \eta $ \\ \hline
$(\varepsilon )$-Sasakian & $\left( (\varepsilon )^{\ell -1}(n-1)^{\ell -1}L-%
\dfrac{2}{n-1}\right) S+\varepsilon g+\varepsilon \eta \otimes \eta $ \\
\hline
para-Sasakian & $\left( (-1)^{\ell -1}(n-1)^{\ell -1}L-\dfrac{2}{n-1}\right)
S-g-\eta \otimes \eta $ \\ \hline
$(\varepsilon )$-para-Sasakian & $\left( (-\varepsilon )^{\ell
-1}(n-1)^{\ell -1}L-\dfrac{2}{n-1}\right) S-\varepsilon g-\varepsilon \eta
\otimes \eta $ \\ \hline
\end{tabular}%
\]
\end{cor}

\begin{cor}
Let $M$ be an $n$-dimensional $\left( {\cal W}_{7},S,S^{\ell }\right) $%
-pseudosymmetric $\left( N(k),\xi \right) $-semi-Riemannian manifold. Then
we have the following table\/{\rm :}~%
\[
\begin{tabular}{|l|l|}
\hline
${\boldmath M}$ & ${\boldmath LS}^{\ell }{=}$ \\ \hline
$N(k)$-contact metric & $\left( k^{\ell -1}(n-1)^{\ell -1}L-\dfrac{2}{n-1}%
\right) S+2kg$ \\ \hline
Sasakian & $\left( (n-1)^{\ell -1}L-\dfrac{2}{n-1}\right) S+2g$ \\ \hline
Kenmotsu & $\left( (-1)^{\ell -1}(n-1)^{\ell -1}L-\dfrac{2}{n-1}\right) S-2g$
\\ \hline
$(\varepsilon )$-Sasakian & $\left( (\varepsilon )^{\ell -1}(n-1)^{\ell -1}L-%
\dfrac{2}{n-1}\right) S+2\varepsilon g$ \\ \hline
para-Sasakian & $\left( (-1)^{\ell -1}(n-1)^{\ell -1}L-\dfrac{2}{n-1}\right)
S-2g$ \\ \hline
$(\varepsilon )$-para-Sasakian & $\left( (-\varepsilon )^{\ell
-1}(n-1)^{\ell -1}L-\dfrac{2}{n-1}\right) S-2\varepsilon g$ \\ \hline
\end{tabular}%
\]
\end{cor}

\begin{cor}
Let $M$ be an $n$-dimensional $\left( {\cal W}_{8},S,S^{\ell }\right) $%
-pseudosymmetric $\left( N(k),\xi \right) $-semi-Riemannian manifold. Then
we have the following table\/{\rm :}~%
\[
\begin{tabular}{|l|l|}
\hline
${\boldmath M}$ & ${\boldmath LS}^{\ell }=$ \\ \hline
$N(k)$-contact metric & $\left( Lk^{\ell -1}(n-1)^{\ell -1}-\dfrac{2}{n-1}%
\right) S+kg+k\eta \otimes \eta $ \\ \hline
Sasakian & $\left( L(n-1)^{\ell -1}-\dfrac{2}{n-1}\right) S+kg+\eta \otimes
\eta $ \\ \hline
Kenmotsu & $\left( L(-1)^{\ell -1}(n-1)^{\ell -1}-\dfrac{2}{n-1}\right)
S+kg-\eta \otimes \eta $ \\ \hline
$(\varepsilon )$-Sasakian & $\left( L(\varepsilon )^{\ell -1}(n-1)^{\ell -1}-%
\dfrac{2}{n-1}\right) S+kg+\eta \otimes \eta $ \\ \hline
para-Sasakian & $\left( L(-1)^{\ell -1}(n-1)^{\ell -1}-\dfrac{2}{n-1}\right)
S+kg-\eta \otimes \eta $ \\ \hline
$(\varepsilon )$-para-Sasakian & $\left( L(-\varepsilon )^{\ell
-1}(n-1)^{\ell -1}-\dfrac{2}{n-1}\right) S+kg-\eta \otimes \eta $ \\ \hline
\end{tabular}%
\]
\end{cor}

\begin{cor}
Let $M$ be an $n$-dimensional $\left( {\cal W}_{9},S,S^{\ell }\right) $%
-pseudosymmetric $\left( N(k),\xi \right) $-semi-Riemannian manifold. Then
we have the following table\/{\rm :}~%
\[
\begin{tabular}{|l|l|}
\hline
${\boldmath M}$ & ${\boldmath LS}^{\ell }=$ \\ \hline
$N(k)$-contact metric & $\left( k^{\ell -1}(n-1)^{\ell -1}L-\dfrac{1}{n-1}%
\right) S+k\eta \otimes \eta $ \\ \hline
Sasakian & $\left( (n-1)^{\ell -1}L-\dfrac{1}{n-1}\right) S+\eta \otimes
\eta $ \\ \hline
Kenmotsu & $\left( (-1)^{\ell -1}(n-1)^{\ell -1}L-\dfrac{1}{n-1}\right)
S-\eta \otimes \eta $ \\ \hline
$(\varepsilon )$-Sasakian & $\left( (\varepsilon )^{\ell -1}(n-1)^{\ell -1}L-%
\dfrac{1}{n-1}\right) S+\eta \otimes \eta $ \\ \hline
para-Sasakian & $\left( (-1)^{\ell -1}(n-1)^{\ell -1}L-\dfrac{1}{n-1}\right)
S-\eta \otimes \eta $ \\ \hline
$(\varepsilon )$-para-Sasakian & $\left( (-\varepsilon )^{\ell
-1}(n-1)^{\ell -1}L-\dfrac{1}{n-1}\right) S-\eta \otimes \eta $ \\ \hline
\end{tabular}%
\]
\end{cor}

\begin{rem-new}
If in the Theorem \ref{GCT-rssss}, we take $M$ be an $n$-dimensional $\left(
{\cal T}_{\!a},S,S\right) $-pseudosymmetric $\left( N(k),\xi \right) $%
-semi-Riemannian manifold. Then the result is same as given in \cite[Theorem
7.6]{TG}.
\end{rem-new}

\begin{cor}
\label{th-T-ric-pseudo-1} Let $M$ be an $n$-dimensional $(R,S_{{\cal T}%
_{a}},S^{\ell })$-pseudosymmetric $(N(k),\xi )$-semi-Riemannian manifold.
Then
\begin{eqnarray*}
&&(a_{0}+na_{1}+a_{2}+a_{3}+a_{5}+a_{6})\times \\
&&(Lk(n-1)S^{\ell }-Lk^{\ell }(n-1)^{\ell }S+kS-k^{2}(n-1)g) \\
&=&Lr(a_{4}+(n-1)a_{7})(S^{\ell }-k^{\ell }(n-1)^{\ell }g).
\end{eqnarray*}%
\medskip
\end{cor}

\medskip

\noindent Department of Mathematics and DST-CIMS\newline
Faculty of Science\newline
Banaras Hindu University\newline
Varanasi-221005\newline
mmtripathi66@yahoo.com \medskip

\noindent Department of Mathematics\newline
Faculty of Science\newline
Banaras Hindu University\newline
Varanasi-221005\newline
punam\_2101@yahoo.co.in


\begin{thebibliography}{99}
\bibitem{Beem-Ehrlich-81} J.K. Beem and P.E. Ehrlich, {\bf Global Lorentzian
geometry}, Marcel Dekker, New York, 1981.

\bibitem{Blair-76} D.E. Blair, {\bf Contact manifolds in Riemannian geometry}%
, Lectures Notes in Mathematics, Springer-Verlag, Berlin, {\bf 509}~(1976),
146.

\bibitem{Blair-Kim-Tripathi-05} D.E. Blair, J.-S. Kim and M.M. Tripathi,
{\em On the concircular curvature tensor of a contact metric manifold}, J.
Korean Math. Soc. {\bf 42}~(2005), no. 5, 883--892.

\bibitem{Deszcz-89} R. Deszcz, {\em On Ricci-pseudosymmetric warped products}%
, Demonstratio Math. {\bf 22}~(1989), 1053--1065.

\bibitem{Deszcz-90} R. Deszcz, {\em On conformally flat Riemannian manifolds
satisfying certain curvature conditions}, Tensor (N.S.) {\bf 49}~(1990),
134-145.

\bibitem{Deszcz-Grycak-87} R. Deszcz and W. Grycak, {\em On some class of
warped product manifolds}, Bull. Inst. Math. Acad. Sinica {\bf 15}~(1987),
311--322.


\bibitem{DVY} R. Deszcz, L. Verstraelen and S. Yaprak, {\em Pseudosymmetric
hypersurfaces in $4$-dimensional spaces of constant curvature}, Bull. Inst.
Math. Acad. Sinica {\bf 22}~(1994), no. 2, 167--179.

\bibitem{Duggal-90-IJMMS} K.L. Duggal, {\em Space time manifolds and contact
structures}, Internat. J. Math. Math. Sci. {\bf 13}(1990), 545-554.

\bibitem{Eisenhart-49} L.P. Eisenhart, {\bf Riemannian Geometry}, Princeton
University Press, 1949.

\bibitem{Ishii-57} Y. Ishii, {\em On conharmonic transformations}, Tensor
(N.S.) {\bf 7}~(1957), 73--80.

\bibitem{Hong-Ozgur-Tripathi-06} S. Hong, C. \"{O}zg\"{u}r and M.M.
Tripathi, {\em On some classes of Kenmotsu manifold}, Kuwait J. Sci. Engrg.
{\bf 33}~(2006), no. 2, 19--32.

\bibitem{Kenmotsu-72} K. Kenmotsu, {\em A class of almost contact Riemannian
manifold}, T\^{o}hoku Math. J. (2) {\bf 24}~(1972), 93--103.

\bibitem{Matsumoto-89} K. Matsumoto, {\em On Lorentzian paracontact manifolds%
}, Bull. Yamagata Univ. Nat. Sci. {\bf 12}~(1989), no. 2, 151--156.

\bibitem{ONeill-83} B. O'Neill, {\bf Semi-Riemannian geometry with
applications to relativity}, Academic Press, New York, London, 1983.

\bibitem{Ozgur-2005} C. \"{O}zg\"{u}r, {\em On a class of para-Sasakian
manifold}, Turk. J. Math. {\bf 29}~(2005), 249--257.

\bibitem{Ozgur-06} C. \"{O}zg\"{u}r, {\em On Kenmotsu manifolds satisfying
certain pseudosymmetry conditions}, World Applied Sciences Journal {\bf 1~}%
(2006), no. 2, 144--149.

\bibitem{Pokhariyal-Mishra-70} G.P. Pokhariyal and R.S. Mishra, {\em %
Curvature tensors and their relativistic significance}, Yokohama Math. J.
{\bf 18}~(1970), 105--108.

\bibitem{Pokhariyal-Mishra-71} G.P. Pokhariyal and R.S. Mishra, {\em %
Curvature tensors and their relativistic significance~II}, Yokohama Math. J.
{\bf 19}~(1971), no. 2, 97--103.

\bibitem{Pokhariyal-82} G.P. Pokhariyal, {\em Relativistic significance of
curvature tensors}, Internat. J. Math. Math. Sci. {\bf 5}~(1982), no. 1,
133--139.

\bibitem{Prasad-2002} B. Prasad, {\em A pseudo projective curvature tensor
on a Riemannian manifold}, Bull. Calcutta Math. Soc. {\bf 94}~(2002), no. 3,
163--166.

\bibitem{Prvanovic-90} M. Prvanovi\'c, {\em On SP-Sasakian manifold
satisfying some curvature conditions}, SUT J. Math. {\bf 26}~(1990),
201--220.

\bibitem{Sasaki-60} S. Sasaki, {\em On differentiate manifolds with certain
structures which are closely related to almost contact structure~I}, T\^{o}%
hoku Math. J. {\bf 12}~(1960), 459-476.

\bibitem{Sato-76} I. Sat\={o}, {\em On a structure similar to the almost
contact structure}, Tensor (N.S.) {\bf 30}~(1976), no. 3, 219-224.

\bibitem{Takahashi-69} T. Takahashi, {\em Sasakian manifold with
pseudo-Riemannian metric}, T\^{o}hoku Math J. {\bf 21}~(1969), 271--290.

\bibitem{Tanno-88} S. Tanno, {\em Ricci curvatures of contact Riemannian
manifolds}, T\^{o}hoku Math J. {\bf 40}~(1988), 441--448.

\bibitem{TKYK-09} M.M. Tripathi, E. K\i l\i \c{c}, S. Y\"{u}ksel Perkta\c{s}
and S. Kele\c{s}, {\em Indefinite almost paracontact metric manifolds},
Internat. J. Math. Math. Sci. {\bf 2010}~(2010), 19 pp. Art. ID 846195.

\bibitem{Tripathi-Gupta} M.M. Tripathi and P. Gupta, {\em ${\cal T}$%
-curvature tensor on a semi-Riemannian manifolds}, Jour. Adv. Math. Stud.
{\bf 4}~(2011), no. 1, 117--129.

\bibitem{TG} M.M. Tripathi and P. Gupta, {\em On $(N(k),\xi )$%
-semi-Riemannian manifolds: Semisymmetry}, Preprint.

\bibitem{Yano-40} K. Yano, {\em Concircular Geometry~I. Concircular
transformations}, Math. Institute, Tokyo Imperial Univ. Proc. {\bf 16}%
~(1940), 195--200.

\bibitem{Yano-Bochner-53} K. Yano and S. Bochner, {\bf Curvature and Betti
numbers}, Annals of Mathematics Studies {\bf 32}, Princeton University
Press, 1953.

\bibitem{Yano-Sawaki-68} K. Yano and S. Sawaki, {\em Riemannian manifolds
admitting a conformal transformation group}, J. Diff. Geom. {\bf 2}~(1968),
161--184.
\end{thebibliography}
\end{document}